\documentclass[11pt]{article}
\usepackage{amsmath,amsbsy,amsfonts,amssymb}

\usepackage[frenchb]{babel}

 \newtheorem{thm}{Th\'eor\`eme}[section]
 \newtheorem{cor}[thm]{Corollaire}
 \newtheorem{lem}[thm]{Lemme}
 \newtheorem{prop}[thm]{Proposition}
 \newtheorem{rem}[thm]{Remarque}

 \title{Comparaison des param\`etres de Langlands et des exposants \`a l'int\'erieur d'un paquet d'Arthur}
\author{C. M{\oe}glin\\
Institut de math\'ematiques de Jussieu, CNRS\\
moeglin@math.jussieu.fr}
\date{}
\begin{document}

\maketitle
\begin{abstract} In this paper, one proves an idea expressed by Clozel: inside an Arthur's packet, one has the representations in the Langlands' packet inside the Arthur's packet and more tempered representations than these representations. 
\end{abstract}
Dans cet article $F$ est un corps local non archim\'edien et $W_{F}$ est son groupe de Weil.

Dans \cite{clozel} (derni\`eres lignes apr\`es 2.4) qui est bas\'e sur \cite{clozelparkcity}, Clozel \'enonce la philosophie des paquets d'Arthur \`a savoir qu'un paquet d'Arthur doit \^etre compos\'e de repr\'esentations dans le paquet de Langlands naturellement associ\'e au param\`etre du paquet d'Arthur et de repr\'esentations ''plus temp\'er\'ees''; c'est aussi ce qui est expliqu\'e dans le cas des places archim\'ediennes par Adams, Barbasch et Vogan mais nous nous limitons ici aux places finies et nos m\'ethodes ne s'\'etendent pas aux places archim\'ediennes. On peut donner plusieurs sens \`a une telle assertion. Clozel en donne un qui ne s'applique qu'\`a certaines repr\'esentations mais qui est particuli\`erement simple \`a expliquer: on rappelle que les paquets d'Arthur aux places p-adiques sont associ\'es \`a des morphismes de $W_{F}\times SL(2,{\mathbb C})\times SL(2,{\mathbb C})$ dans le $L$-groupe ayant de bonnes propri\'et\'es. On ne consid\`ere que des groupes, $G$, pourlequels on peut voir de tels morphismes comme des repr\'esentations de $W_{F}\times SL(2,{\mathbb C})\times SL(2,{\mathbb C})$ \`a valeurs dans un groupe $GL(m^*_{G},{\mathbb C})$ (ce qui d\'efinit $m^*_{G}$) se factorisant par un sous-groupe de similitudes orthogonales, symplectiques ou unitaires; cela revient \`a dire que $G$ est un groupe classique tel que les repr\'esentations automorphes irr\'eductibles de carr\'e int\'egrable de sa forme int\'erieure quasid\'eploy\'ee sont accessibles par la th\'eorie de l'endoscopie tordue \`a la Arthur; on reprend le yoga expliqu\'e par Arthur auquel on n'ajoute rien (cf. par exemple \cite{arthurseriediscrete}, \cite{arthurgsp4}). Tous les groupes $G$ consid\'er\'es sont en particulier tels que les sous-groupes paraboliques ont des Levi isomorphes \`a un produit de groupes lin\'eaires par un groupe de m\^eme type que $G$ et ceci est aussi vrai pour le groupe dual. Pour \'eviter une multiplication de notations, le plus simple est de consid\'erer que $G$ est un groupe classique usuel.

La propri\'et\'e de base (\cite{arthur}) des paquets de repr\'esentations associ\'ees \`a un tel morphisme $\psi$ est que ce paquet contient toutes les repr\'esentations du paquet de Langlands qu'Arthur a naturellement associ\'e \`a $\psi$ en composant, l'inclusion de $W_{F}\times SL(2,{\mathbb C})$ dans $W_{F}\times SL(2,{\mathbb C})\times SL(2,{\mathbb C})$ avec $\psi$ o\`u la premi\`ere inclusion est  $$(w,g)\mapsto (w, g,(\begin{matrix} \vert w\vert^{1/2}&0\\ 0&\vert w\vert^{-1/2}\end{matrix})).$$
Le fait que le paquet de Langlands est effectivement inclus dans le paquet d'Arthur est certainement d\'emontr\'e par Arthur et se trouve aussi dans \cite{pourshahidi} section {\bf 6}.

 Dans tout l'article on note $\Pi(\psi)$ le paquet de repr\'esentations associ\'ees \`a $\psi$ (on en rappelle la d\'efinition et la construction en \ref{rappelconstruction} ci-dessous) et on note $\Pi(\psi_{L})$ le paquet de Langlands \`a l'int\'erieur du paquet d'Arthur.
.

La remarque de Clozel la plus simple \`a exprimer et que si $\pi\in \Pi(\psi)$ et s'il existe un autre morphisme $\psi'$ tel que $\pi\in \Pi(\psi'_{L})$, alors l'orbite unipotente de $GL(m_{G}^*,{\mathbb C})$ d\'efinie par la restriction de $\psi$ \`a la deuxi\`eme copie de $SL(2,{\mathbb C})$ contient dans sa fermeture l'orbite unipotente d\'efinie par la restriction de $\psi'$ \`a la deuxi\`eme copie de $SL(2,{\mathbb C})$. Nous  d\'emontrons cette remarque en \ref{preuveassertionclozel} ci-dessous.  Bien s\^ur ici, il est important que $\pi$ soit dans le paquet de Langlands associ\'e \`a $\psi'$. Et ce r\'esultat ne s'applique donc qu'a peu d'\'el\'ements de $\Pi(\psi)$ puisqu'en g\'en\'eral si $\pi\in \Pi(\psi)-\Pi(\psi_{L})$, il n'existe pas de morphisme $\psi'$ tel que $\pi\in \Pi(\psi'_{L})$.

On donne donc un sens \`a la remarque de Clozel qui s'applique \`a toute repr\'esentation dans $\Pi(\psi)$. Puisque gr\^ace aux r\'esultats d'Arthur, on a  maintenant, pour les groupes auxquels les r\'esultats d'Arthur s'appliquent, la classification de Langlands des s\'eries discr\`etes (cf \cite{paquetdiscret} et \cite{pourshahidi}), pour toute repr\'esentation irr\'eductible $\pi$, on peut d\'efinir \`a l'aide du groupe dual (au moins th\'eoriquement)  les param\`etres de Langlands de $\pi$ et on peut donc associer \`a $\pi$ un morphisme de $W_{F}\times SL(2,{\mathbb C})$ dans le $L$-groupe de $G$; on voit encore ce morphisme comme une repr\'esentation, $\phi_{\pi}$ de $W_{F}\times SL(2,{\mathbb C})$. On associe alors \`a $\pi$ une orbite unipotente en restreignant $\phi_{\pi}$ \`a $SL(2,{\mathbb C})$; on note $O^L_{\pi}$ cette orbite. On suppose maintenant que $\pi\in \Pi(\psi)$ (avec $\psi$ comme ci-dessus). On montre alors (cf. \ref{comparaisonorbite}) que $O^L_{\pi}$ contient dans sa fermeture l'orbite unipotente de $GL(m_{G}^*,{\mathbb C})$,  d\'efinie par la restriction de $\psi$ \`a la premi\`ere copie de $SL(2,{\mathbb C})$; on v\'erifie ais\'ement que si $\pi\in \Pi(\psi_{L})$, les 2 orbites co\"{\i}ncident. Ce r\'esultat r\'epond \`a la m\^eme philosophie que le pr\'ec\'edent mais diff\`ere par le fait qu'il fait intervenir la premi\`ere copie de $SL(2,{\mathbb C})$ et non la 2e. Il est donc de nature purement locale. Comme application imm\'ediate de ce r\'esultat (sugg\'er\'ee par \cite{clozelimrn}) on obtient que si $\psi$ est  un morphisme d\'eterminant un paquet de repr\'esentations, ce paquet de repr\'esentations contient une repr\'esentation non ramifi\'ee  seulement si la restriction de $\psi$ \`a la premi\`ere copie de $SL(2,{\mathbb C})$ est triviale. De plus, facilement, on v\'erifie que la restriction de $\psi$ \`a $W_{F}$ est alors non ramifi\'ee et  que quand ces 2 conditions sont satisfaites,  le paquet contient une unique repr\'esentation non ramifi\'ee (cf.  \ref{nonramifie}); pour les groupes orthogonaux pairs il faut bien consid\'erer tout le groupe orthogonal.

Cette relation n'est donc sans doute pas celle que Clozel avait en vue. Ce qui est plus int\'eressant d'un point de vue global, est de comparer les exposants des repr\'esentations, car les exposants locaux donnent quelques informations sur les exposants globaux. La propri\'et\'e ici est une majoration des exposants locaux \`a l'aide de l'orbite unipotente associ\'ee \`a la restriction de $\psi$ \`a la 2e copie de $SL(2,{\mathbb C})$ pour tout $\pi\in \Pi(\psi)$. Plus pr\'ecis\'ement, soit $\pi$ une repr\'esentation irr\'eductible, gr\^ace au th\'eor\`eme du quotient de Langlands, et au choix particulier des groupes que l'on a fait ici, il existe une inclusion
$$
\pi\hookrightarrow \times_{(\rho,a,x)\in \mathcal{L}_{\pi}}St(\rho,a)\vert\,\vert^{-x}\times \pi_{temp},
$$
o\`u $\mathcal{L}_{\pi}$ est un ensemble de triplets form\'e chacun d'une repr\'esentation cuspidale unitaire, $\rho$, d'un entier $a$ et d'un r\'eel strictement positif $x$ et o\`u $\pi_{temp}$ est une repr\'esentation temp\'er\'ee irr\'eductible d'un groupe de m\^eme type que $G$ mais de rang en g\'en\'eral plus petit (cf. \ref{definitionorbitelanglands}). L'ensemble $\{x; (\rho,a,x)\in \mathcal{L}_{\pi}\}$ vu comme ensemble avec multiplicit\'e de r\'eels strictement positifs est bien d\'efini et on l'appelle ensemble des exposants de $\pi$; on le note $Exp(\pi)$. Il y a un d\'efaut \`a cette d\'efinition, c'est qu'elle ne tient pas compte de la s\'erie discr\`ete qui porte l'exposant, donc ci-dessus dans le triplet $(\rho,a,x)$, on peut tenir compte du $\rho$ (cf. ci-dessous) mais on ne tient pas compte du $a$. Si $\pi$ est dans le paquet de Langlands associ\'e \`a $\psi$, on conna\^{\i}t parfaitement $Exp(\pi)$: pour l'obtenir, il faut d\'ecomposer $\psi$ en sous-repr\'esentations irr\'eductibles; une repr\'esentation irr\'eductible de $W_{F}\times SL(2,{\mathbb C})\times SL(2,{\mathbb C})$ est la donn\'ee d'un triplet $(\rho,a,b)$ o\`u $\rho$ est une repr\'esentation irr\'eductible de $W_{F}$ et $a,b$ sont des entiers qui d\'eterminent chacun une repr\'esentation irr\'eductible de $SL(2,{\mathbb C})$. Pour les paquets d'Arthur, on ne consid\`ere que le cas o\`u les repr\'esentations $\rho$ intervenant sont unitaires (pour palier \`a l'absence de conjecture de Ramanujan il faudrait accepter des torsions par des caract\`eres de la forme $\vert\,\vert^{x}$ avec $x\in ]-1/2,1/2[$; c'est une variante facile que l'on ne traite pas ici). On note $Jord(\psi)$ l'ensemble avec multiplicit\'e des sous-repr\'esentations irr\'eductibles de $\psi$. On appelle $Exp(\psi)$ la collection de demi-entiers, $\cup_{(\rho,a,b)\in Jord(\psi)}(b-1)/2, \cdots, \delta_{b},$ o\`u pour tout entier $b$, $\delta_{b}$ vaut $1/2$ si $b$ est pair et $1$ si $b$ est impair. On v\'erifie ais\'ement (cf. \ref{description}) que $Exp(\psi)=Exp(\pi)$ si $\pi\in \Pi(\psi_{L})$. On suppose simplement que $\pi\in \Pi(\psi)$ et on montre en \ref{comparaisonexposants}  que $Exp(\pi)\leq Exp(\psi)$ au sens des inclusions d'orbites unipotentes, c'est-\`a-dire que pour tout entier $t$ la somme des $t$ plus grands \'el\'ements de $Exp(\pi)$ est inf\'erieure ou \'egal \`a la somme des $t$ plus grands \'el\'ements de $Exp(\psi)$, on ajoute des 0 si n\'ecessaire. L'assertion est moins jolie que l'assertion sur les orbites unipotentes et elle ne la g\'en\'eralise pas compl\`etement par la remarque d\'ej\`a faite que dans cette assertion on oublie compl\`etement la taille de la repr\'esentation de Steinberg qui ''porte'' l'exposant.

 L'article d\'emontre successivement les 3 propri\'et\'es expliqu\'ees dans cette introduction; la premi\`ere propri\'et\'e est d'une nature diff\'erente des autres. On a essay\'e d'expliquer cela dans l'article: elle
est  cons\'equence d'une propri\'et\'e des modules de Jacquet des repr\'esentations des groupes lin\'eaires associ\'ees aux repr\'esentations de $W_{F}\times SL(2,{\mathbb C})\times SL(2,{\mathbb C})$ consid\'er\'ees; c'est la propri\'et\'e de \ref{lemmesurjacquet}. Les repr\'esentations dans les paquets de Langlands \`a l'int\'erieur d'un paquet d'Arthur sont extr\^emement particuli\`eres (cf. par exemple \ref{description} ci-dessous): elles h\'eritent du fait que les composantes locales des formes automorphes de carr\'e int\'egrable pour les groupes lin\'eaires sont particuli\`erement simples, ces composantes locales sont des induites n\'ecessairement irr\'eductibles de modules de Speh. Pour les groupes consid\'er\'es ici, la situation est un peu plus compliqu\'ee, ces repr\'esentations ne sont pas des induites irr\'eductibles mais peu s'en faut quand m\^eme. En particulier leur module de Jacquet ont des propri\'et\'es assez remarquables et c'est cela que l'on utilise. Le fait qu'une non nullit\'e de modules de Jacquet pour une repr\'esentation dans un paquet d'Arthur entra\^{\i}ne aussi une non nullit\'e de module de Jacquet pour la repr\'esentation du groupe lin\'eaire associ\'e est imm\'ediate si tous les signes dans le paquet d'Arthur sont les m\^emes. Quand ceci n'est pas vrai, la situation est plus compliqu\'ee et la preuve n\'ecessite alors d'entrer dans la structure fine des constructions pour montrer cela (cf. \ref{proprietesimpledejac} ci-dessous). Quand on a cette non nullit\'e des modules de Jacquet, la remarque de Clozel est une propri\'et\'e combinatoire simple qui se r\'esume \`a: soit $S$ et $S'$ 2 segments d\'ecroissants de m\^eme d\'ebut et tels que $S'$ est inclus dans $S$ alors le milieu de $S$ est inf\'erieur ou \'egal au milieu de $S'$ (cf. \ref{preuveassertionclozel}).

Les 2 autres propri\'et\'es ont une d\'emonstration analogue; les assertions sont vraies si le paquet ne contient que des repr\'esentations temp\'er\'ees et la preuve consiste donc \`a se ramener \`a ce cas. Pour cela il faut utiliser des descriptions par r\'ecurrence des \'el\'ements de $\Pi(\psi)$ qui sont donn\'ees en \ref{blocpositif} et \ref{blocnegatif} et la d\'emonstration est purement technique.

Dans l'article on d\'emontre une version un peu plus pr\'ecise des propri\'et\'es annonc\'ees ici, on fixe la repr\'esentation irr\'eductible de $W_{F}$ que l'on identifie via la correspondance de Langlands \`a une repr\'esentation cuspidale d'un groupe lin\'eaire convenable et toutes les propri\'et\'es \'enonc\'ees se font, c\^ot\'e $\psi$, dans la composante isotypique pour l'action de $W_{F}$ et pour les exposants, on d\'emontre la propri\'et\'e en fixant le support cuspidal des repr\'esentations de Steinberg tordues comme \'etant une collection de repr\'esentation de la forme $\rho\vert\,\vert^{y}$ pour $y\in {\mathbb R}$ avec une repr\'esentation cuspidale $\rho$ fix\'ee d'un groupe lin\'eaire convenable. Cela ne co\^ute pas plus cher et simplifie m\^eme les notations.
\pagebreak
\tableofcontents
 
 \pagebreak
 
 \bf Conventions importantes:\rm
 
 on suppose que la classification de Langlands des s\'eries discr\`etes est connues pour les groupes consid\'er\'es, ou plus exactement notre article ne s'applique qu'aux groupes classiques (ou m\^eme de similitudes) pour lesquels cette classification est connue en les termes de \cite{arthurseriediscrete}. Pour le moment le bilan est plut\^ot maigre.

Pour la clart\'e des notations, on a pr\'ef\'er\'e oubli\'e dans le texte le fait que  dans le cas des groupes unitaires,  les groupes lin\'eaires qui interviennent ne sont pas sur le corps de base $F$ mais sur l'extension quadratique  de $F$  qui sert \`a d\'efinir la forme lin\'eaire et au lieu de consid\'erer la contragr\'ediente des repr\'esentations de ces groupes lin\'eaires il faut consid\'erer le dual hermitien.

\section{Description des repr\'esentations de Langlands \`a l'int\'erieur d'un paquet d'Arthur}
\subsection{Description\label{description}}
On fixe $\psi$ comme dans l'introduction, d'o\`u $Jord(\psi)$ d\'efini aussi dans l'introduction et qui donne la d\'ecomposition de $\psi$ en repr\'esentations irr\'eductibles. Dans l'introduction on a d\'efini le morphisme de Langlands associ\'e \`a $\psi$ par la m\'ethode d'Arthur; on le note  $\phi_{\psi}$ et sa d\'ecomposition en sous-repr\'esentations irr\'eductibles est explictement donn\'ee en fonction de celle de $\psi$ par:
$$
\phi_{\psi}=\oplus_{(\rho,a,b)\in Jord(\psi)}\oplus_{c\in [-(b-1)/2,(b-1)/2]}\rho\vert\,\vert^{c}\otimes rep_a$$
comme repr\'esentation de $W_{F}\otimes SL(2,{\mathbb C})$, o\`u pout $a\in {\mathbb N}$, $rep_a$ est la repr\'esentation irr\'eductible de dimension $a$ de $SL(2,{\mathbb C})$.  Les repr\'esentations dans le paquet de Langlands associ\'es \`a $\phi_{\psi}$ sont en bijection avec les repr\'esentations temp\'er\'ees dans le paquet associ\'e \`a
$$\psi_{temp}:=
\oplus_{(\rho,a,b)\in Jord(\psi); b\equiv 1[2]}\rho\otimes rep_a.
$$
A tout $\pi_{temp} \in \Pi(\psi_{temp})$, on associe la repr\'esentation qui est le sous-quotient de Langlands de l'induite:
$$
\times_{(\rho,a,b)\in Jord(\psi); b>1}\times_{c\in [\delta_{b},(b-1)/2]}St(\rho,a)\vert\,\vert^{-c}\times \pi_{temp},\eqno(1)
$$o\`u $\delta_{b}=1/2$ si $b$ est pair et $1$ si $b$ est impair et o\`u l'intervalle est form\'e de demi-entier non entier ou d'entier suivant que $b$ est pair ou impair.
Ceci veut dire que l'on conjugue les premiers termes de l'induite ci-dessus pour que les exposants soient dans la chambre de Weyl n\'egative et que  l'on prend l'unique sous-module irr\'eductible de l'induite ainsi obtenue.

On peut construire ce sous-module irr\'eductible de fa\c{c}on beaucoup plus agr\'eable avec la d\'efinition suivante: soit $\rho$ comme ci-dessus et $[x,y]$ un segment croissant, on note $J(St(\rho,a),x,y)$ l'unique sous-module irr\'eductible pour le $GL$ convenable inclus dans l'induite
$$
\times_{c\in [x,y]}St(\rho,a)\vert\,\vert^{c}.
$$
On a montr\'e en \cite{mw} 1.6.3 que pour $\rho,x,y,\rho',x',y'$, l'induite pour le GL convenable $J(St(\rho,a),x,y)\times J(St(\rho',a'),x',y')$ est irr\'eductible par exemple si $\vert y-y'\vert \leq 1/2$. 

\begin{prop}Les \'el\'ements de $\Pi(\psi_{L})$ sont en bijection avec les \'el\'ements de $\Pi(\psi_{temp})$ l'application $\pi\mapsto \pi_{temp}$ est caract\'eris\'ee par le fait que
$$
\pi \hookrightarrow \times_{(\rho,a,b)\in Jord(\psi); b>1}J(St(\rho,a),-(b-1)/2,-\delta_{b})\times \pi_{temp}.\eqno(2)
$$
De plus dans (2) on peut bouger les premiers facteurs comme on veut.
\end{prop}
C'est une cons\'equence facile (montr\'ee en \cite{seriedeisenstein} 3.5.1) de la propri\'et\'e d'irr\'eductibilit\'e d\'emontr\'ee en \cite{mw} et rappel\'ee avant l'\'enonc\'e. On ne refait pas la d\'emonstration ici car les id\'ees sont celles que l'on va utiliser (et \'ecrire) dna s\ref{jacquetpourlanglands} ci-dessous.

\begin{rem}Cette forme tr\`es particuli\`ere des param\`etres de Langlands des repr\'esentations dans un paquet du type $\Pi(\psi_{L})$ prouve que pour $\pi$ g\'en\'eral, il n'existe pas $\psi$ tel que $\pi\in \Pi(\psi_{L})$. On donnera ci-dessous la d\'efinition des \'el\'ements de $\Pi(\psi)$ mais pour un \'el\'ement de $\Pi(\psi)$ il n'existe pas en g\'en\'eral de morphisme $\psi'$ tel que cet \'el\'ement soit aussi dans $\Pi(\psi'_{L})$ sauf \'evidemment pour les \'el\'ements de $\Pi(\psi_{L})$.
\end{rem}

\subsection{Notations concernant les modules de Jacquet et les repr\'esentations associ\'ees \`a des segments\label{notationjac}}
Soit $\pi$ une repr\'esentation de longueur finie de $G$ et soit $\rho$ une repr\'esentation cuspidale irr\'eductible (qui sera toujours unitaire ici) d'un groupe lin\'eaire $GL(d_{\rho},F)$. Pour $x\in {\mathbb R}$, on d\'efinit $Jac_{x}\pi$; pour cela, il faut que $G$ soit de rang d\'eploy\'e sup\'erieur ou \'egal \`a $d_{\rho}$ (sinon on pose $Jac_{x}\pi=0$) et  cela suppose aussi que l'on a fix\'e un sous-groupe parabolique minimal une fois pour toute de fa\c{c}on \`a avoir la notion de sous-groupe parabolique standard; on fait un tel choix dans tout l'article mais on le fait de fa\c{c}on compatible pour tous les groupes intervenant par exemple quand il est  fait pour $G$ il est fait pour le groupe $G'$ qui arrive imm\'ediatement dans la phrase suivante; pour \^etre plus correct, on fixe un ''tr\`es gros'' groupe $H$ de m\^eme type que $G$ mais de rang bien plus grand; on fixe un parabolique minimal standard pour $H$ et on consid\`ere chaque groupe de m\^eme type que $G$ qui intervient dans cet article comme sous-groupe de $H$, en consid\'erant $GL(d_{H,G})\times G$ comme un sous-groupe de Levi standard dans le parabolique standard convenable pour $H$, o\`u $d_{H,G}$ est la diff\'erence des rangs entre $H$ et $G$.

On consid\`ere alors le sous-groupe parabolique standard de Levi isomorphe \`a $GL(d_{\rho},F)\times G'$, o\`u $G'$ est un groupe de m\^eme type que $G$ mais de rang $d_{\rho}$ de moins et $Jac_{x}\pi$ est par d\'efinition telle que dans la groupe de Grothendieck de $GL(d_{\rho},F)\times G'$, le module de Jacquet de $\pi$ pour le radical unipotent du parabolique d\'ecrit est de la forme $$\rho\vert\,\vert^{x}\times Jac_{x}\pi\oplus_{\sigma,\pi_{\sigma}}\sigma\otimes \pi_{\sigma},$$ o\`u $\sigma$ parcourt l'ensemble des repr\'esentations irr\'eductibles de $GL(d_{\rho},F)$ non isomorphes \`a $\rho\vert\,\vert^x$.

On a besoin d'une notation analogue pour les repr\'esentations des groupes lin\'eaires eux-m\^emes. Soit $\pi^{GL}$ une repr\'esentation  d'un groupe lin\'eaire $GL(m,F)$ ; ici on d\'efinit $Jac^g_{x}\pi^{GL}$, nul si $m<d_{\rho}$ et sinon  de telle sorte que le module de Jacquet de $\pi^{GL}$ pour le sous-groupe parabolique standard usuel de $GL(m,F)$ de Levi $GL(d_{\rho})\times GL(m-d_{\rho},F)$ soit dans le groupe de Grothendieck convenable, de la forme $$\rho\vert\,\vert^{x}\times Jac^{g}_{x}(\pi^{GL})\oplus \sigma\otimes \pi'_{\sigma},$$
o\`u $\sigma$ parcourt l'ensemble des repr\'esentations irr\'eductibles de $GL(d_{\rho},F)$ qui ne sont pas isomorphe \`a $\rho\vert\,\vert^{x}$. On d\'efinit $Jac^d_{-x}\pi^{GL}$ de fa\c{c}on sym\'etrique en utilisant la droite au lieu de la gauche, on a chang\'e $x$ en $-x$ et en utilisant $\rho^*$ au lieu de $\rho$. On pose alors
$$
Jac^\theta_{x}\pi^{GL}=Jac_{-x}^dJac^g_{x}\pi^{GL}
$$
ce qui vaut 0 si $m<2d_{\rho}$.

\

Ces applications peuvent s'it\'erer en particulier on d\'efinit $Jac_{x, \cdots, y}:= Jac_{y}\circ \cdots Jac_{x}$. La propri\'et\'e cl\'e qui vient de la r\'eciprocit\'e de Frobenius est que $$Jac_{x}\pi\neq 0\Rightarrow \exists \sigma, \, \pi\hookrightarrow \rho\vert\,\vert^{x}\times \sigma.$$

De plus $Jac_{x}\circ Jac_{y}= Jac_{y}\circ Jac_{x}$ si $\vert x-y\vert >1$.

\

On fixe toujours $\rho$ et on fixe aussi un segment $[x,y]$, c'est-\`a-dire que $x,y\in {\mathbb R}$ et $x-y\in {\mathbb Z}$ (ici dans l'article, le segment sera form\'e de demi-entiers et sera d\'ecroissant); on note $<x, \cdots, y>_{\rho}$ l'unique sous-repr\'esentation irr\'eductible du groupe lin\'eaire $GL(d_{\rho}(\vert x-y\vert+1),F)$ incluse dans l'induite $\times_{z\in [x,y]}\rho\vert\,\vert^{z}$.

\subsection{Propri\'et\'es des modules de Jacquet des \'el\'ements de $\Pi(\psi_{L})$\label{jacquetpourlanglands}}
On fixe $\rho$.
Soit $\mathcal{I}$ un sous-ensemble de $Jord(\psi)$ form\'e de triplets de la forme $(\rho,a,b)$ le m\^eme $\rho$, v\'erifiant tous $b>1$; on note $\mathcal{E}$ l'union des segments d\'ecroissants $[(a-b)/2, -(a+b)/2+1]$ pour $(\rho,a,b)$ parcourant $\mathcal{I}$ et on r\'eordonne $\mathcal{E}$ de sorte que l'ordre sur $\mathcal{E}$ soit l'ordre d\'ecroissant pour les demi-entiers.
\begin{lem}
Soit $\pi\in \Pi(\psi_{L})$; alors $Jac_{x\in \mathcal{E}}\pi\neq 0$.
\end{lem}
On ordonne $\mathcal{I}$ de sorte que si $(\rho,a,b)<(\rho,a',b')$ alors $(a-b)/2\geq (a'-b')/2$ et on r\'ecrit \ref{description}(2), en utilisant le fait que $\delta_{b}$ est bien d\'efini pour tout $(\rho,a,b)\in \mathcal{I}$ puisque $b>1$ par hypoth\`ese:
$$
\pi\hookrightarrow \times_{(\rho,a,b)\in \mathcal{I}}J(St(\rho,a),-(b-1)/2,-\delta_{b})\times_{(\rho',a',b')\in Jord(\psi)-\mathcal{I}}J(St(\rho',a'),-(b'-1)/2,-\delta_{b'})\times \pi_{temp}.
$$
Il est plus simple de remplacer $\mathcal{I}$ par un segment $[1,v]$ pour $v$ un entier convenalble et d'\'ecrire les \'el\'ements de $\mathcal{I}$ sous la forme $(\rho,a_{j},b_{j}); j\in [1,v]$. On a la suite de morphismes: $\pi\hookrightarrow$
$$
 J(St(\rho,a_{1}),-(b_{1}-1)/2,-\delta_{b_{1}})\times_{j\in [2,v]}J(St(\rho,a_{j}),-(b_{j}-1)/2,-\delta_{b_{j}})$$
 $$\times _{(\rho',a',b')\in Jord(\psi)-\mathcal{I}}J(St(\rho',a'),-(b'-1)/2,-\delta_{b'})\times \pi_{temp}$$
$$
\hookrightarrow
St(\rho,a_{1})\vert\,\vert^{-(b_{1}-1)/2}\times J(St(\rho,a_{1}),-(b_{1}-3)/2,-\delta_{b_{1}})\times 
 _{j\in [2,v]}J(St(\rho,a_{j}),-(b_{j}-1)/2,-\delta_{b_{j}})$$
 $$\times _{(\rho',a',b')\in Jord(\psi)-\mathcal{I}}J(St(\rho',a'),-(b'-1)/2,-\delta_{b'})\times \pi_{temp}$$
 $$
 \simeq 
 St(\rho,a_{1})\vert\,\vert^{-(b_{1}-1)/2}\times 
 _{j\in [2,v]}J(St(\rho,a_{j}),-(b_{j}-1)/2,-\delta_{b_{j}})\times J(St(\rho,a_{1}),-(b_{1}-3)/2,-\delta_{b_{1}})$$
 $$\times _{(\rho',a',b')\in Jord(\psi)-\mathcal{I}}J(St(\rho',a'),-(b'-1)/2,-\delta_{b'})\times \pi_{temp}$$
 $$
 \hookrightarrow 
 St(\rho,a_{1})\vert\,\vert^{-(b_{1}-1)/2}\times St(\rho,a_{2})\vert\,\vert^{-(b_{2}-1)/2}\times  J(St(\rho,a_{2}),-(b_{2}-3)/2,-\delta_{b_{2}})\times$$
 $$ 
 _{j\in [3,v]}J(St(\rho,a_{j}),-(b_{j}-1)/2,-\delta_{b_{j}})\times J(St(\rho,a_{1}),-(b_{1}-3)/2,-\delta_{b_{1}})$$
 $$\times _{(\rho',a',b')\in Jord(\psi)-\mathcal{I}}J(St(\rho',a'),-(b'-1)/2,-\delta_{b'})\times \pi_{temp}
 $$
 et de proche en proche on montre l'existence d'une repr\'esentation $\sigma$ et d'une inclusion
 $$
 \pi\hookrightarrow \times_{j\in [1,v]}St(\rho,a_{j})\vert\,\vert^{-(b_{j}-1)/2}\times \sigma.
 $$
 On va v\'erifier que
 $$
 \times_{j\in [1,v]}St(\rho,a_{j})\vert\,\vert^{-(b_{j}-1)/2}\hookrightarrow \times_{x\in \mathcal{E}}\rho\vert\,\vert^{x}.
 $$
 Ceci utilise la propri\'et\'e de l'ordre mis sur $\mathcal{I}$; d'abord on simplifie l'\'ecriture. Soit  pour $i\in [1,v]$ un ensemble de segments d\'ecroissant $[d_{i},f_{i}]$ rang\'e tels que pour $i<i'$, $d_{i}\geq d_{i'}$. On note $\mathcal{E}$ l'ensemble avec multiplicit\'e des r\'eels $\cup_{i}\{x\in [d_{i},f_{i}]\}$ rang\'e dans l'ordre d\'ecroissant et on veut montrer que 
 $$
 \times_{i\in [1,v]}<d_{i},f_{i}>_{\rho}\hookrightarrow \times_{x\in \mathcal{E}}\rho\vert\,\vert^x.
 $$
 Disons que l'on fait une r\'ecurrence sur $\sum_{i}(d_{i}-f_{i})$. Si ce nombre est nul chaque segment est r\'eduit \`a un terme et l'assertion est cons\'equence du fait que ${E}=\cup_{i\in [1,v]}d_{i}$ dans cet ordre. Soit $i_{0}$ le plus petit indice tel que $d_{i_{0}}>f_{i_{0}}$ et notons $i_{1}$ le plus grand indice tel que $d_{i_{1}}=d_{i_{0}}$. On pose $d=d_{i_{0}}$ a:
 $$
 \times_{i\in [i_{0},i_{1}]}<d_{i},f_{i}>_{\rho}\hookrightarrow \times_{i\in [i_{0},i_{1}[}<d,f_{i}>_{\rho}\times \rho\vert\,\vert^{d} \times <d-1,f_{i_{1}}>_{\rho}$$
 $$
 \simeq
 \rho\vert\,\vert^{d}\times _{i\in [i_{0},i_{1}[}<d,f_{i}>_{\rho}\times <d-1,f_{i_{1}}>_{\rho}
 $$
 et de proche en proche, on montre
 $$
 \times_{i\in [i_{0},i_{1}]}<d_{i},f_{i}>_{\rho}\hookrightarrow \rho\vert\,\vert^{d} \times \cdots \times \rho\vert\,\vert^{d}\times_{i\in [i_{0},i_{1}]}<d-1,f_{i}>_{\rho},
 $$
o\`u il y a $i_{1}-i_{0}+1$ facteurs $\rho\vert\,\vert^{d}$. On obtient:
$$
\times_{i\in [1,v]}<d_{i},f_{i}>_{\rho}\hookrightarrow \times_{i\in [1,i_{0}[}\rho\vert\,\vert^{d_{i}}\times \rho\vert\,\vert^{d}\times \cdots \times \rho\vert\,\vert^{d}\times_{i\in [i_{0},i_{1}]} <d-1,f_{i}>_{\rho} \times_{i>i_{1}}<d_{i},f_{i}>_{\rho}.
$$
Il suffit d'appliquer la r\'ecurrence \`a $\times_{i\in [i_{0},i_{1}]} <d-1,f_{i}>_{\rho} \times_{i>i_{1}}<d_{i},f_{i}>_{\rho}$ pour obtenir le r\'esultat annonc\'e. On a donc montr\'e qu'il existe une repr\'esentation $\sigma$ et une inclusion
$$
\pi\hookrightarrow \times_{x\in \mathcal{E}}\rho\vert\,\vert^x \times \sigma.
$$
Par r\'eciprocit\'e de Frobenius le module de Jacquet de $\pi$ contient un quotient de la forme $\otimes_{x\in \mathcal{E}}\rho\vert\,\vert^{x}\otimes \sigma$ et cela force la  non nullit\'e de $Jac_{x\in \mathcal{E}}\pi$ comme annonc\'e.

\begin{rem}Le lemme pr\'ec\'edent est une propri\'et\'e extr\^emement particuli\`ere car les \'el\'ements de $\mathcal{E}$ sont rang\'es dans l'ordre d\'ecroissant.
\end{rem}
\section{Premi\`ere d\'efinition des paquets d'Arthur  et comparaison de certaines orbites unipotentes}
\subsection{D\'efinition des paquets \`a la Arthur\label{premieredefinition} }
Ici on rappelle comment Arthur a d\'efini ses paquets. On consid\`ere un morphisme $\psi$ comme pr\'ec\'edemment  et on a donc d\'ej\`a d\'efini le morphisme de Langlands $\phi_{\psi}$ que l'on voit comme une repr\'esentation semi-simple de $W_{F}\times SL(2,{\mathbb C})$ de dimension $m_{G}^*$. Ce morphisme de Langlands d\'efinit une repr\'esentation irr\'eductible de $GL(m_{G}^*,F)$ via la correspondance de Langlands (\cite{harris}, \cite{henniart}); on la note $\pi^{GL}(\psi)$. Elle est tr\`es facile \`a d\'ecrire en terme explicite: pour $(\rho,a,b)\in Jord(\psi)$, on note $Speh (St(\rho,a),b)$ la repr\'esentation not\'ee $J(St(\rho,a),(b-1)/2,-(b-1)/2)$ en \ref{jacquetpourlanglands}.

Pour les $G$ que l'on consid\`ere ici, sauf les groupes orthgonaux pairs et leurs variantes, Arthur a remarqu\'e que ce groupe se r\'ealise comme le groupe endoscopique principal pour le produit semi-direct de $GL(m^*_{G})\times \{1,\theta\}$ (o\`u $\theta$ est un automorphisme ext\'erieur) pour les groupes unitaires, l'\'ecriture est un peu diff\'erente, il s'agit alors de la situation du changement de base. Dans le cas des groupes orthogonaux pairs, deux difficult\'es apparaissent, la non connexit\'e, prise en compte par Arthur en controlant l'action de l'automorphisme ext\'erieur et le fait que ce ne soit plus qu'un groupe endoscopique elliptique ce qui introduit des facteurs de transfert (compris explicitement). Dans tous les cas, Arthur annonce/d\'emontre l'existence d'un ensemble fini de repr\'esentations $\Pi(\psi)$ irr\'eductibles (Arthur l'exprime en demandant de longueur finie mais c'est \'equivalent en d\'ecomposant) tel que la trace tordue de $\pi^{GL}(\psi)$ soit un transfert d'une combinaison lin\'eaire convenable des traces des \'el\'ements de $\Pi(\psi)$, cette combinaison lin\'eaire \'etant stable; ce r\'esultat est cons\'equence de toute la stabilisation des formules des traces pour le groupe $GL$ tordu et tous ces groupes endoscopiques elliptiques r\'ealis\'ee par Arthur.  Ainsi cette formulation d'Arthur donne le calcul de la trace de cette combinaison lin\'eaire et d\'etermine uniquement les \'el\'ements de $\Pi(\psi)$; ici on n'a pas besoin de conna\^{\i}tre les coefficients. 

On a montr\'e en \cite{paquetdiscret} et \cite{pourshahidi} qu'en fait, pour d\'eterminer $\Pi(\psi)$ il suffit d'avoir le r\'esultat d'Arthur pour les repr\'esentations temp\'er\'ees c'est-\`a-dire pour les $\psi$ triviaux sur la 2e copie de $SL(2,{\mathbb C})$. On revient sur ces constructions dans \ref{rappelconstruction}. 
\subsection{Un lemme sur les modules de Jacquet\label{lemmesurjacquet}}
On fixe $\psi$ d'o\`u $\pi^{GL}(\psi)$. On fixe aussi un ensemble de demi-entiers $\mathcal{E}$ ordonn\'es par l'ordre d\'ecrois\-sant. Soit $\pi\in \Pi(\psi)$.
\begin{lem}On suppose que $Jac_{x\in \mathcal{E}}\pi\neq 0$ alors $Jac^\theta_{x\in \mathcal{E}}\pi^{GL}(\pi)\neq 0$.
\end{lem}
On ne d\'emontre pas ce lemme ici car pour l'avoir en toute g\'en\'eralit\'e, il faut quand m\^eme avoir quel\-ques renseignements sur les \'el\'ements de $\Pi(\psi)$. Toutefois, il y a un cas o\`u le lemme r\'esulte imm\'ediatement des r\'esultats d'Arthur sans passer par une description plus pr\'ecise. En effet, on sait que le transfert commute \`a la restriction y compris \`a la restriction partielle telle que d\'efinit ici; cela est v\'erifi\'e en \cite{selecta} 4.2. Si on sait que l'hypoth\`ese du lemme entra\^{\i}ne que $Jac_{x\in \mathcal{E}} $ de la distribution stable associ\'ee \`a $\psi$ et $G$ est non nulle, on aura bien que $Jac^\theta_{x\in \mathcal{E}}\pi^{GL}(\psi)\neq 0$. Or un cas \'evident entra\^{\i}ne cela, c'est le cas o\`u il ne peut pas y avoir de simplification c'est \`a dire le  cas est  o\`u tous les coefficients de la combinaison lin\'eaire stable se transferant en la trace tordue de $\pi^{GL}(\psi)$ ont  m\^eme signe. Arthur a d\'emontr\'e que ceci se produit si pour tout $(\rho,a,b)\in Jord(\psi)$, $b$ est impair. On ne d\'etaille pas puisque l'on fera une d\'emonstration plus g\'en\'erale plus loin. Mais j'esp\`ere que cela montre que le lemme n'est pas tr\`es profond et c'est ce lemme qui va donner le r\'esultat cherch\'e. Dans cette d\'emonstration (qui n'est que partielle), on n'a pas utilis\'e le fait que les \'el\'ements de $\mathcal{E}$ sont ordonn\'es par ordre d\'ecroissant; la d\'emonstration g\'en\'erale utilise l'hypoth\`ese mais je ne suis pas s\^ure qu'elle soit indispensable. Mais elle est indispensable pour le corollaire ci-dessous qui est ce que nous utiliserons.
\begin{cor}On suppose que $Jac_{x\in \mathcal{E}}\pi\neq 0$, alors dans le groupe lin\'eaire convenable $$Jac_{x\in \mathcal{E}}\times_{(\rho,a,b)\in Jord(\psi)}St(\rho,a)\vert\,\vert^{-(b-1)/2}\neq 0,$$ ici on consid\`ere bien toute l'induite qui n'est pas irr\'eductible en g\'en\'eral.
\end{cor}
En effet, on rappelle la forme particuli\`erement simple des repr\'esentations $\pi^{GL}(\psi)$; pour cela on utilise la notation matricielle pour des multisegments:
$$
Speh(St(\rho,a),b)= \bigg(\begin{matrix}
(a-b)/2 & (a-b)/2-1 &\cdots &-(a+b)/2+2 &-(a+b)/2+1\\
\vdots &\vdots &\cdots &\vdots &\vdots\\
(a+b)/2-1&(a+b)/2-2 &\cdots & -(a-b)/2+1 &-(a-b)/2
\end{matrix}\biggr)
_{\rho}
$$Ce qui est important est que les lignes sont des segments d\'ecroissants et les colonnes des segments croissants.
Et $\pi^{GL}(\psi)$ est l'induite irr\'eductible de ces repr\'esentations quand $(\rho,a,b)$ parcourt $Jord(\psi)$. Soit $\mathcal{Y}$ un ensemble de demi-entiers ordonn\'es de fa\c{c}on d\'ecroissante et supposons que $Jac_{y\in \mathcal{Y}}Speh(St(\rho,a),b)\neq 0$; alors n\'ecessairement $\mathcal{Y}$ est un sous-segment de $[(a-b)/2, -(a+b)/2+1]$: en effet \'ecrivons $\mathcal{Y}=\cup_{j\in [1,\ell]}y_{j}$ avec $y_{j'}\leq y_{j}$ si $j'<j$. Alors $Jac_{y_{1}}Speh(St(\rho,a),b)\neq 0$ entra\^{\i}ne que $y_{1}=(a-b)/2$; ensuite $Jac_{y_{2}}Jac_{y_{1}}Speh(St(\rho,a),b)\neq 0$ entra\^{\i}ne que $y_{2}=(a-b)/2-1$ ou $(a-b)/2+1$ mais la propri\'et\'e d'ordre assure que seule la premi\`ere propri\'et\'e est vraie et on obtient le r\'esultat annonc\'e de proche en proche. Soit $\mathcal{E}$ comme dans l'\'enonc\'e. Comme $Jac_{x\in \mathcal{E}}\pi^{GL}(\psi)\neq 0$, il existe un d\'ecoupage de $\mathcal{E}=\cup \mathcal{E}_{\rho,a,b}$ en sous-ensembles ordonn\'es par l'ordre induit et index\'es par les \'el\'ements de $Jord(\psi)$ tel que pour tout $(\rho,a,b)\in Jord(\psi)$, on ait
$$
Jac_{x\in \mathcal{E}_{\rho,a,b}}Speh(St(\rho,a),b)\neq 0.
$$
D'apr\`es la remarque pr\'ec\'edente cela veut dire que $\mathcal{E}_{\rho,a,b}$ est un sous-segment de $[(a-b)/2,-(a+b)/2+1]$ de m\^eme d\'ebut. Une autre fa\c{c}on de dire les choses est de dire que 
$$
Jac_{x\in \mathcal{E}}\times_{(\rho,a,b)\in Jord(\psi)}St(\rho,a)\vert\,\vert^{-(b-1)/2}\neq 0.$$
D'o\`u le corollaire.
\subsection{Comparaison des orbites unipotentes\label{preuveassertionclozel}}
On fixe $\psi$ et $\psi'$ des morphismes comme pr\'ec\'edemment; on fixe aussi une repr\'esentation irr\'eductible $\rho$ de $W_{F}$. On note $O^{unip}_{\psi,\rho}$ l'orbite unipotente du groupe lin\'eaire $GL(\sum_{(\rho,a,b)\in Jord(\psi)}ab,{\mathbb C})$ dont les blocs de Jordan sont exactement $\cup_{(\rho,a,b)\in Jord(\psi)}\underbrace{b, \cdots, b}_{a}$. On d\'efinit de m\^eme $O^{unip}_{\psi',\rho}$. 
\begin{rem}On suppose que $\Pi(\psi)\cap \Pi(\psi')\neq \emptyset$; on d\'emontrera en \ref{supportcuspidal} que les restrictions de $\psi$ et de $\psi'$ \`a $W_{F}$ fois la diagonale de $SL(2,{\mathbb C})\times SL(2,{\mathbb C})$ sont conjugu\'es. Donc sous cette hypoth\`ese, pour tout $\rho$ comme ci-dessus, $O^{unip}_{\psi,\rho}$ et $O^{unip}_{\psi',\rho}$ sont relatives au m\^eme groupe lin\'eaire. 
\end{rem}

\begin{thm}On suppose que $\Pi(\psi)\cap \Pi(\psi'_{L})\neq \emptyset$. Alors, pour toute repr\'esentation $\rho$ irr\'eductible de $W_{F}$, l'orbite $O^{unip}_{\psi,\rho}$ contient dans sa fermeture l'orbite $O^{unip}_{\psi',\rho}$.
\end{thm}

On consid\`ere $Jord(\psi')$ et on ordonne les $(\rho,a',b')\in Jord(\psi')$ tel que $b'>1$ par l'ordre d\'ecroissant sur $b'$; c'est \`a dire que l'on pose $t:=\vert \{(\rho,a',b')\in Jord(\psi'); b'>1\}\vert$ et que l'on \'ecrit les \'el\'ements de $Jord(\psi')$ contenant $\rho$, sous la forme $(\rho,a'_{j},b'_{j})$ pour $j\in [1,t]$ de sorte que pour $j<j'$, $b_{j}\geq b_{j'}$. Pour tout $j\in [1,t]$, on pose $\mathcal{E}_{\leq j}$, l'ensemble des \'el\'ements inclus dans l'union de segments $\cup_{i\leq j}[(a'_{i}-b'_{i})/2,-(a'_{i},b'_{i})/2+1]$, ensemble consid\'er\'e avec multiplicit\'e et r\'eordonn\'e pour que les \'el\'ements apparaissent par ordre d\'ecroissant. On a v\'erifi\'e en \ref{jacquetpourlanglands} que l'hypoth\`ese de \ref{lemmesurjacquet} est v\'erifi\'e pour $\mathcal{E}=\mathcal{E}_{\leq j}$ pour tout $j\in [1,t]$. On fixe $j$ et on applique le corollaire de \ref{lemmesurjacquet}. On sait que
$$
Jac_{x\in \mathcal{E}_{\leq j}}\times_{(\rho,a,b)\in Jord(\psi)}St(\rho,a)\vert\,\vert^{-(b-1)/2}\neq 0,
$$
dans le groupe lin\'eaire convenable (c'est-\`a-dire celui de rang $(d_{\rho}(\sum_{(\rho,a,b)\in Jord(\psi)}ab)$). Cela force l'existence d'un d\'ecoupage de chaque segment $[(a-b)/2,-(a+b)/2+1]$ en 2 sous-segments, $D'_{\rho,a,b}\cup D''_{\rho,a,b}$ dont l'un peut \^etre vide et dont le premier, s'il est non vide, d\'ebute par $(a-b)/2$ et tels que l'on ait l'\'egalit\'e d'ensembles avec multiplicit\'e mais non ordonn\'es
$$
\mathcal{E}_{\leq j}=\cup_{\rho,a,b}D'_{\rho,a,b}.\eqno(1)
$$On fixe un tel d\'ecoupage.
On calcule $2\sum_{x\in \mathcal{E}_{\leq j}}x$ de 2 fa\c{c}ons en utilisant l'\'egalit\'e ci-dessus. D'abord on utilise le fait $\mathcal{E}_{\leq j}$ est l'union des segments $[(a'_{i}-b'_{i})/2,-(a'_{i}+b'_{i})/2+1]$ pour $i\leq j$, chacun de ces segments est de milieu $-(b'_{i}-1)/2$, d'o\`u
$$
2\sum_{x\in \mathcal{E}_{\leq j}}x=\sum_{i\leq j}-(b'_{i}-1)a'_{i}.\eqno(2)
$$
On utilise maintenant le membre de droite de (1); fixons $(\rho,a,b)$ tel que $D'_{\rho,a,b}\neq \emptyset$ et notons $m(\rho,a,b,\leq j)$ le nombre d'\'el\'ements de cet intervalle. Ce nombre est s\^urement inf\'erieur ou \'egal \`a $a$. Comme $D'_{\rho,a,b}$ est un sous-segment du segment d\'ecroissant $[(a-b)/2,-(a+b)/2+1]$ commen\c{c}ant par $(a-b)/2$, on a encore que 2 fois la somme des \'el\'ements de ce sous-segment est  sup\'erieur ou \'egal au nombre d'\'el\'ements du sous-segment multipli\'e par 2 fois la valeur du milieu du segment c'est-\`a-dire:
$$
2\sum_{x\in D'_{\rho,a,b}}x\geq -m(\rho,a,b,\leq j)(b-1).\eqno(3)
$$
D'o\`u l'in\'egalit\'e:
$$
-\sum_{i\leq j}a'_{i}(b'_{i}-1)\geq -\sum_{\rho,a,b}m(\rho,a,b,\leq j)(b-1).
$$
Comme n\'ecessairement $$\sum_{i\leq j}a'_{i}=\sum_{\rho,a,b}m(\rho,a,b,\leq j)\eqno(*)$$, ce qui est une propri\'et\'e non \'evidente a priori et qui r\'esulte de (1), on obtient
$$
\sum_{\rho,a,b}m(\rho,a,b,\leq j)b\geq \sum_{i\leq j}a'_{i}b'_{i}.
$$
A droite on a la somme des $\sum_{i\leq j}a'_{i}$ plus grands blocs de Jordan de l'orbite $O^{unip}_{\psi',\rho}$; \`a gauche on a la somme de $\sum_{i\leq j}a'_{i}$ blocs de Jordan de l'orbite $O^{unip}_{\psi,\rho}$. On ordonne les blocs de Jordan de l'orbite $O^{unip}_{\psi,\rho}$ sous la forme $b_{1}\geq \cdots \geq b_{T}\geq 0$, chacun intervenant avec sa multiplicit\'e et on ajoute autant de $0$ qu'il faudra (ci-dessous on aura besoin que $T\geq T'$ o\`u $T'$ est d\'efini apr\`es (4));  on a donc, a fortiori, pour tout $j\in [1,t]$
$$
\sum_{k\leq \sum_{i\leq j}a'_{i}}b_{k}\geq \sum_{i\leq j}a'_{i}b'_{i}.\eqno(4)
$$
Ces in\'egalit\'es suffisent pour prouver le th\'eor\`eme. En effet, on pose $T':=\sum_{i\in [1,t]}a'_{i}+\vert \{(\rho,a',b')\in Jord(\psi');b'=1\}\vert$; c'est le nombre de blocs de Jordan de $O^{unip}_{\psi',\rho}$. On doit d'abord montrer que pour tout $\ell \leq T'$, la somme, not\'ee $S'_{ \ell}$ des $\ell$ plus grands blocs de Jordan de $O^{unip}_{\psi',\rho}$ est inf\'erieure ou \'egale \`a $\sum_{k\leq \ell}b_{k}$. On fixe $j$ tel que $\ell\in ]\sum_{i<j}a'_{i},\sum_{i\leq j}a'_{i}]$ ou $j=t+1$ si $\ell\in ]\sum_{i\leq t}a'_{t},T']$. Si $k=\sum_{i\leq j}a'_{i}$, on conna\^{\i}t le r\'esultat; on conna\^{\i}t aussi le r\'esultat si $k=T'$, dans ce cas il y a \'egalit\'e puisque $O^{unip}_{\psi,\rho}$ et $O^{unip}_{\psi',\rho}$ sont relatives au m\^eme groupe. Pour $\ell$ g\'en\'eral et $j\leq t$, on  a
$$
S'_{\ell}=(\sum_{i<j}a'_{i}b'_{i})+ (a'_{j}-r)b'_{j}=(\sum_{i\leq j}a'_{i}b'_{i})-r b'_{j}.
$$
D'o\`u,  $S'_{\ell}-\sum_{k\leq \ell}b_{k}=\biggl((\sum_{i<j}a'_{i}b'_{i})-(\sum_{k\leq \sum_{i<j}a'_{i}}b_{k}
)\biggr)+\sum_{k\in ]\sum_{i<j}a'_{i},\ell]}(b'_{j}-b_{k}).
$
Si $b'_{j}\leq b_{\ell}$, a fortiori, on a $b'_{j}\leq b_{k}$ pour tout $k\leq \ell$ et le terme de droite est n\'egatif ou nul d'apr\`es (4) appliqu\'e \`a $j-1$ et non $j$. Mais on a aussi, si $j\leq t$
$$
S'_{\ell}-\sum_{k\leq \ell}b_{k}=\biggl((\sum_{i\leq j}a'_{i}b'_{i})-(\sum_{k\leq \sum_{i\leq j}a'_{i}}b_{k}
)\biggr)-\sum_{k\in ]\ell,\sum_{i\leq j}a'_{i}]}(b'_{j}-b_{k});
$$
Supposons que $b'_{j}\geq b_{\ell}$, a fortiori, $b'_{j}\geq b_{k}$ pour tout $k\geq \ell$ et le terme de droite est encore n\'egatif ou nul d'apr\`es (4) appliqu\'e \`a $j$. 
Il reste le cas o\`u $j=t+1$; on a encore  2 \'egalit\'es du m\^eme type. La premi\`ere est
$$
S'_{\ell}-\sum_{k\leq \ell}b_{k}=\biggl((\sum_{i\leq t}a'_{i}b'_{i})-(\sum_{k\leq \sum_{i\leq t}a'_{i}}b_{k}
)\biggr)+\sum_{k\in ]\sum_{i<j}a'_{i},\ell]}(b'_{j}-b_{k}),
$$
qui r\`egle le cas o\`u $b'_{\ell}=1\leq b_{\ell}$, c'est-\`a-dire $b_{\ell}>0$. Si $b_{\ell}=0$, on a $\sum_{k\leq \ell}b_{k}=\sum_{(\rho,a',b')\in Jord(\psi')}a'b'\geq S'_{\ell}$ puisque l'on a vu que $O^{unip}_{\psi,\rho}$ est une orbite du m\^eme groupe que $O^{unip}_{\psi',\rho}$. Cela termine la preuve.
\subsection{Remarque}
La d\'emonstration de \ref{preuveassertionclozel} est une cons\'equence imm\'ediate du r\'esultat de \ref{lemmesurjacquet} qui est en fait plus puissant et en particulier vrai pour toute repr\'esentation dans $\Pi(\psi)$. Toutefois \ref{lemmesurjacquet} est technique alors que \ref{preuveassertionclozel} est plus simple \`a comprendre. Voil\`a les autres corollaires que l'on peut tirer de \ref{lemmesurjacquet} pour toute repr\'esentation $\pi$ dans $\Pi(\psi)$.
On fixe seulement $\psi$ ici et soit $\pi\in \Pi(\psi)$. Soient $v$ un entier et pour $i\in [1,v]$ des couples d'entiers, $(\alpha_{i},\beta_{i})$. On suppose que $1<\beta_{1}\leq \cdots \leq \beta_{v}$. On suppose aussi que pour tout $j\in [1,v]$, il existe une repr\'esentation $\sigma_{j}$ convenable, on a une inclusion
$$
\pi\hookrightarrow \times_{i\in [j,v]}St(\rho,\alpha_{i})\vert\,\vert^{-(\beta_{i}-1)/2}\times \sigma_{j}.\eqno(1)
$$
On remarque que les exposants $-(\beta_{i}-1)/2$ vont donc dans l'ordre d\'ecroissant ce qui est la condition oppos\'ee \`a celle des Langlands. Toutefois, si $\pi$ n'est pas temp\'er\'e, on peut trouver une telle donn\'ee au moins pour $v=1$. En particulier, on appelle exposant $\rho$ maximal de $\pi$, le plus grand entier strictement sup\'erieur \`a 1 (s'il existe) tel qu'il existe $\alpha$ un entier et $\sigma$ une repr\'esentation avec une inclusion
$$
\pi\hookrightarrow St(\rho,\alpha)\vert\,\vert^{-(\beta-1)/2}\times \sigma.
$$
On revient \`a la situation ci-dessus.
On consid\`ere l'orbite unipotente dont les blocs de Jordan sont exactement les $\beta_{i}$ avec multiplicit\'e $\alpha_{i}$. Alors, $\pi$ \'etant dans $\Pi(\psi)$, on a:
\begin{rem}L'orbite unipotente d\'efinie ci-dessus est une orbite d'un groupe $GL(m',{\mathbb C})$ avec $m'\leq \sum_{(\rho,a,b)\in Jord(\psi)}$. On la consid\`ere par inclusion comme une orbite du groupe $GL(\sum_{(\rho,a,b)\in Jord(\psi)}ab)$ et elle est dans la fermeture de l'orbite unipotente $O^{unip}_{\psi,\rho}$. En particulier, si $\pi$ n'est pas temp\'er\'ee, l'exposant $\rho$ maximal de $\pi$, $\beta$ v\'erifie $\beta\leq sup_{(\rho,a,b)\in Jord(\psi)}b$.
\end{rem}
C'est exactement la d\'emonstration qui a \'et\'e faite, en particulier la premi\`ere assertion r\'esulte de (*) ci-dessus. On remarque d'ailleurs que la condition forte que (1) doit \^etre vrai pour tout $j\in [1,v]$ est indispensable: supposons par exemple que $(\rho,a,b)\in Jord(\psi)$ avec $b$ maximum pour cette propri\'et\'e. On suppose que pour ce choix  $a<b$ et soit $\pi\in \Pi(\psi_{L})$. Il existe alors $\sigma$ et une inclusion $\pi\hookrightarrow St(\rho,a)\vert\,\vert^{(b-1)/2}\times \sigma$. Cette inclusion donne a fortiori
$$
\pi \hookrightarrow \times_{x\in [(b-a)/2, (a+b)/2-1]}\rho\vert\,\vert^{-x}\times \sigma;
$$
mais l'orbite de blocs de Jordan $b-a+1, \cdots, a+b-1$ et autant de $1$ que n\'ecessaire n'est pas dans la fermeture de l'orbite $O^{unip}_{\psi,\rho}$ dont le plus grand bloc de Jordan est $b$.

\section{Rappel de la construction des repr\'esentations dans un paquet d'Arthur\label{rappelconstruction}}
La construction des repr\'esentations dans les paquets d'Arthur se fait par r\'ecurrence; on en rappelle ici les grandes lignes pour rendre l'article ind\'ependant d'autres r\'ef\'erences mais les r\'ef\'erences pour les preuves sont \cite{elementaire}, \cite{paquetdiscret}, \cite{holomorphie}. On se ram\`ene d'abord au cas o\`u pour tout $(\rho,a,b)\in Jord(\psi)$ la repr\'esentation $\rho\otimes rep_a\otimes rep_b$ est \`a valeurs dans un groupe de m\^eme type que $G^*$; exactement on note $\psi_{bp}$ la somme des repr\'esentations incluses dans $\psi$ et ayant cette propri\'et\'e et on note $\psi_{mp}$ la somme des autres; on v\'erifie que $\psi_{mp}$ s'\'ecrit sous la forme $\psi_{1/2,mp}\oplus \theta^*(\psi_{1/2,mp})$ pour un choix non unique d'une sous-repr\'esentation, $\psi_{1/2,mp}$ dans $\psi_{mp}$. On admet savoir construire $\Pi(\psi_{bp})$ et on a montr\'e en \cite{general} repris en \cite{pourshahidi} 3.2, que pour tout $\pi_{pb}\in \Pi(\psi_{bp})$, l'induite $\pi^{GL}(\psi_{1/2,mp})\times \pi_{bp}$ est irr\'eductible et que $\Pi(\psi)$ est constitu\'e exactement de ces induites. En prouvant cette irr\'eductibilit\'e on donne aussi les param\`etres de Langlands de $\pi$ en fonction de ceux de $\pi_{bp}$; pr\'ecis\'ement, on montre que 
$$
\pi\hookrightarrow_{(\rho,a,b)\in Jord(\psi_{1/2,mp}),b>1} \biggl(J(St(\rho,a), -(b-1)/2,-\delta_{b})\times J(St(\tilde{\rho},a),-(b-1)/2,-\delta_{b})\biggr)$$
$$\times_{(\rho,a,b)\in Jord(\psi_{1/2,mp}); b\equiv 1[2]}St(\rho,a)\times \pi_{bp},\eqno(1)
$$
o\`u la recette pour obtenir $\tilde{\rho}$ en fonction de $\rho$ d\'epend de $G$ (pour les groupes classiques standard $\tilde{\rho}$ est la contragr\'ediente de $\rho$). On montre en plus que l'induite $\times_{(\rho,a,b)\in Jord(\psi_{1/2,mp}); b\equiv 1[2]}St(\rho,a)\times \pi_{bp}$ est irr\'eductible elle aussi et que les facteurs $St(\rho,a)$ qui y interviennent commutent avec les $St(\rho_{i},\alpha_{i})\vert\,\vert^{x}$ qui interviennent comme param\`etre de Langlands de $\pi_{bp}$.  

\

Ensuite on construit $\Pi(\psi)$ quand $\psi=\psi_{bp}$. On fait cette construction en 2 \'etapes. D'abord on consid\`ere le cas o\`u pour tout $(\rho,a,b),(\rho,a',b')\in Jord(\psi)$, (le m\^eme $\rho$), $[\vert (a-b)\vert+1,a+b-1]\cap [\vert (a'-b')\vert +1,a'+b'-1]=\emptyset$.  On dit que $\psi$ est de restriction discr\`ete \`a la diagonale car c'est le cas o\`u la restriction de $\psi$ \`a $W_{F}$ fois la diagonale de $SL(2,{\mathbb C})\times SL(2,{\mathbb C})$ est sans multiplicit\'e. La construction est faite dans \cite{paquetdiscret}; on rappelle qu'elle se fait par r\'ecurrence sur $\ell(\psi):=\sum_{(\rho,a,b)\in Jord(\psi)}(inf(a,b)-1)$.
Quand $\ell(\psi)=0$ la construction a \'et\'e faite en \cite{elementaire} en se ramenant via une g\'en\'eralisation de l'involution d'Iwahori-Matsumoto au cas des s\'eries discr\`etes; le cas des s\'eries discr\`etes est d\^u \`a Arthur et compl\'et\'e pour les propri\'et\'es dont nous avons besoin comme expliqu\'e dans \cite{pourshahidi}. Le cas g\'en\'eral est fait en \cite{paquetdiscret}. Les \'el\'ements de $\Pi(\psi)$ sont en bijection avec les couples d'applications $\underline{t},\underline{\eta}$ de $Jord(\psi)$ \`a valeurs dans ${\mathbb N}\times \{\pm 1\}$ soumises aux conditions suivantes:
$$
\underline{t}(\rho,a,b)\in [0,[inf(a,b)/2]]
$$
et $\underline{\eta}(\rho,a,b)=+$ si $\underline{t}(\rho,a,b)=[inf(a,b)/2]$ et le signe$$\times_{(\rho,a,b)\in Jord(\psi)}(\underline{\eta}(\rho,a,b)^{inf(a,b)}(-1)^{[inf(a,b)/2]+\underline{t}(\rho,a,b)}
)$$
est d\'etermin\'e par la forme de $G$.
On rappelle ici un cas simple, le seul o\`u on conna\^{\i}t vraiment les param\`etres de Langlands, o\`u $\psi$ en plus d'\^etre de restriction discr\`ete \`a la diagonale,  v\'erifie que pour tout $(\rho,a,b)\in Jord(\psi)$, $a\geq b$. On a donc les param\`etres $\underline{t},\underline{\eta}$; on construit d'abord un morphisme, $\psi_{temp}$ de $W_{F}\times SL(2,{\mathbb C})$ dont la d\'ecomposition en repr\'esentations irr\'eductibles est $$\oplus_{(\rho,a,b)\in Jord(\psi)} \oplus_{c\in [a-b+1+2t,a+b-1-2t]; (-1)^c=(-1)^{a+b-1}}\rho \otimes rep_c.$$ A l'aide de $\underline{\eta}$, on d\'ecrit un caract\`ere $\eta_{temp}$ du centralisateur de $\psi_{temp}$ (on ne le fait pas ici, car on n'en aura pas besoin); on note $\pi_{temp}$ la repr\'esentation temp\'er\'ee (ici une s\'erie discr\`ete) associ\'ee \`a ces donn\'ees. Alors la repr\'esentation de $\Pi(\psi)$ associ\'ee \`a $\underline{t}$ et $\underline{\eta}$ est l'unique sous-module irr\'eductible de l'induite
$$
\pi\hookrightarrow \times_{(\rho,a,b)\in Jord(\psi); \underline{t}(\rho,a,b)>0}St(\rho,a)\vert\,\vert^{-(b-1)/2}\times \cdots \times St(\rho,a)\vert\,\vert^{-(b-1)/2+\underline{t}(\rho,a,b)}\times \pi_{
temp}.
$$
Il n'est pas utile de mettre d'ordre sur $Jord(\psi)$ pour construire l'induite, c'est la condition de restriction discr\`ete \`a la diagonale qui assure que l'on peut prendre n'importe quel ordre sans changer le r\'esultat; c'est d'ailleurs aussi cette condition qui assure que l'induite a un unique sous-module irr\'eductible.

Dans le cas o\`u $\psi$ est de restriction discr\`ete \`a la diagonale, on a donc une tr\`es bonne param\'etrisation des \'el\'ements de $\Pi(\psi)$ m\^eme si en g\'en\'eral, on ne conna\^{\i}t pas leur param\`etre de Langlands. La construction est sym\'etrique entre les 2 copies de $SL(2,{\mathbb C})$ alors que la param\'etrisation de Langlands ne l'est \'evidemment pas.

\

Pour passer au cas g\'en\'eral, il faut fixer un ordre total, $>$, sur $Jord(\psi)$ qui doit satisfaire \`a la condition:

pour tout $(\rho,a,b),(\rho,a',b')\in Jord(\psi)$ (m\^eme $\rho$) tels que $(a-b)(a'-b')>0$, si $\vert (a-b)\vert /2>\vert (a'-b')\vert /2$ et si $(a+b)/2-1>(a'+b')/2-1$ alors $(\rho,a,b)>(\rho,a',b')$.

Un ordre qui v\'erifie cette condition sera dit {\bf un bon ordre}; on fait en plus un choix de signe $\zeta_{\rho,a,b}$ tel que $\zeta_{\rho,a,b}(a-b)\geq 0$; le choix n'est donc que pour les $(\rho,a,b)$ tel que $a-b=0$. On accepte le fait que $Jord(\psi)$ ait de la multiplicit\'e et l'ordre distingue donc entre des \'el\'ements de $Jord(\psi)$ qui sont \'egaux. La param\'etrisation de $Jord(\psi)$ d\'epend de l'ordre fix\'e mais \'evidemment pas l'ensemble $\Pi(\psi)$ lui-m\^eme (cf. \cite{holomorphie} 2.8). Fixons donc un bon ordre sur $Jord(\psi)$; on dit qu'un morphisme $\psi_{>}$ domine $Jord(\psi)$ si $Jord(\psi_{>})$ est lui aussi muni d'un bon ordre et qu'il existe une bijection de $Jord(\psi_{>})$ sur $Jord(\psi)$ respectant l'ordre (une telle bijection est uniquement d\'etermin\'e) et telle que pour tout $(\rho,a,b)\in Jord(\psi)$, il existe $T_{\rho,a,b}\in {\mathbb N}_{\geq 0}$ tel que l'image r\'eciproque de $(\rho,a,b)$ dans $Jord(\psi_{>})$ soit $(\rho,a+2T_{\rho,a,b},b)$ si $\zeta_{\rho,a,b}=+$ et $(\rho,a,b+2T_{\rho,a,b})$ si $\zeta_{\rho,a,b}=-$. On ne consid\`ere que le cas o\`u il existe $(\rho,a,b)\in Jord(\psi)$ tel que $T_{\rho,a',b'}=0$ pour tout $(\rho,a',b')\leq (\rho,a,b)$ et $T_{\rho,a',b'}>> T_{\rho,a'',b''}$ pour tout $(\rho,a',b')>(\rho,a'',b'')\geq (\rho,a,b)$. En particulier si $(\rho,a,b)$ est l'\'el\'ement minimal de $Jord(\psi)$ un tel morphisme sera dit tr\`es dominant et en particulier, il est n\'ecessairement de restriction discr\`ete \`a la diagonale.

Fixons un tel morphisme tr\`es dominant, not\'e $\psi_{>>}$. On sait donc d\'efinir $\Pi(\psi_{>>})$ et on a montr\'e d'abord en \cite{paquetdiscret} (avec un choix d'ordre particulier) et pour tous les bons ordres en \cite{holomorphie} que pour tout $\pi_{>>}\in \Pi(\psi_{>>})$, la repr\'esentation
$$
\circ_{(\rho,a,b)\in Jord(\psi)}\circ_{\ell \in [1,T_{\rho,a,b]}}Jac_{(a-b)/2+\ell, \cdots, \zeta_{\rho,a,b,>>}((a+b)/2-1+\ell)}\pi_{>>},
$$o\`u l'on prend les \'el\'ements de $Jord(\psi)$ dans l'ordre croissant (c'est-\`a-dire que le plus petit est le plus \`a gauche),
est irr\'eductible ou nulle. Et on a montr\'e que l'ensemble de ces repr\'esentations (quand on a enlev\'e celles qui donnent 0) d\'efinit $\Pi(\psi)$. Pour r\'esumer, on choisit l'ordre et les signes $\zeta_{\rho,a,b}$ quand c'est n\'ecesaire. Et alors pour tout $\pi\in \Pi(\psi)$, il existe des fonctions $\underline{t}(\rho,a,b)$ et $\underline{\eta}(\rho,a,b)$ qui v\'erifient exactement les propri\'et\'es ci-dessus mais qui a priori sont d\'efinies sur l'ensemble $Jord(\psi)$ vu avec multiplicit\'e pour obtenir toutes les fonctions possibles sur $Jord(\psi_{>>})$ (on a montr\'e en \cite{general} que l'on peut se limiter aux fonctions sur $Jord(\psi)$)); on les remonte en des fonctions sur $Jord(\psi_{>>})$ via la bijection de $Jord(\psi_{>>})$ sur $Jord(\psi)$, on construit $\pi_{>>}$ avec ces fonctions et on retrouve $\pi$ par module de Jacquet, \`a partir de $\pi_{>>}$; ce qui est important est que $\underline{t}$ et $\underline{\eta}$ ne d\'ependent que de l'ordre et des signes choisis pas de $\psi_{>>}$ dominant $\psi$ pour cet ordre.

Le point important dans cette construction est sa transitivit\'e: soit $\psi,\psi_{>}, \psi_{>>}$ des morphismes tels que $Jord(\psi)$ soit muni d'un bon ordre et que $\psi_{>}$ domine $\psi$ (d'o\`u $Jord(\psi_{>})$ est muni d'un bon ordre compatible) et tel que $\pi_{>>}$ domine $\psi_{>}$; ainsi $\psi_{>>}$ domine $\psi$ et en suivant les d\'efinitions on v\'erifie que  le compos\'e de l'application de $\Pi(\psi_{>>})$ sur $\Pi(\psi_{>})$ avec l'application de $\Pi(\psi_{>})$ sur $\Pi(\psi)$ est l'application de $\Pi(\psi_{>>})$ sur $\Pi(\psi)$; le seul point qui sert est que $Jac_{x}Jac_{y}=Jac_{y}Jac_{x}$ si $\vert x-y\vert >1$. En particulier on obtient la remarque suivante:
\begin{rem}On fixe $\psi$ d'o\`u $Jord(\psi)$ et $(\rho,a,b)\in Jord(\psi)$ et $\zeta_{a,b}=\pm 1$ tel que $\zeta_{a,b}(a-b)\geq 0$. On suppose que pour tout $(\rho,a',b')\in Jord(\psi)$ avec $\zeta_{a,b}(a'-b')>0$, $\vert(a'-b')/2\vert >\vert (a-b)/2\vert$ et $(a'+b')/2-1>(a+b)/2-1$, on a en plus $(a'-b')/2\geq (a+b)/2+1$ alors on note $\psi_{+}$ le morphisme qui se d\'eduit de $\psi$ en rempla\c{c}ant $(\rho,a,b)$ par $(\rho,a+2,b)$ si $\zeta_{a,b}=+$ et $(\rho,a,b+2)$ si $\zeta_{a,b}=-$. Pour tout $\pi_{+}\in \Pi(\psi_{+})$, $Jac_{(a-b)/2+1, \cdots, \zeta_{a,b}(a+b)/2}\pi_{+}$ est soit 0 soit un \'el\'ement irr\'eductible de $\Pi(\psi)$ et tout \'el\'ement de $\Pi(\psi)$ est obtenu ainsi pour un unique choix de $\pi_{+}$.
\end{rem}
On note $(\rho,a_{+},b_{+})$ l'\'el\'ement de $Jord(\psi_{+})$ qui vaut $(\rho,a+2,b)$ si $\zeta_{a,b}=+$ et $(\rho,a,b+2)$ sinon. 
On fixe un bon ordre sur $Jord(\psi_{+})$ tel que les \'el\'ements de $Jord(\psi_{+})$ strictement plus grand que $(\rho,a_{+},b_{+})$ sont exactement les \'el\'ements de la forme $(\rho,a',b')$ avec $(a'-b')/2 \geq (a+b)/2+1$; ceci est possible \'etant donn\'e les hypoth\`eses. Cet ordre induit un bon ordre sur $Jord(\psi)-\{(\rho,a,b)\}$ qui devient un bon ordre sur $Jord(\psi)$ en disant que $(\rho,a,b)<(\rho,a',b')$ si et seulement si $(a'-b')/2\geq (a+b)/2+1$. On fixe alors $\psi_{>>}$ dominant $\psi_{+}$ pour l'ordre fix\'e et on remarque que $\psi_{>>}$ domine alors aussi $\psi$ pour l'ordre fix\'e. Soit $\pi_{>>}\in \Pi(\psi_{>>})$; gr\^ace \`a cet \'el\'ement on construit un \'el\'ement de $\Pi(\psi)$ (ou 0) et un \'el\'ement de $\Pi(\psi_{+})$ (ou 0) en prenant des modules de Jacquet convenable. On commence par faire redescendre les \'el\'ements de $Jord(\psi_{>>})$ qui dominent les \'el\'ements de $Jord(\psi)$ strictement inf\'erieurs\`a $(\rho,a,b)$; ce sont aussi ceux qui dominent les \'el\'ements de $Jord(\psi_{+})$ strictement inf\'erieurs \`a $(\rho,a_{+},b_{+})$; l'op\'eration est la m\^eme pour $\psi$ et $\psi_{+}$; puis on fait redescendre l'\'el\'ement dominant $(\rho,a,b)$ et $(\rho,a_{+},b_{+})$; ici il y a une diff\'erence car on est arriv\'e \`a $(\rho,a_{+},b_{+})$ il faut encore prendre $Jac_{(a-b)/2+1, \cdots, \zeta_{a,b}(a+b)/2}$ du r\'esultat pour arriver \`a $(\rho,a,b)$. Puis on fait redescendre les \'el\'ements de $Jord(\psi_{>>})$ strictement sup\'erieur \`a $(\rho,a,b)$ et $(\rho,a_{+},b_{+})$ ce qui se fait via un module de Jacquet de la forme $Jac_{x\in {\mathcal E}}$ pour un ensemble ${\mathcal E}$ totalement ordonn\'e. Mais d'apr\`es les hypoth\`eses pour tout $x\in {\mathcal E}$ on a s\^urement $\vert x\vert >(a+b)/2+1$; d'o\`u pour tout $y\in [(a-b)/2+1,\zeta_{a,b}((a+b)/2)]$, on a $\vert x-y\vert >1$. Ainsi $Jac_{y\in [(a-b)/2+1,\zeta_{a,b}((a+b)/2)] }\circ Jac_{x\in \mathcal{E}}=Jac_{x\in \mathcal{E}}\circ Jac_{y\in [(a-b)/2+1,\zeta_{a,b}((a+b)/2)] }$. Ainsi l'\'el\'ement de $\Pi(\psi)$ obtenu gr\^ace \`a $\pi_{>>}$ s'il est non nul est exactement $Jac_{y\in [(a-b)/2+1,\zeta_{a,b}((a+b)/2)] }\pi_{+}$ o\`u $\pi_{+}\in \Pi(\psi_{+})$ est obtenu gr\^ace \`a $\pi_{>>}$; et cet \'el\'ement est nul exactement si soit $\pi_{>>}$ donne d\'ej\`a 0 dans $\Pi(\psi_{+})$ soit $Jac_{y\in [(a-b)/2+1,\zeta_{a,b}((a+b)/2)] }\pi_{+}=0$. Cela termine la preuve de la remarque.

\section{Support cuspidal \'etendu\label{supportcuspidal}}
\subsection{D\'efinition et calcul\label{definitionsupportcuspidal}}
On a d\'ej\`a utilis\'e cette d\'efinition de support cuspidal \'etendu dans plusieurs papiers  et on la rappelle ici.
Soit $\pi_{cusp}$ une repr\'esentation cuspidale d'un groupe de m\^eme type que $G$ \'eventuellement de rang plus petit. A cette repr\'esentation on associe un morphisme, $\phi_{cusp}$, de $W_{F}\times SL(2,{\mathbb C})$ dans un groupe $GL(m^*_{cusp},{\mathbb C})$ pour $m^*_{cusp}$ un entier convenable. On d\'ecompose cette repr\'esentation en sous-repr\'esentation irr\'eductibles $\oplus_{(\rho,a)\in Jord(\pi_{cusp}))}\rho\otimes rep_a$, o\`u $\rho$ parcourt  un ensemble (avec multiplicit\'e) de repr\'esentations irr\'eductibles de $W_{F}$ et $a$ parcourt un sous-ensemble (avec multiplicit\'e) de ${\mathbb N}$. Ceci d\'efinit $Jord(\pi)$. On appelle support cuspidal \'etendu de $\pi_{cusp}$ l'ensemble des repr\'esentations cuspidales $\cup_{(\rho,a)\in Jord(\pi))}\cup_{x\in [-(a-1)/2,(a-1)/2]}\rho\vert\,\vert^{x}$. Soit $\pi$ une repr\'esentation irr\'eductible de $G$; on \'ecrit $\pi$ comme sous-quotient irr\'eductible d'une induite de la forme $\times _{(\rho',z')}\rho'\vert\,\vert^{z'}\times \pi_{cusp}$, o\`u $\pi_{cusp}$ est une repr\'esentation cuspidale convenable et o\`u $(\rho',z')$ parcourt un ensemble de couples form\'es d'une repr\'esentation cuspidale unitaire irr\'eductible, $\rho'$, d'un groupe lin\'eaire convenable et $z'$ est un r\'eel. On appelle alors support cuspidal \'etendu de $\pi$ l'ensemble union du support cuspidal \'etendu de $\pi_{cusp}$ et de l'ensemble des repr\'esentations $\rho'\vert\,\vert^{z'},\rho^{'*}\vert\,\vert^{-z'}$ qui interviennent ci-dessus; le sens de $\rho^{',*}$ d\'epend du groupe, c'est essentiellement la repr\'esentation contragr\'ediente. Le support cuspidal \'etendu de $\pi$ est donc bien d\'efini comme ensemble non ordonn\'e (et avec multiplicit\'e).

Soit maintenant un morphisme, $\psi_{\vert}$ de $W_{F}\times SL(2,{\mathbb C})$ dans $GL(m_{G}^*,{\mathbb C})$ semi-simple; comme ci-dessus on le d\'ecompose en repr\'esentations irr\'eductibles $\oplus_{(\rho,c)\in Jord(\psi_{\vert})}\rho\otimes rep_c$, ce qui d\'efinit $Jord(\psi_{\vert})$ comme ensemble avec multiplicit\'e. On appelle support cuspidal de $\psi_{\vert}$ l'ensemble $\cup_{(\rho,c)Jord(\psi_{\vert})}\cup_{x\in [-(c-1)/2,(c-1)/2]}\rho\vert\,\vert^{x}$. Il est extr\^emement facile de retrouver $\psi_{\vert}$ \`a conjugaison pr\`es (\'evidemment) quand on conna\^{\i}t son support cuspidal. On fait remarquer au lecteur qu'un tel support cuspidal a une propri\'et\'e de sym\'etrie forte puisque dans l'union on consid\`ere des segments centr\'es en 0.

Soit $\psi$ ici un morphisme de $W_{F}\times SL(2,{\mathbb C})\times SL(2,{\mathbb C})$ dans $GL(m_{G}^*,{\mathbb C})$. On appelle support cuspidal de $\psi$ le support cuspidal de la restriction de $\psi$ \`a $W_{F}$ fois la diagonale de $SL(2,{\mathbb C})$. 

\begin{prop} Soit $\psi$ un morphisme de $W_{F}\times SL(2,{\mathbb C})\times SL(2,{\mathbb C})$ dans $GL(m_{G}^*,{\mathbb C})$ se factorisant par $G^*$. Soit  $\pi\in\Pi(\psi)$ alors le support cuspidal \'etendu de $\pi$ est le support cuspidal de $\psi$.
\end{prop}
On suit les d\'efinitions. Consid\'erons d'abord le cas de restrictrion discr\`ete \`a la diagonale; on traite d'abord le cas des paquets temp\'er\'ees qui ici sont n\'ecessairement des paquets de s\'eries discr\`etes.  Donc par hypoth\`ese pour tout $(\rho,a,b)\in Jord(\psi)$, $inf(a,b)=b=1$. On a montr\'e en \cite{europe} que si $\pi\in \Pi(\psi)$, trois cas sont possibles:

soit il existe $(\rho,a,b)$ avec $a>2 $ avec $(\rho,a-2,b)\notin Jord(\psi)$ et alors on notant $\psi'$ le morphisme qui se d\'eduit de $\psi$ en rempla\c{c}ant $(\rho,a,b)$ par $(\rho,a-2,b)$, il existe $\pi'\in \Pi(\psi')$ avec $\pi\hookrightarrow \rho\vert\,\vert^{(a-1)/2}\times \pi'$;

soit le cas pr\'ec\'edent n'est pas v\'erifi\'e mais il existe $(\rho,a,1)\in Jord(\psi)$ avec $a\geq 2$ tel qu'en notant $\psi'$ le morphisme tel que $Jord(\psi')$ se d\'eduit de $Jord(\psi)$ en enlevant $(\rho,a,1)$ et, si $a>2$, $(\rho,a-2,1)$, il existe $\pi'\in \Pi(\psi')$ tel que $\pi\hookrightarrow <(a-1)/2, \cdots, -(a'-1)/2>_{\rho}\times \pi'$;

soit $\pi$ est cuspidale.

On d\'emontre donc la proposition dans le cas des morphismes temp\'er\'es de restriction discr\`ete \`a la diagonale par r\'ecurrence sur le rang du groupe.

On passe du cas pr\'ec\'edent au cas des morphismes \'el\'ementaires en utilisant les d\'efinitions explicites de \cite{elementaire}, on y avait obtenu les repr\'esentations cherch\'ees en appliquant une g\'en\'eralisation de l'involution d'Iwahori Matsumoto qui ne change ni le support cuspidal des repr\'esentations (donc pas non plus le support cuspidal \'etendu) ni le support cuspidal des morphismes car on garde constant leur restriction \`a $W_{F}$ fois la diagonale de $SL(2,{\mathbb C})$. 

On passe au cas d'un morphisme de restriction discr\`ete \`a la diagonale; ce que l'on a rappel\'e en \ref{rappelconstruction} n'est pas tout \`a fait suffisant pour cela, il faut expliquer un peu plus. On fixe donc $\pi\in \Pi(\psi)$ d'o\`u ses param\`etres $\underline{t}$ et $\underline{\eta}$. On montre la proposition par r\'ecurrence sur $\ell(\psi):=\sum_{(\rho,a,b)\in Jord(\psi)}(inf(a,b)-1)$. Si $\ell(\psi)=0$, on est dans le cas \'el\'ementaire qui vient d'\^etre vu. Sinon, on prend $(\rho,a,b)\in Jord(\psi)$ tel que $inf(a,b)>1$; on a alors 2 possibilit\'es, 

\

\noindent
soit $\pi\in \Pi(\psi')$ o\`u $\psi'$ se d\'eduit de $\psi$ en rempla\c{c}ant $(\rho,a,b)$ par $\cup_{c\in [\vert a-b\vert+1, a+b-1]; (-1)^{c}=(-1)^{a+b-1}}(\rho,c,1)$, si $a\geq b$, $(\rho,1,c)$ si $a<b$;

\noindent
soit, on note $\psi'$ le morphisme tel que $Jord(\psi')$ se d\'eduit de $Jord(\psi)$ en rempla\c{c}ant $(\rho,a,b)$ par $(\rho,a,b-2)$ si $a\geq b$ et par $(\rho,a-2,b)$ si $a<b$ et il existe $\pi'\in \Pi(\psi')$ avec une inclusion  $\pi\hookrightarrow <(a-b)/2, \cdots, -\zeta(a+b)/2+1>_{\rho}\times \pi'$, o\`u $\zeta=+$ si $a\geq b$ et $-$ sinon.

\

Le premier cas est particuli\`erement facile cas on a remplac\'e la repr\'esentation associ\'e \`a $(\rho,a,b)$ par essentiellement sa restriction \`a $W_{F}$ fois la diagonale de $SL(2,{\mathbb C})\times SL(2,{\mathbb C})$ sans changer (\'evidemment) le support cuspidal de $\pi$;  dans le deuxi\`eme cas $$\psi_{\vert W_{F}\times \Delta_{SL(2,{\mathbb C})\times SL(2,{\mathbb C})}}= \psi'_{\vert W_{F}\times \Delta_{SL(2,{\mathbb C})\times SL(2,{\mathbb C})}} \oplus \rho\otimes [\vert a-b\vert+1] \oplus \rho\otimes [a+b-1]
$$ et
 $\{\pm x; x\in [(a-b)/2,-\zeta ((a+b)/2-1)]\}= [(a-b)/2,-(a-b)/2] \cup [(a+b)/2-1,-(a+b)/2+1]$; l'assertion pour $\pi'$ et $\psi'$ l'entra\^ine donc pour $\pi$ et $\psi$.

\

On traite maintenant le cas g\'en\'eral. On \'ecrit la d\'ecomposition en repr\'esentations irr\'eductibles de $\psi=\oplus_{(\rho,a,b)\in Jord(\psi)}\rho\otimes rep_a\otimes rep_b$. On suppose d'abord qu'il existe $(\rho,a,b)\in Jord(\psi)$ avec soit $\rho\not\simeq \rho*$ ou $(\rho,a,b)$ de mauvaise parit\'e et on fixe un tel triplet. On note $\psi'$ le morphisme qui se d\'eduit de $\psi$ en enlevant $\rho\otimes rep_a\otimes rep_b\oplus \rho^*\otimes rep_a\otimes rep_b$. On sait que tout \'el\'ement  de $\Pi(\psi)$ est sous-quotient d'une induite de la forme $Speh(St(\rho,a),b)\times \pi'$ o\`u $\pi'\in \Pi(\psi')$. Ainsi le support cuspidal \'etendu de $\pi$ s'obtient en ajoutant \`a celui de $\pi'$ l'ensemble $\cup_{d\in [-(b-1)/2,(b-1)/2, c\in [-(a-1)/2,(a-1)/2 } (\rho\vert\,\vert^{d+c},\rho^{*}\vert\,\vert^{d+c})$. Clairement $\cup_{d\in [-(b-1)/2,(b-1)/2),c\in [-(a-1)/2,(a-1)/2)}d+c=\cup_{e\in [\vert a-b\vert+1,a+b-1]_{2}}\cup_{x\in [-(e-1)/2,(e-1)/2}x$; l'indice 2 signifie que l'on ne prend que les entiers dans l'intervalle ayant m\^eme parit\'e que les bornes . Ainsi si l'on sait que le support cuspidal \'etendu de $\pi'$ co\"{\i}ncide avec le support cuspidal de $\psi'$, on obtient le m\^eme r\'esultat pour $\pi$ et $\psi$ puisque la restriction de la repr\'esentation $rep_a\otimes rep_b$ de $SL(2,{\mathbb C})\times SL(2,{\mathbb C})$ \`a sa diagonale est $\oplus_{e\in [\vert a-b\vert+1,a+b-1]_{2}}[e]$. On est donc ramen\'e au cas o\`u le morphisme $\psi$ est de bonne parit\'e. On fixe $\psi_{>>}$ dominant $\psi$, de restriction discr\`ete \`a la diagonale. On a fix\'e $\pi\in \Pi(\psi)$; on sait qu'il existe $\pi_{>>}\in \Pi(\psi_{>>})$ et une inclusion
$$
\pi_{>>}\hookrightarrow \times_{(\rho',x')\in {\mathcal E}}\rho'\vert\,\vert^{x'}\times \pi,
$$
o\`u l'ensemble ${\mathcal E}$ est de la forme $\cup_{(\rho,a,b)\in Jord(\psi)}\cup_{\ell\in [1,T_{\rho,a,b}}\cup_{x\in [(a-b)/2+\zeta_{\rho,a,b}\ell, \zeta_{\rho,a,b}((a+b)/2-1+\ell)]}(\rho,x)$ o\`u $\zeta_{\rho,a,b}$ est un signe essentiellement celui de $a-b$ et o\`u  $T_{\rho,a,b}$ est l'entier fix\'e pour chaque $(\rho,a,b)$ tel que $Jord(\psi_{>>})=\cup_{(\rho,a,b)\in Jord(\psi)}(\rho,a+(1+\zeta_{\rho,a,b})T_{\rho,a,b},b+(1-\zeta_{\rho,a,b})T_{\rho,a,b})$. Le support cuspidal \'etendu de $\pi$ s'obtient donc en enlevant au support cuspidal \'etendu de $\pi_{>>}$ l'ensemble $\cup_{(\rho',x')\in {\mathcal E}} (\rho'\vert\,\vert^{x'},\rho'\vert\,\vert^{-x'})$. Et d'apr\`es la description des $T_{\rho,a,b}$ c'est exactement la m\^eme op\'eration qui fait passer du support cuspidal de $\psi_{>>}$ au support cuspidal de $\psi$: pour s'en convaincre, on regarde le cas d'une repr\'esentation $\rho\otimes rep_a\otimes rep_b$ apparaissant dans $\psi$ et de $\rho\otimes [a+(1+\zeta_{\rho,a,b})T_{\rho,a,b}]\otimes [b+(1-\zeta_{\rho,a,b})T_{\rho,a,b}]$ apparaissant dans $\psi_{>>}$. En restriction \`a la diagonale de $SL(2,{\mathbb C})\times SL(2,{\mathbb C})$, la premi\`ere repr\'esentation se d\'ecompose en $\oplus_{c\in [\vert a-b\vert +1,a+b-1]_{2}}\rho\otimes rep_c$ et la deuxi\`eme en 
$\oplus_{c'\in [\vert a-b\vert+1+T_{\rho,a,b},a+b-1+T_{\rho,a,b}]_{2}}\rho\otimes [c']$. Ainsi le support cuspidal de la premi\`ere repr\'esentation s'obtient \`a partir de celui de la deuxi\`eme repr\'esentation en enlevant $\cup_{\ell \in [1,T_{\rho,a,b}}\cup _{x\in [\vert a-b\vert/2+\ell, (a+b)/2-1+\ell]}\rho\vert\,\vert^{x},\rho\vert\,\vert^{-x}$; pour obtenir l'assertion annonc\'e on utilise le fait que $\zeta_{\rho,a,b}$ est le signe de $a-b$ si $a-b\neq 0$ et on a donc pour tout $\ell$ comme ci-dessus:
$$
[\vert a-b\vert/2+\ell, (a+b)/2-1+\ell] \cup -[\vert a-b\vert/2+\ell, (a+b)/2-1+\ell]=$$
$$[(a-b)/2+\zeta_{\rho,a,b}\ell, \zeta_{\rho,a,b}((a+b)/2-1+\ell)]\cup -[(a-b)/2+\zeta_{\rho,a,b}\ell, \zeta_{\rho,a,b}((a+b)/2-1+\ell)].
$$
Cela termine la preuve de la proposition.
\subsection{Intersection de 2 paquets d'Arthur\label{instersection}}

La proposition pr\'ec\'edente a comme corollaire imm\'ediat:

\begin{cor} Soient $\psi,\psi'$ des morphismes de $W_{F}\times SL(2,{\mathbb C})\times SL(2,{\mathbb C})$ dans $GL(m_{G}^*,{\mathbb C})$ se factorisant par $G^*$. On suppose que $\Pi(\psi)\cap \Pi(\psi')\neq \emptyset$ alors les restrictions de $\psi$ et de $\psi'$ \`a $W_{F}$ fois la diagonale de $SL(2,{\mathbb C})$ sont conjugu\'ees.
\end{cor}

En effet le support cuspidal \'etendu de $\pi$ fixe la restriction de $\psi$ et de $\psi'$ \`a $W_{F}$ fois la diagonal de $SL(2,{\mathbb C})\times SL(2,{\mathbb C})$.

\

En g\'en\'eral, il semble peu clair de trouver des conditions suffisantes pour que $\Pi(\psi)\cap \Pi(\psi')\neq \emptyset$.

\subsection{Exemples}
Dans cette section, on donne des exemples de repr\'esentations $\pi\in \Pi(\psi)$ telles qu'il existe un autre morphisme $\psi'$ avec $\pi\in \Pi(\psi')$ et plus pr\'ecis\'ement $\pi$ dans le paquet de Langlands associ\'e \`a $\psi'$.  Supposons d'abord que $\psi$ soit de restriction discr\`ete \`a la diagonale et que pour tout $(\rho,a,b)\in Jord(\psi)$, $a\geq b$. On fixe $\pi\in \Pi(\psi)$ et on note $\underline{t}$ et $\underline{\eta}$ le param\`etre de $\psi$. 
\begin{prop}Sous les hypoth\`eses ci-dessus, il existe un morphisme $\psi'$ tel que $\pi\in \Pi(\psi'_{L})$ si et seulement si pour tout $(\rho,a,b)\in Jord(\psi)$, $\underline{t}(\rho,a,b)\in \{0, [b/2]\}$. De plus, si ceci est r\'ealis\'e, $\psi'$ se d\'eduit de $\psi$ en rempla\c{c}ant tous les $(\rho,a,b)\in Jord(\psi)$ tel que $\underline{t}(\rho,a,b)=0$ par $\cup_{c\in [a-b+1,a+b-1]; c\equiv{a+b-1}[2]}(\rho,c,1)$.
\end{prop}
On fixe $\pi\in \Pi(\psi)$ et donc ses param\`etres $\underline{t}$ et $\underline{\eta}$. Ceux-ci permettent de d\'eterminer les param\`etres de Langlands de $\pi$ et plus pr\'ecis\'ement:
$$
\pi\hookrightarrow \times_{(\rho,a,b)\in Jord(\psi); \underline{t}(\rho,a,b)>0} J(St(\rho,a),-(b-1)/2, -(b-1)/2+\underline{t}(\rho,a,b)) \times \pi_{temp},\eqno(1)
$$
o\`u $\pi_{temp}$ est une repr\'esentation temp\'er\'ee (ici d'ailleurs une s\'erie discr\`ete) dans le paquet associ\'e au morphisme $\psi_{temp}$ tel que $Jord(\psi_{temp})=\cup_{(\rho,a,b)\in Jord(\psi)}\cup_{c\in [a-b+1+2\underline{t}(\rho,a,b),a+b-1-2\underline{t}(\rho,a,b)]; (-1)^{c}=(-1)^{a+b-1}}(\rho,c,1)$.

Soit $\psi'$ un morphisme et $\pi'\in \Pi(\psi'_{L})$; on conna\^{\i}t aussi les param\`etres de Langlands de $\pi'$:
$$
\pi'\hookrightarrow \times_{(\rho',a',b')\in Jord(\psi'); b'>1}J(St(\rho',a'),-(b'-1)/2),-\delta_{b'})\times \pi'_{temp}\eqno(2)
$$
o\`u pour tout $(\rho',a',b')\in Jord(\psi')$ avec $b'>1$, $\delta_{b'}=1/2$ si $b'$ est pair et $1$ si $b'$ est impair et o\`u $\pi'_{temp}$ est une repr\'esentation temp\'er\'ee dans le paquet associ\'e au morphisme $\psi'_{temp}$ dont les blocs de Jordan sont $\cup_{(\rho,a',b')\in Jord(\psi'); b'\equiv 1[2]}(\rho,a',1)$.
Il est clair que si $\underline{t}(\rho,a,b)\in \{0,[b/2]\}$ pour tout $(\rho,a,b)\in Jord(\psi)$, le morphisme $\psi'$ d\'ecrit dans l'\'enonc\'e est tel que (2) soit v\'erifi\'e pour $\pi=\pi'$ et pour $\pi'_{temp}=\pi_{temp}$. R\'eciproquement supposons qu'il existe $\psi'$ tel que $\pi\in \Pi(\psi'_{L})$; l'unicit\'e des param\`etres de Langlands assurent que pour tout $(\rho,a,b)\in Jord(\psi)$ tel que $\underline{t}(\rho,a,b)>0$, on a $(b-1)/2-\underline{t}(\rho,a,b)=\delta_{b}$ ce qui veut dire que $\underline{t}(\rho,a,b)=[b/2]$. Et on a bien la condition de l'\'enonc\'e. Cela termine la preuve.

\subsection{Unicit\'e des repr\'esentations non ramifi\'ees dans un paquet non ramifi\'e\label{unicite}}
Ici on suppose que $\psi$ est trivial sur la premi\`ere copie de $SL(2,{\mathbb C})$ et que la restriction de $\psi$ \`a $W_{F}$ est non ramifi\'ee; c'est ce que Clozel dans \cite{clozelimrn} appelle un param\`etre d'Arthur non ramifi\'e. Le r\'esultat ci-dessous est assez ''\'evident'', on ne l'inclut que par souci de compl\'etude.

\begin{prop} Dans un paquet d'Arthur associ\'e \`a un param\`etre non ramifi\'e, il y a exactement une repr\'esentation non ramifi\'ee.
\end{prop}
On fixe donc $\psi$ non ramifi\'e et $\pi\in \Pi(\psi)$. On suppose que $\pi$ est non ramifi\'e; on conna\^{\i}t alors le support cuspidal de $\pi$ gr\^ace \`a \ref{definitionsupportcuspidal}. En particulier avec un support cuspidal fix\'e, il y a au plus une repr\'esentation non ramifi\'e par la th\'eorie g\'en\'erale. Ainsi $\Pi(\psi)$ a au plus une repr\'esentation non ramifi\'ee. R\'eciproquement comme $\psi$ est non ramifi\'e, le paquet de Langlands \`a l'int\'erieur de $\Pi(\psi)$ contient exactement une repr\'esentation d\'ecrite dans \ref{description} et cette repr\'esentation est non ramifi\'ee. D'o\`u la proposition puisque c'est le quotient de Langlands d'une induite non ramifi\'ee.

\section{Propri\'et\'es des modules de Jacquet}
\subsection{Preuve de \ref{lemmesurjacquet}\label{proprietesimpledejac}}
On fixe $\psi$ et $\pi\in \Pi(\psi)$. On fixe aussi $\mathcal{E}$ un ensemble de demi-entiers rang\'es par ordre d\'ecroissant. On va montrer \ref{lemmesurjacquet}, c'est-\`a-dire que 
$$Jac_{x\in \mathcal{E}}\pi\neq 0 \Rightarrow Jac_{x\in \mathcal{E}}^\theta\pi^{GL}(\psi)\neq 0.\eqno(1)$$
On montre cette propri\'et\'e par r\'ecurrence d'abord sur $\vert \mathcal{E}\vert$ puis sur $\ell_{\psi}:=\sum_{(\rho,a,b)\in Jord(\psi)}(inf(a,b)-1)$.  

D'abord on se ram\`ene au cas o\`u $Jord(\psi)$ est sans multiplicit\'e et de bonne parit\'e: en effet soit $(\rho,a,b)\in Jord(\psi)$ et supposons que cet \'el\'ement soit y apparaisse avec multiplicit\'e soit ne soit pas de bonne parit\'e. On note $\psi'$ le morphisme qui se d\'eduit de $Jord(\psi)$, en enlevant $(\rho,a,b)$ et $\theta^*(\rho,a,b)$ ($\theta^*$ est le dual de l'automorphisme, $\theta$, qui sert pour l'endoscopie tordue, en g\'en\'eral la contragr\'ediente). On sait (cf. \cite{holomorphie} 2.6) qu'il existe $\pi'\in \Pi(\psi')$ et une inclusion
$$
\pi\hookrightarrow Speh(St(\rho,a),b)\times \pi'.\eqno(2)
$$De plus $\pi^{GL}(\psi)\simeq Speh(St(\rho,a),b)\times \theta(Speh(St(\rho,a),b))\times \pi^{GL}(\psi')$.
On applique les formules standard des modules de Jacquet \`a (2): il existe une d\'ecomposition de $\mathcal{E}$ en 3 sous-ensembles (dont certains peuvent \^etre vides), $\mathcal{E}=\cup_{i\in [1,3]}\mathcal{E}_{i}$ totalement ordonn\'e par l'ordre induit tels que
$$
Jac_{x\in \mathcal{E}_{1}}Speh(St(\rho,a),b)\neq 0, Jac_{x\in \mathcal{E}_{2}}\theta^*(Speh(St(\rho,a),b))\neq 0, Jac_{x\in \mathcal{E}_{3}}\pi'\neq 0,\eqno(3)
$$
o\`u les 2 premiers calculs de modules de Jacquet se font dans le groupe lin\'eaire convenable.  Il y a une condition suppl\'ementaire aux non nullit\'e ci-dessus: si $\rho\not\simeq \theta^*(\rho)$, $\mathcal{E}_{2}=\emptyset$ et il n'y a rien \`a ajouter. Sinon, il faut aussi que $\mathcal{E}_{1}\cup -\mathcal{E}_{2}\subset \cup_{\ell \in [(a-b)/2,(a+b)/2-1]}[\ell, \ell-a]$. On conclut donc facilement que l'implication (1) pour $\pi$ r\'esulte de son analogue avec $\mathcal{E}$ remplac\'e par $\mathcal{E}_{3}$ pour $\pi'$ mais on va d\'etailler en utilisant les simplifications d\^ues au fait que  $\mathcal{E}$ est rang\'e par ordre d\'ecroissant; on suppose que $\rho\simeq \theta^*\rho$ pour qu'il y ait quelque chose de non trivial \`a montrer. La non nullit\'e des 2 premiers modules de Jacquet dans (3) assure que $\mathcal{E}_{1}$ est un sous-intervalle de $[(a-b)/2, -(a+b)/2+1]$ et que $\mathcal{E}_{2}$ est aussi un sous-intervalle de $[(a-b)/2, -(a+b)/2+1]$ et la condition suppl\'ementaire ne joue que si $b=1$ et dit que $\mathcal{E}_{1}\cup -^t\mathcal{E}_{2}\subset [(a-1)/2,-(a-1)/2]$, o\`u $^t$ est l'inversion de l'ordre. On a donc $Jac^g_{x\in \mathcal{E}_{1}}Jac^d_{x\in -^t\mathcal{E}_{2}}Speh(St(\rho,a),b)\neq 0$ et de m\^eme $Jac^g_{x\in \mathcal{E}_{2}}Jac^d_{x\in -^t\mathcal{E}_{1}}Speh(St(\rho,a),b)\neq 0$. D'o\`u certainement $$Jac^\theta_{x\in \mathcal{E}_{1}\cup\mathcal{E}_{2}}\biggl(Speh(St(\rho,a),b)\times Speh(St(\rho,a),b)\biggr)\neq 0$$
et si on sait que $Jac^\theta_{x\in \mathcal{E}_{3}}\pi^{GL}(\psi')\neq 0$ on aura certainement (1). On est donc ramen\'e au cas o\`u $Jord(\psi)$ a bonne parit\'e et n'a pas de multiplicit\'e.

Soit $x_{0}$ l'\'el\'ement maximal de $\mathcal{E}$.

On fixe un bon ordre total sur $Jord(\psi)$  tel que $(\rho,a,b)< (\rho,a',b')$ si $(a-b)/2 \leq x_{0}$ et $(a'-b')/2>sup(-1/2,x_{0})$. On fixe $\psi_{>}$ un morphisme dominant $\psi$, tel que

(a)$Jord(\psi_{>})$ contient tous les \'el\'ements de $Jord(\psi)$ de la forme $(\rho,a,b)$ avec $(a-b)/2\leq x_{0}$

(b) pour tout $(\rho,a_{>},b_{>})\in Jord(\psi_{>})$ v\'erifiant $(a_{>}-b_{>})/2>x_{0}$ on a m\^eme $>>x_{0}$.

On sait qu'il existe $\pi_{>}$ tel que $\pi$ s'obtient \`a partir de $\psi_{>}$ en prenant des modules de Jacquet de la forme $Jac_{y\in \mathcal{Y}}$; comme $Jord_{\psi_{>}}$ v\'erifie (a), on sait que les \'el\'ements de $\mathcal{Y}$ sont tous strictement sup\'erieur \`a $x_{0}+1$. On a donc aussi
$$
Jac_{x\in \mathcal{E}}\pi=Jac_{x\in \mathcal{E}}Jac_{y\in \mathcal{Y}}\pi_{>}=Jac_{y\in \mathcal{Y}}Jac_{x\in \mathcal{E}}\pi_{>}.
$$En particulier si $Jac_{x\in \mathcal{E}}\pi\neq 0$, n\'ecessairement $Jac_{x\in \mathcal{E}}\pi_{>}\neq 0$.  On pose $$\pi^{GL}_{\leq }:=\times_{(\rho,a',b')\in Jord(\psi); (a'-b')/2\leq x_{0}}Speh(St(\rho,a),b)$$ et on a $\pi^{GL}(\psi)=\pi^{GL}_{\leq}\times \sigma$ o\`u $\sigma$ est explicite et v\'erifie $Jac_{x\in \mathcal{E}'}\sigma=0$ pour tout sous-ensemble non vide $\mathcal{E}'$ de $\mathcal{E}$. De plus $\pi^{GL}(\psi_{>})$ est une induite analogue $\pi^{GL}_{\leq }\times \sigma_{>}$. Ainsi si $Jac_{x\in \mathcal{E}}^\theta\pi ^{GL}(\psi_{>})\neq 0$, n\'ecessairement $Jac_{x\in \mathcal{E}}^\theta \pi^{GL}_{\leq}\neq 0$ et $Jac^\theta_{x\in \mathcal{E}}\pi^{GL}(\psi)\neq 0$. Il suffit donc de montrer (1) pour $\pi_{>}$. On suppose donc que $\pi=\pi_{>}$, c'est-\`a-dire que pour tout $(\rho,a,b)\in Jord(\psi)$ si $(a-b)/2>x_{0}$, on a $(a-b)/2>>x_{0}$.

Puisqu'en particulier, $Jac_{x_{0}}\pi\neq 0$, on sait (cf. \cite{holomorphie} 2.7), qu'il existe $(\rho,a,b)\in Jord(\psi)$ avec $x_{0}=(a-b)/2$. Si $sup(a,b)>1$, on a alors certainement $Jac^\theta_{x_{0}}Speh(St(\rho,a),b)\neq 0$ et $Jac^\theta(\pi^{GL})(\psi)\neq 0$. Supposons que pour tout tel choix de $(\rho,a,b)$, on a $a=b=1$, d'o\`u n\'ecessairement $x_{0}=0$; on a suppos\'e que $Jord(\psi)$ n'a pas de multiplicit\'e et on v\'erifie que l'on ne peut alors avoir $Jac_{x_{0}}\pi\neq 0$; en effet, on peut consid\'erer que $(\rho,1,1)$ est le plus petit \'el\'ement de $Jord(\psi)$. Par les hypoth\`eses, on sait que pour tout $(\rho,a',b')\in Jord(\psi)$ diff\'erent de cet \'el\'ement, $\vert (a'-b')/2\vert\geq 1$. On obtient donc $\pi$ comme module de Jacquet \`a partir d'une repr\'esentation $\pi_{>>}$, en ayant $\pi=Jac_{y\in \mathcal{Y}}\pi_{>>}$; tous les \'el\'ements de $\mathcal{Y}$ sont de valeurs absolue au moins 2. On v\'erifie sur la d\'efinition de $\pi_{>>}$ que $Jac_{0}\pi_{>>}=0$ et comme $Jac_{0}$ commute \`a $Jac_{y\in \mathcal{Y}}$, on a aussi $Jac_{0}\pi=0$. Ainsi on vient de d\'emontrer le cas o\`u $\mathcal{E}$ est r\'eduit \`a $x_{0}$, ce qui initie la r\'ecurrence sur $\mathcal{E}$.

On suppose donc que $\mathcal{E}$ n'est pas r\'eduit \`a $x_{0}$ et on consid\`ere d'abord le cas o\`u $x_{0}\geq 0$.
On va utiliser aussi la r\'ecurrence sur $\ell(\psi)$ que l'on va initialiser.

Supposons pour commencer que pour tout $(\rho,a,b)\in Jord(\psi)$ si $x_{0}=(a-b)/2$ alors $inf(a,b)=1$. Cela se produit par exemple si $\ell(\psi)=0$. On a vu ci-dessus, que n\'ecessairement $sup(a,b)>1$. Avec l'hypoth\`ese que $Jord(\psi)$ n'a pas de multiplicit\'e, on sait donc qu'il existe un unique $(\rho,a,b)$ v\'erifiant la propri\'et\'e $x_{0}=(a-b)/2$; on rappelle qu'ici $x_{0}\geq 0$, d'o\`u $a>b=1$. On note $\psi'$ le morphisme qui se d\'eduit de $\psi$ en rempla\c{c}ant $(\rho,a,b)$ par $(\rho,a-2,b)$ ou en enlevant $(\rho,a,b)$ si $a=2$. On v\'erifie que $Jac_{x_{0}}\pi\in \Pi(\psi')$ et que $Jac^\theta_{x_{0}}\pi^{GL}(\psi)=\psi'$; pour la premi\`ere assertion, on applique la d\'efinition des \'el\'ements de $\Pi(\psi')$, c'est exactement comme cela qu'on les obtient. La deuxi\`eme assertion est un calcul simple (qui en fait est utilis\'e dans la d\'efinition de $\Pi(\psi')$):
$$
\pi^{GL}(\psi)=\times_{(\rho',a',b')\neq (\rho,a,b)}Speh(St(\rho',a'),b')\times St(\rho,a).
$$Or
$
Jac^g_{x_{0}}Speh(St(\rho',a'),b')=0=Jac^d_{-x_{0}}Speh(\rho',a',b')$ pour tout $(\rho,a',b')\neq (\rho,a,1)$ car $(a'-b')/2\neq x_{0}$, et $Jac^\theta_{(a-1)/2)}St(\rho,a)=St(\rho,a-2)$, d'o\`u a fortiori $$Jac^\theta_{x_{0}}\pi^{GL}(\psi)=\times_{(\rho',a',b')\neq (\rho,a,b)}Speh(St(\rho',a'),b')\times St(\rho,a-2)=\pi^{GL}(\psi').$$
Il suffit maitenant d'appliquer le r\'esultat \`a $\psi'$ et $\mathcal{E}-\{x_{0}\}$ ce qui est loisible par r\'ecurrence sur $\vert \mathcal{E}\vert$ et on obtient le r\'esultat. D'o\`u en particulier l'initialisation de la r\'ecurrence sur $\ell(\psi)$.

On suppose maintenant qu'il existe $(\rho,a,b)\in Jord(\psi)$ tel que $(a-b)/2=x_{0}$ et $b>1$; on fixe un tel \'el\'ement. On va alors diminuer $\ell(\psi)$.  On fixe un bon ordre sur $Jord(\psi)$ de sorte que $(\rho,a,b)$ soit l'\'el\'ement maximal parmi les $(\rho,a',b')\in Jord(\psi)$ v\'erifiant $(a'-b')/2\leq x_{0}$. Ceci est possible. On note $\underline{t}$ et $\underline{\eta}$ les param\`etres permettant de d\'efinir $\psi$. On fixe $T$ un entier ''grand'' et on note  $\psi_{T}$ le morphisme qui s'obtient \`a partir de $\psi$ en rempla\c{c}ant $(\rho,a,b)$ par $(\rho,a+2T,b)$ pour $T$ grand. On rappelle que l'on a fait l'hypoth\`ese que pour tout $(\rho,a',b')\in Jord(\psi)$ avec $(a'-b')/2>x_{0}$, on a $a'>>b'$; on peut donc supposer que $\psi_{T}$ domine $\psi$ pour l'ordre fix\'e. Ainsi $\Pi(\psi)$ s'obtient \`a partir de $\Pi(\psi_{T})$ en prenant les modules de Jacquet qui font descendre $(\rho,a+2T,b)$ vers $(\rho,a,b)$. Pr\'ecis\'ement, on sait qu'il existe $\pi_{T}\in \Pi(\psi_{T})$ tel que 
$$
\pi=\circ_{\ell\in [1,T]}Jac_{(a-b)/2+T-\ell+1, \cdots, (a+b)/2+T-\ell}\pi_{T}.
$$
On pose $t:=\underline{t}(\rho,a,b)$ et on note $\psi'_{T}$ le morphisme qui se d\'eduit de $\psi_{T}$ en rempla\c{c}ant $(\rho,a+2T,b)$ par $\cup_{c\in [a-b+1+t, a+b-1-t],(-1)^{c}=(-1)^{a+b-1}}(\rho,c+2T,1)$. On sait qu'il existe une repr\'esentation $\pi'_{T}\in \Pi(\psi'_{T})$ et une inclusion
$$
\pi_{T}\hookrightarrow <\times_{j\in [1,t]}St(\rho,a+2T)\vert\,\vert^{(b-1)/2+j-1}> \times \pi'_{T}.\eqno(4)
$$Les crochets indiquent que l'on prend l'unique sous-module irr\'eductible pour l'induite \'ecrite. Pour tout $T'\in [0,T[$, on d\'efinit par r\'ecurrence descendante $\pi_{T'}:=Jac_{(a-b)/2+T'+1, \cdots, (a+b)/2+T' }\pi_{T'+1}$. On sait que le r\'esultat est une repr\'esentation irr\'eductible dans $\Pi(\psi_{T'})$ analogue de $\psi_{T}$ quand $T=T'$. On d\'efinit, $\pi'_{T'}$ de fa\c{c}on analogue en appliquant $Jac_{(a-b)/2+t+T'+1, \cdots, (a+b)/2-t+T'}\pi'_{T'+1}$. On sait aussi que si cette repr\'esentation est non nulle, c'est un \'el\'ement de $\Pi(\psi'_{T'})$ obtenu en rempla\c{c}ant $T$ par $T'$ dans la d\'efinition ci-dessus. On montre par r\'ecurrence descendante que
$$
\pi_{T'}\hookrightarrow <\times_{j\in [1,t]}St(\rho,a+2T')\vert\,\vert^{-(b-1)/2+j-1}>\times \pi'_{T'};
$$
il suffit d'oublier que dans (4), $T$ est grand et de calculer, $Jac_{(a-b)/2+T, \cdots, (a+b)/2-1+T}$ aux 2 membres de (4). En fait on calcule d'abord $Jac_{(a-b)/2+T, \cdots, (a-b)/2+T+t}$ du membre de droite; pour tout $x\in [(a-b)/2+T,(a-b)/2+T+t-1]$, on certainement $Jac_{x}\pi'_{T}=0$ car pour tout $(\rho,a'',b'')\in Jord(\psi'_{T})$, ou bien $(a''-b'')/2<x_{0}$ ou bient $(a''-b'')/2>>x_{0}$ ou bien $(a''-b'')/2\in [(a-b)/2+t,(a+b)/2-1-t]$ et donc il n'existe aucun \'el\'ement de $Jord(\psi'_{T})$, $(\rho,a'',b'')$, tel que $(a''-b'')/2=x$. On \'ecrit $<\times_{j\in [1,t]}St(\rho,a+2T')\vert\,\vert^{(b-1)/2+j-1}>$ comme la repr\'esentation attach\'ee \`a $\rho$ et aux multisegments:
$$
\begin{matrix}
(a-b)/2+T &\cdots &-(a+b)/2+1-T\\
\vdots &\vdots &\vdots\\
(a-b)/2+t-1+T &\cdots &-(a+b)/2+t-T
\end{matrix}
$$
Quand on applique $Jac_{(a-b)/2+T, \cdots, (a-b)/2+T+t}$, on enl\`eve la premi\`ere colonne. Ensuite on applique $Jac_{(a-b)/2+t+T, \cdots, (a+b)/2-1-t+T}$ qui lui n'agit que sur $\pi'_{T}$ pour donner par d\'efinition $\pi'_{T-1}$. Ensuite on applique $Jac_{(a+b)/2-t+T, \cdots, (a+b)/2+T}$; on v\'erifie encore que pour tout $x\in [(a+b)/2-t+T, (a+b)/2+T]$, $Jac_{x}\pi'_{T-1}=0$ et le r\'esultat est donc d'enlever la derni\`ere colonne de la matrice ci-dessus. En mettant ces 3 calculs ensemble on voit qu'appliquer le module de Jacquet au membre de droite de (4) revient \`a remplacer $T$ par $T-1$. Progressivement, on obtient
$$
\pi\hookrightarrow <\times_{j\in [1,t]}St(\rho,a)\vert\,\vert^{-(b-1)/2+j-1}>\times \pi'_{0},\eqno(5)
$$
o\`u $\pi'_{0}\in \Pi(\psi'_{0})$. On applique $Jac_{x\in \mathcal{E}}$ aux 2 membres de (5). Il faut diff\'erentier suivant que $t=0$ ou non. Dans le premier cas, $\pi=\pi'_{0}$; on a alors
$$
\pi^{GL}(\psi'_{0})=\times_{c\in [(a-b)+1,(a+b)-1], (-1)^c=(-1)^{a+b-1}}St(\rho,c)\times_{(\rho',a',b')\neq (\rho,a,b)}Speh (St(\rho',a'),b').
$$
$$
Jac^\theta_{x\in \mathcal{E}}\pi^{GL}(\psi'_{0})=Jac_{x\in \mathcal{E}_{g}}^gJac_{-x\in \mathcal{E}_{d}}^dSt(\rho,a-b+1)\times_{c\in ](a-b)+2,(a+b)-1],(-1)^c=(-1)^{a+b-1}}St(\rho,c)
$$
$$\times Jac_{x\in \mathcal{E}-\mathcal{E}_{g}-\mathcal{E}_{d}}\times_{(\rho',a',b')\neq (\rho,a,b)}Speh (St(\rho',a'),b'),\eqno(6)
$$o\`u $\mathcal{E}_{g}$ est soit vide soit r\'eduit \`a $x_{0}$ et $\mathcal{E}_{d}$ est soit vide soit r\'eduit \`a $x_{0}$. De fa\c{c}on analogue, on a
$$
\pi^{GL}(\psi)=Speh(St(\rho,a),b)\times _{(\rho',a',b')\neq (\rho,a,b)}Speh (St(\rho',a'),b').
$$
Quand on calcule $Jac^\theta_{x\in \mathcal{E}}\pi^{GL}(\psi)$ on obtient des termes index\'es par 2 sous-segments, de $[(a-b)/2,-(a+b)/2]$ chacun \'etant vide ou commen\c{c}ant par $(a-b)/2$ que l'on note $\mathcal{E}_{g}$ et $\mathcal{E}_{d}$ par analogie au calcul pr\'ec\'edent et le r\'esultat est comme ci-dessus. A fortiori si (6) est non nul, $Jac^\theta_{x\in \mathcal{E}}\pi^{GL}(\psi)\neq 0$. Or on conna\^{\i}t le r\'esultat pour $\pi'_{0}$ puisqu'ici on ne change pas $\mathcal{E}$ mais que l'on baisse $\ell(\psi)$ et que l'on a d\'ej\`a amorc\'e la r\'ecurrence pour $x_{0}\geq 0$.

On traite maintenant le cas o\`u $t>0$.
On utilise tout de suite le fait que les \'el\'ements de $\mathcal{E}$ sont rang\'es par ordre d\'ecroissant; le module de Jacquet du membre de droite de (5) a une filtration dont les termes du gradu\'e associ\'e sont index\'es par les  sous-segments $\mathcal{E}_{1}$ de $[(a-b)/2, -(a+b)/2+1]$ soit vide soit commen\c{c}ant  par $(a-b)/2$ et le terme correspondant est
$$
Jac_{x\in \mathcal{E}_{1}}<\times_{j\in [1,t]}St(\rho,a)\vert\,\vert^{-(b-1)/2+t-1}>\times Jac_{x\in \mathcal{E}_{2}}\pi'_{0},
$$
o\`u $\mathcal{E}_{2}$ est le compl\'ementaire de $\mathcal{E}_{1}$ dans $\mathcal{E}$. On peut appliquer le r\'esultat par r\'ecurrence \`a $\pi'_{0}$ et $\mathcal{E}_{2}$ car soit $\mathcal{E}_{2}$ est strictement plus petit que $\mathcal{E}$ et on utilise la r\'ecurrence sur le cardinal de $\mathcal{
E}$ soit on a \'egalit\'e et on utilise la r\'ecurrence sur $\ell(\psi)$.

Pour conclure, il suffit de remarquer que $$Jac^\theta_{x\in \mathcal{E}_{2}}\pi^{GL}(\psi'_{0})\neq 0\Rightarrow
Jac^{\theta}_{x\in \mathcal{E}_{2}}\times_{(\rho',a',b')\neq (\rho,a,b)}Speh(St(\rho',a'),b')\neq 0$$et que $Jac^\theta_{x\in \mathcal{E}}\pi^{GL}(\psi)$ contient $Jac^\theta_{x\in \mathcal{E}_{1}}Sp(St(\rho,a),b)\times Jac^{\theta}_{x\in \mathcal{E}_{2}}\times_{(\rho',a',b')\neq (\rho,a,b)}Speh(St(\rho',a'),b')$ qui est alors non nul.

\

Il reste le cas o\`u $x_{0}<0$. Ici on fait une r\'ecurrence sur $\sum_{(\rho,a,b)\in Jord(\psi); a<b}b$. Si ce nombre vaut 0,  on a $Jac_{x_{0}}\pi=0$ et il n'y a rien \`a d\'emontrer. On a d\'ej\`a r\'eduit au cas o\`u $Jord(\psi)$ est de bonne parit\'e,  sans multiplicit\'e et o\`u pour tout $(\rho,a,b)\in Jord(\psi)$ si $a\geq b$ alors $a>>b$ cela ne change pas le nombre sur lequel on fait une r\'ecurrence. On fixe $(\rho,a,b)$ tel que $\vert(a-b)/2\vert$ soit minimal, d'o\`u certainement $(a-b)/2\geq x_{0}$ pour qu'il y ait quelque chose \`a d\'emontrer. Le choix n'est pas unique mais on en fait un. On fixe un bon ordre sur $Jord(\psi)$ tel que $(\rho,a,b)$ soit l'\'el\'ement minimal; ceci est possible. D'abord on montre que l'on peut supposer que $a=b$ ou $a=b-1$. En effet, s'il n'en est pas ainsi, on note $\psi'$ le morphisme qui se d\'eduit de $\psi$ en changeant $(\rho,a,b)$ par $(\rho,a,b-2)$ et on v\'erifie qu'il existe $\pi'\in \Pi(\psi')$ avec une inclusion
$$
\pi\hookrightarrow St(\rho,a)\vert\,\vert^{-(b-1)/2}\times \pi';\eqno(7)
$$
en effet ceci est montr\'e en \cite{paquetdiscret} dans le cas ou $\psi$ est de restriction discr\`ete \`a la diagonale et cela s'applique donc ici pour $\psi_{>>}$ dominant $\psi$ et donc $\psi'_{>>}$ dominant $\psi'$; d'o\`u
$$
\pi_{>>}\hookrightarrow St(\rho,a)\vert\,\vert^{-(b-1)/2}\times \pi'_{>>};
$$
Pour obtenir $\pi$ \`a partir de $\pi_{>>}$, il faut appliquer des $Jac_{y, y\in {\mathcal Y}}$ mais avec nos hypoth\`eses sur $\psi$, on est s\^ur que ${\mathcal Y}$ est form\'e d'\'el\'ements n\'egatifs tous strictement inf\'erieur \`a $(a-b)/2$; ainsi $Jac_{y,y\in {\mathcal Y}}$ appliqu\'e au terme de droite commute \`a l'induction par $St(\rho,a)\vert\,\vert^{-(b-1)/2}$ et on obtient (7).

On peut appliquer l'hypoth\`ese de r\'ecurrence \`a $\pi'$ pour $\mathcal{E}$ ou pour tout sous-ensemble de $\mathcal{E}$. Soit $\mathcal{F}$ un sous-ensemble de $\mathcal{E}$ h\'eritant de l'ordre de $\mathcal{E}$. On v\'erifie encore que $Jac^\theta_{x\in \mathcal{F}}\pi^{GL}(\psi')\neq 0$ n\'ecessite que $Jac_{x\in \mathcal{F}}^{\theta}\times _{(\rho',a',b')\neq (\rho,a,b)}Speh(St(\rho',a'),b')\neq 0$ puisque $(\rho,a,b)$ n'est plus dans $Jord(\psi')$ et qu'il est remplac\'e par $(\rho,a,b-2)$ avec certainement $(a-b+2)/2=(a-b)/2+1>x_{0}$. On conclut ensuite exactement comme dans la fin de la preuve du cas $x_{0}\geq 0$. Le cas $a=b$, a d\'ej\`a \'et\'e trait\'e puisqu'on la ramen\'e \`a $a>>b$. Il faut voir le cas o\`u $a=b-1$. On note $\underline{t}$ et $\underline{\eta}$ les param\`etres de $\pi$ pour l'ordre choisi. Et on pose $t:=\underline{t}(\rho,a,b)$ et $\eta:=\underline{\eta}(\rho,a,b)$.

On termine la preuve en admettant momentan\'ement le lemme ci-dessous: en effet dans le cas (i) de ce lemme, on applique l'hypoth\`ese de r\'ecurrence \`a $\psi'$ et on obtient facilement le r\'esultat pour $\psi$. Dans le cas (ii), on raisonne comme ci-dessus.

\begin{lem} (i)
On suppose  que $t=0$ et que ${\eta}=-$. Ce cas est tr\`es simple: on note $\psi'$ le morphisme qui se d\'eduit de $\psi$ en changeant $(\rho,a,a+1)$ en $(\rho,a+1,a)$ alors $\pi\in \Pi(\psi')$.

(ii) On suppose que $t\neq 0$ ou que ${\eta}\neq -$. On note $\psi'$ le morphisme qui se d\'eduit de $\psi$ en rempla\c{c}ant $(\rho,a,a+1)$ par $(\rho,a,a-1)$. Alors, il existe $\pi'\in \Pi(\psi')$ et une inclusion
$$
\pi\hookrightarrow St(\rho,a)\vert\,\vert^{-a/2}\times\pi'.
$$
\end{lem}
Le (i) est vrai par d\'efinition si $\psi$ est de restriction discr\`ete \`a la diagonale; on l'applique \`a $\psi_{>>}$ dominant $\psi$. Et (i)  reste vrai quand on a appliqu\'e les modules de Jacquet pour redescendre de $\psi_{>>}$ \`a $\psi$ puisque l'on ne modifie pas le plus petit bloc de Jordan qui est $(\rho,a,a+1)$ par notre choix.

On suppose encore que $t=0$ mais que $\eta=+$. On fixe encore $\psi_{>>}$ dominant $\psi$ et contenant $(\rho,a,a+1)$ comme plus petit \'el\'ement. Soit $\pi_{>>}\in \Pi(\psi_{>>})$ d\'etermin\'e par $\underline{t},\underline{\eta}$. Par construction, $\pi_{>>}\in \Pi(\psi'_{>>})$ ou $\psi'_{>>}$ se d\'eduit de $\psi_{>>}$ en rempla\c{c}ant $(\rho,a,a+1)$ par $(\rho,1,2j)$ pour $j\in [1,a]$. Et comme le signe alterne sur ces blocs en prenant la valeur $+$ sur $(\rho,1,2)$, on a construit $\pi_{>>}$ comme sous-module irr\'eductible de l'induite
$$
\pi_{>>}\hookrightarrow <-1/2, \cdots, -a/2>_{\rho}\times \pi''_{>>},
$$
o\`u $\pi''_{>>}\in \Pi(\psi'_{>>})$ avec $\psi''_{>>}$ qui se d\'eduit de $\psi'_{>>}$ en rempla\c{c}ant $(\rho,1,2j)$ par $(\rho,1,2j-2)$ pour tout $j\in ]1,a]$ ($(\rho,1,2)$ dispar\^{\i}t). Le signe est la restriction du signe pour $\psi'_{>>}$ de fa\c{c}on naturelle; il alterne sur ces blocs en commen\c{c}ant maintenant par $-$; on peut donc remplacer $(\rho,1,2j-2)$ par $(\rho,2j-2,1)$ (c'est dans les d\'efinitions). On note $\psi''$ le morphisme qui se d\'eduit de $\psi$ en rempla\c{c}ant $(\rho,a,a+1)$ par $\cup_{j\in [1,a-1]}(\rho,2j,1)$ et on montre qu'il existe $\pi''\in \Pi(\psi'')$ avec une inclusion
$$
\pi\hookrightarrow <-1/2, \cdots, -a+1/2>_{\rho}\times \pi''.
$$
Pour arriver \`a ce r\'esultat on applique les modules de Jacquet qui font passer de $\psi_{>>}$ \`a $\psi$ et ils ne touchent pas aux blocs $(\rho,2j,1)$ ni \`a la repr\'esentation $<-1/2,\cdots,-a/2>_{\rho}$ car ils sont de la forme $Jac_{y\in \mathcal{Y}}$ avec $y\leq -3/2$. Pour avoir exactement l'\'enonc\'e du lemme on regroupe les $\cup_{j\in [1,a-1]}(\rho,2j,1)$ en $(\rho,a,a-1)$ (ce qui est loisible); les param\`etres pour $\pi''$ sont $\underline{t}'', \underline{\eta}''$; ils co\"{\i}ncident avec $\underline{t},\underline{\eta}$ sur $Jord(\psi)-\{(\rho,a,a+1)\}$ et valent
$$
\underline{t}''(\rho,a,a-1)=0, \underline{\eta}''(\rho,a,a-1)=-.
$$
On ne suppose plus que $t=0$; comme ci-dessus, on peut supposer que $\psi$ est de restriction discr\`ete \`a la diagonale. On d\'efinit $\psi''$ en changeant dans $Jord(\psi)$, $(\rho,a,a+1)$ en $(\rho,a,a-1)$. On d\'efinit $\underline{t}'', \underline{\eta}''$ sur $Jord(\psi'')$ en demandant qu'ils co\"{\i}ncident avec $\underline{t}$ et $\underline{\eta}$ sur $Jord(\psi)-\{(\rho,a,a+1)\}$ et valent

si $\eta=+$, 
$\underline{t}''(\rho,a,a-1)=t-1$ et  $\underline{\eta}''(\rho,a,a-1)=-$;

si $\eta=-$, $\underline{t}''(\rho,a,a-1)=t$ et $\underline{\eta}''(\rho,a,a-1)=+$.

On note $\psi_{1}$ le morphisme qui se d\'eduit de $\psi$ en changeant $(\rho,a,a+1)$ en $(\rho,a-2,a+1)$ et $\underline{t}_{1}, \underline{\eta}_{1}$ qui se d\'eduisent naturellement de $\underline{t}$ et $\underline{\eta}$ avec comme seul changement que $\underline{t}_{1}(\rho,a-2,a+1)=t-1$. On note $\pi_{1}$ la repr\'esentation dans $\Pi(\psi_{1})$ correspondante et par d\'efinition
$$
\pi\hookrightarrow <-1/2, \cdots, a-1/2>_{\rho}\times \pi_{1}.
$$
On peut passer de $(\rho,a-2,a+1)$ \`a $(\rho,a-2,a-1)$ sans changer les param\`etres, ce qui donnent une repr\'esentation $\pi'_{1}$ et une inclusion
$$
\pi_{1}\hookrightarrow <-3/2, \cdots, -a+1/2>_{\rho}\times \pi'_{1}.
$$
On utilise l'inclusion de la repr\'esentation de Steinberg tordue $<-1/2, \cdots, a-1/2>_{\rho}$ dans l'induite $\rho\vert\,\vert^{-1/2}\times <1/2, \cdots, a-3/2>_{\rho}\times \rho\vert\,\vert^{a-1/2}$ et on compose avec l'inclusion ci-dessus. On obtient$$
\pi \hookrightarrow \rho\vert\,\vert^{-1/2}\times <-3/2, \cdots, -a+3/2>_{\rho}\times <1/2, \cdots, a-3/2>_{\rho}\times \rho\vert\,\vert^{a-1/2}\times \pi'_{1}.
$$
Mais $\rho\vert\,\vert^{a-1/2}\times \pi'_{1}$ est irr\'eductible  donc isomorphe \`a $\rho\vert\,\vert^{-a+1/2}\times \pi'_{1}$. D'o\`u une inclusion
$$
\pi\hookrightarrow \rho\vert\,\vert^{-1/2}\times <-3/1, \cdots, -a+3/2>_{\rho}\times \rho\vert\,\vert^{-a+1/2}\times <1/2, \cdots, a-3/2>_{\rho}\times \pi'_{1}.
$$
Comme $Jac_{x}\pi=0$ pour $x\in ]-1/2,-a+1/2]$, l'inclusion se factorise par
$$
\pi\hookrightarrow <-1/2, \cdots, -a+1/2>_{\rho}\times <1/2, \cdots, a-3/2>_{\rho}\times \pi'_{1}.
$$
On v\'erifie que l'induite $<1/2, \cdots, a-3/2>_{\rho}\times \pi'_{1}$ a un unique sous-module irr\'eductible, c'est la seule repr\'esentation dont le $Jac_{1/2, \cdots, a-3/2}\neq 0$. L'induite ci-dessus se factorise donc par ce sous-module (\`a cause de cette propri\'et\'e du module de Jacquet). On peut caract\'eriser ce sous-module irr\'eductible; c'est essentiellement la propri\'et\'e inverse de celle que l'on cherche \`a montrer en \'echangeant les r\^oles des 2 copies de $SL(2,{\mathbb C})$ et en rempla\c{c}ant $a$ par $a-1$. Nos constructions sont sym\'etriques en les 2 copies de $SL(2,{\mathbb C})$ on peut donc caract\'eriser le sous-module irr\'eductible comme \'etant la repr\'esentation dans $\Pi(\psi'')$ correspondant aux par\`etres donn\'ees. Cela termine la preuve.
\subsection{Un calcul approximatif de modules de Jacquet\label{approximationjac}}
La difficult\'e de la description de $\Pi(\psi)$ est qu'elle fait intervenir des calculs de modules de Jacquet qui sont difficile \`a expliciter. On peut en trouver des approximations, on veut en fait g\'en\'eraliser (5), (7) et le lemme de \ref{proprietesimpledejac}, ce qui sera r\'ealis\'e en \ref{blocpositif} et \ref{blocnegatif}. La situation typique dans laquelle on se trouvera est la suivante. On fixe $\psi$ et un morphisme $\tilde{\psi}$ tel que $Jord(\psi)$ et $Jord(\tilde{\psi})$ sont \'egaux \`a un \'el\'ement pr\`es qui vaut $(\rho,\alpha,1)$ dans $Jord(\psi)$ et $(\rho,\alpha+2T,1)$ dans $Jord(\tilde{\psi})$ avec $T>0$ et il n'existe pas $(\rho,a,b)\in Jord(\psi)$ avec $$\alpha<a-b+1\leq a+b-1<\alpha+2T.$$ En particulier on saura qu'il existe $\tilde{\pi}\in \Pi(\tilde{\psi})$ et $\pi\in \Pi(\psi)$ avec $\pi=Jac_{(\alpha-1)/2+T, \cdots, (\alpha+1)/2}\tilde{\pi}$. En g\'en\'eral $\tilde{\pi}$ sera mieux connu que $\pi$ comme sous-module irr\'eductible d'une certaine induite et on veut d\'eduire de cette connaissances des renseignements sur $\pi$, donc en fait montrer que $\pi$ est aussi sous-module irr\'eductible d'une certaine induite; les difficult\'es viennent de ce que l'on ne peut \'eliminer plusieurs choix pour l'induite et qu'il n'y a aucun moyen de d\'emontrer que l'induite en question a un unique sous-module irr\'eductible. 

On fixe donc $\psi$, $\pi\in \Pi(\psi)$, $(\rho,a,1)\in Jord(\psi)$, $\tilde{\psi}$, $\tilde{\pi}\in \Pi(\tilde{\psi})$, $T$  comme ci-dessus, en particulier $\pi=Jac_{(\alpha-1)/2+T, \cdots, (\alpha+1)/2}\tilde{\pi}$. On suppose qu'il existe $\tilde{\psi}'$ un autre morphisme et $\tilde{\pi}'\in \Pi(\tilde{\psi}')$ tel que
$$
\tilde{\pi}\hookrightarrow \times_{i\in [1,\ell]}<x_{i},y_{i}>_{\rho}\times \tilde{\pi}', \eqno(1)
$$
o\`u les $[x_{i},y_{i}]$ sont des segments d\'ecroissants avec $x_{i}\leq -y_{i}$, \'eventuellement $\ell=0$ auquel cas il n'y a pas d'induite mais le point alors est que $\tilde{\psi}'$ peut \^etre diff\'erent de $\tilde{\psi}$. On suppose que pour pour tout $i\in [1,\ell[$,  soit $x_{i}>(\alpha-1)/2+T$ soit $x_{i}<(\alpha-1)/2$ et $-y_{i}\notin [(\alpha+1)/2,(\alpha-1)/2+T]$; en particulier, puisque $x_{i}\leq -y_{i}$,  pour tout $x\in [(\alpha-1)/2+T,(\alpha+1)/2]$ l'induite $<x_{i},y_{i}>_{\rho}\times \rho\vert\,\vert^{x}$ est irr\'eductible.   On suppose que l'on a  aussi $x_{\ell}\leq -y_{\ell}$.

On aura en plus que pour tout $(\rho,a',b')\in Jord(\tilde{\psi}')$ tel que $(a'-b')/2\in [(\alpha-1)/2+T,(\alpha+1)/2]$, $b'=1$ et $(\rho,a',1)$ intervient avec multiplicit\'e 1 dans $Jord(\tilde{\psi}')$ et de plus si $x_{\ell}=-y_{\ell}$ alors $(\rho,2x_{\ell}+1,1)\notin Jord(\tilde{\psi}')$. On aura aussi dans les applications une hypoth\`ese  qui simplifie les d\'emonstrations, alors on la fait: si $x_{\ell}= -y_{\ell}$ alors $x_{\ell}=(\alpha-1)/2+T$. Dans le cas o\`u $x_{\ell}=-y_{\ell}$, on pourrait remplacer dans (1), $<x_{i},y_{i}>_{\rho}\times \tilde{\pi}'$ par une repr\'esentation $\tilde{\pi}''\in \Pi(\tilde{\psi}'')$, o\`u $\tilde{\psi}''$ se d\'eduit de $\tilde{\psi}'$ en ajoutant 2 fois $(\rho,2x_{\ell}+1,1)$; donc ce que l'on veut c'est une multiplicit\'e inf\'erieure ou \'egal \`a 2 pour les \'el\'ements de la forme $(\rho,a',1)$ avec $(a'-1)/2 \in [(\alpha-1)/2+T,(\alpha+1)/2]$,  avec multiplicit\'e 2 uniquement \'eventuellement pour $(\rho,\alpha+2T,1)$.

Les premiers facteur $\times_{i\in [1,\ell[}<x_{i},y_{i}>_{\rho}$ jouent un r\^ole muet mais  ils sont l\`a dans les applications.

\begin{lem} On suppose que $T=1$; $Jac_{(\alpha+1)/2}\tilde{\pi}\neq 0$ n\'ecessite que l'une des 3 hypoth\`eses suivante soit r\'ealis\'ee (elles ne sont pas exclusives l'une de l'autre), $x_{\ell}=(\alpha+1)/2$, $-y_{\ell}=(\alpha+1)/2$, $(\rho,\alpha+2,1)\in Jord(\tilde{\psi}')$. Et $\pi$ v\'erifie alors l'une des inclusions ci-dessous:

(i) $
\pi\hookrightarrow \times_{i\in [1,\ell[}<x_{i},y_{i}>_{\rho}\times <x_{\ell}-1,y_{\ell}>_{\rho}\times \tilde{\pi}'
$
ce cas n\'ecessite que $x_{\ell}=(\alpha+1)/2$;

(ii) 
$
\pi\hookrightarrow \times _{i\in [1,\ell[}<x_{i},y_{i}>_{\rho}\times \pi''
$ 
o\`u $\pi''\in \Pi({\psi}'')$ avec ${\psi}''$ qui se d\'eduit de $\tilde{\psi}'$ en ajoutant $(\rho,2x_{\ell},1)$ et $(\rho,2x_{\ell}-1,1)$, et ce cas n\'ecessite $x_{\ell}=-y_{\ell}=(\alpha+1)/2$;

(iii) $\pi\hookrightarrow \times_{i\in [1,\ell[}<x_{i},y_{i}>_{\rho}\times <x_{\ell},y_{\ell}+1>_{\rho}\times \tilde{\pi}'
$ ce cas n\'ecessite  $-y_{\ell}=(\alpha+1)/2>x_{\ell}$;

(iv) si $(\rho,\alpha+2,1)\in Jord(\tilde{\psi}')$, en notant $\tilde{\psi}''$ le morphisme qui se d\'eduit de $\tilde{\psi'}$ en rempla\c{c}ant $(\rho,\alpha+2,1)$ par $(\rho,\alpha,1)$, $\pi\hookrightarrow \times_{i\in [1,\ell]}<x_{i},y_{i}>_{\rho}\times \tilde{\pi}'',$ avec $\tilde{\pi}''\in \Pi(\tilde{\psi}'')$.
\end{lem}
Avant de faire la preuve on remarque que l'on voudra appliquer le lemme par induction pour obtenir le cas $T>1$; le lemme appliqu\'e \`a $\pi_{T-1}:=Jac_{(\alpha-1)/2+T}\tilde{\pi}$ montre que cette repr\'esentation v\'erifie \`a peu pr\`es les  hypoth\`eses n\'ecessaires mais pas tout \`a fait; on ajoute tout de suite une hypoth\`ese que l'on aura: si $x_{\ell}=-y_{\ell}$ alors $x_{\ell}=(\alpha-1)/2+T $. Pour la suite, il est utile de montrer tout de suite ce qu'il manque \`a $\pi_{T-1}$ pour pouvoir lui r\'eappliquer le lemme: 

avec (i) et (ii), il n'y a aucun probl\`eme;

avec (iii), il manque si $x_{\ell}=-y_{\ell}-1$ le fait que $(\rho,2x_{\ell}+1,1)\notin Jord(\tilde{\psi}')$; on v\'erifiera dans les applications que si (iii) se produit on n'a pas $x_{\ell}=-y_{\ell}-1$;

avec (iv), il y a un probl\`eme si $(\rho,\alpha+2T-2,1)\in Jord(\tilde{\psi}')$; dans ce cas, il faut changer $\ell$ en $\ell+1$, car on r\'ealise dans ce cas $\tilde{\pi}''$ comme sous-module d'une induite de la forme $St(\rho,\alpha+2T-2)\times \pi''$ avec $\pi''$ dans le paquet qui se d\'eduit de $\tilde{\psi}''$ en enlevant les 2 copies de $(\rho,\alpha+2T-2,1)$; donc il faut que le segment $[x_{\ell},y_{\ell}]$ v\'erifie les hypoth\`eses des segments $[x_{i},y_{i}]$ pour $i\in [1,\ell[$. On v\'erifiera ce point dans les applications.

\

On montre maintenant le lemme. On applique \'evidemment les formules standard pour les modules de Jacquet. On obtient une filtration du module de Jacquet et par fonctorialit\'e $\pi$ est sous-module d'un des gradu\'es et on ne sait \'evidemment pas lequel. Vues les hypoth\`eses que l'on a mises, la filtration a au plus 2 termes et on obtient directement que l'une des inclusions suivantes est r\'ealis\'ee:

soit (i), soit (iii) mais \'eventuellement m\^eme si $x_{\ell}=y_{\ell}$ soit (iv) avec $\tilde{\pi}'':=Jac_{(\alpha+1)/2}\tilde{\pi}$ dont on sait que cette repr\'esentation est irr\'eductible si non nulle gr\^ace \`a l'hypoth\`ese que pour tout $(\rho,a',b')\in Jord(\tilde{\psi}')$ si $(a'-b')/2=(\alpha+1)/2$ alors $b'=1$ et cet \'el\'ement n'intervient qu'avec mulitplicit\'e 1 dans $Jord(\tilde{\psi}')$. Donc le seul point est de montrer que (iii) avec $x_{\ell}=-y_{\ell}$ est en fait (ii).

On suppose donc que $x_{\ell}=-y_{\ell}=(\alpha+1)/2$ et que 
$$
\pi\hookrightarrow \times_{i\in [1,\ell[}<x_{i},y_{i}>_{\rho}\times <(\alpha+1)/2,-(\alpha-1)/2>_{\rho}\times \tilde{\pi}'.\eqno(2)
$$
On utilise l'inclusion $<(\alpha+1)/2,-(\alpha-1)/2>_{\rho}\hookrightarrow \rho\vert\,\vert^{(\alpha+1)/2}\times St(\rho,\alpha)$. On sait que $St(\rho,\alpha)\times \tilde{\pi}'$ est semi-simple form\'e d'\'el\'ements dans le paquet associ\'e \`a $\psi'$ le morphisme qui se d\'eduit de $\tilde{\psi}'$ en ajoutant 2 copies de $(\rho,\alpha,1)$. Donc il existe $\pi'$ dans $\Pi(\psi')$ tel que l'inclusion (2) donne en fait une inclusion
$$
\pi\hookrightarrow \times_{i\in [1,\ell[}<x_{i},y_{i}>_{\rho}\times\rho\vert\,\vert^{(\alpha+1)/2}\times \pi'.\eqno(3)
$$
On peut remplacer dans l'induite de (3), la repr\'esentation $\rho\vert\,\vert^{(\alpha+1)/2}\times \pi'$ par un sous-quotient irr\'educ\-tible ${\pi}''$, d'o\`u
$$
\pi\hookrightarrow \times_{i\in [1,\ell[}<x_{i},y_{i}>_{\rho}\times{\pi}''.\eqno(4)
$$ et on va montrer que n\'ecessairement ${\pi}''$ est dans $\Pi({\psi}'')$ o\`u ${\psi}''$ s'obtient comme dans l'\'enonc\'e de (ii). On commence par remarquer que les hypoth\`eses sur les $x_{i},y_{i}$ pour $i\in [1,\ell[$ entra\^{\i}nent que dans (3), on peut pousser $\rho\vert\,\vert^{(\alpha+1)/2}$, en premi\`ere position; d'o\`u $Jac_{(\alpha+1)/2}\pi\neq 0$. Quand on reporte cette propri\'et\'e dans (4), on voit que n\'ecessairement $Jac_{(\alpha+1)/2}\pi''\neq 0$ car aucun des $x_{i}$ ou des $-y_{i}$ ne vaut $(\alpha+1)/2$ si $i\in [1,\ell[$. Or $Jac_{(\alpha+1)/2}$ appliqu\'e \`a l'induite $\rho\vert\,\vert^{(\alpha+1)/2}\times \pi'$ est irr\'eductible (\'egal \`a $\pi'$) car $Jac_{(\alpha+1)/2}{\pi'}=0$ car $Jord({\psi'})$ ne contient aucun \'el\'ement de la forme $(\rho,a',b')$ avec $(a'-b')/2=(\alpha+1)/2$ par hypoth\`ese dans le cas o\`u $x_{\ell}=-y_{\ell}$. Donc $\pi''$ est l'unique sous-module irr\'eductible de l'induite $\rho\vert\,\vert^{(\alpha+1)/2}\times {\pi}'$. Or l'application $Jac_{(\alpha+1)/2}$ d\'efinit une surjection de $\Pi(\psi'')$ sur $\Pi(\psi')$ (cf. la remarque de \ref{rappelconstruction}); on fait remarquer au lecteur que la d\'emonstration de ce fait utilise la formule des traces via les travaux d'Arthur et que je n'en conna\^{\i}t pas de d\'emonstration uniquement \`a l'aide de module de Jacquet; cela me semble donc un r\'esultat non trivial. Ainsi l'\'el\'ement $\tilde{\pi}''\in \Pi(\psi'')$ qui s'envoit sur $\pi'$ est l'unique sous-module irr\'eductible de l'induite $\rho\vert\,\vert^ {(\alpha+1)/2}\times \pi'$ et co\"{\i}ncide donc avec $\pi''$. Cela termine la preuve.

\subsection{Description partielle de certaines repr\'esentations\label{blocpositif}}
En \ref{proprietesimpledejac} (5), on a montr\'e comment on pouvait remplacer $(\rho,a,b)\in Jord(\psi)$ tel que $a\geq b>1$ par soit $(\rho,a,b-2)$ soit par $\cup_{c\in [(a-b)+1,a+b-1]; (-1)^c=(-1)^{a+b-1}}(\rho,c,1)$. Mais on avait une hypoth\`ese forte que pour tout $(\rho,a',b')\in Jord(\psi)$,  soit $(a'-b')/2\leq (a-b)/2$ soit $a'>>b'$. Il faut lever cette hypoth\`ese, mais on ne peut plus alors qu'esp\'erer un r\'esultat approximatif comme expliqu\'e en \ref{approximationjac}. C'est donc un r\'esultat technique qui permet de d\'emontrer des propri\'et\'es des repr\'esentations consid\'er\'ees mais ne permet pas de les construire.

On fixe $\psi$ un morphisme et $(\rho,a,b)\in Jord(\psi)$. On suppose que pour tout $(\rho,a',b')\in Jord(\psi)$ v\'erifiant $(a'-b')/2\in ](a-b)/2,(a+b)/2-1[$, $b'=1$ et qu'un tel \'el\'ement n'intervient qu'avec multiplicit\'e 1. On fixe un ordre sur $Jord(\psi)$ tel que les \'el\'ements de $Jord(\psi)$ strictement plus grands que $(\rho,a,b)$ soient exactement les \'el\'ements $(\rho,a',b')\in Jord(\psi)$ avec $(a'-b')/2>(a-b)/2$. On note $\mathcal{E}$ l'ensemble des entiers $a'$ tel que $(a'-1)/2\in ](a-b)/2,(a+b)/2-1[$ et $(\rho,a',1)\in Jord(\psi)$. Soit $\pi\in \Pi(\psi)$; gr\^ace \`a l'ordre mis sur $Jord(\psi)$, on d\'efinit les param\`etres $\underline{t}, \underline{\eta}$ qui permettent de construire $\pi$ et on pose $t:=\underline{t}(\rho,a,b)$.
On note $\psi'$ le morphisme qui se d\'eduit de $\psi$ en rempla\c{c}ant $(\rho,a,b)$ par $(\rho,a,b-2)$.

\begin{lem}(i) On suppose que $t=0$. Il existe $v$ un entier positif  ou nul,  un morphisme $\psi''$, $\pi''\in \Pi(\psi'')$ et des \'elements $a_{i}$ pour $i\in [1,v]$ de $\mathcal{E}$ rang\'es dans l'ordre d\'ecroissant avec une inclusion

si $v$ est pair
$$
\pi\hookrightarrow \times_{j\in [1,[v/2]]}<(a_{2j}-1)/2,-(a_{2j-1}-1)/2>_{\rho}\times\pi'',
$$
o\`u $\psi''$ s'obtient \`a partir de $\psi'$ en ajoutant $(\rho,a-b+1,1)$ et $(\rho,a+b-1,1)$ et en enlevant tous les $(\rho,a_{j},1)$, pour $j\in [1,v]$; si $v=0$, il n'y a pas d'induite, $\pi=\pi''$ tout simplement

si $v$ est impair
$$
\pi\hookrightarrow \times_{j\in [1,[v/2]]}<(a_{2j+1}-1)/2,-(a_{2j}-1)/2>_{\rho}\times <(a_{1}-1)/2,-(a+b)/2+1>_{\rho},
$$
o\`u $\psi''$ s'obtient \`a partir de $\psi'$ en ajoutant $(\rho,a-b+1,1)$ et en enlevant tous les $(\rho,a_{j},1)$ pour $j\in [1,v]$.

(ii) On suppose que $t>0$. Il existe $v$ un entier positif ou nul, un 
morphisme $\psi''$, $\pi''\in \Pi(\psi'')$ et des \'elements $a_{i}$ pour $i\in [1,v]$ de $\mathcal{E}$ rang\'es dans l'ordre d\'ecroissant avec une inclusion

si $v$ est impair
$$
\pi\hookrightarrow <(a-b)/2,-a_{v}>\times_{j\in [1,[v/2]]} <(a_{2j}-1)/2,-(a_{2j-1}-1)/2>_{\rho}\times \pi'',
$$
o\`u $\psi''$ s'obtient \`a partir de $\psi'$ en ajoutant $(\rho,a+b-1,1)$ et en enlevant tous les $(\rho,a_{j},1)$ pour $j\in [1,v]$;

si $v$ est pair

$$
\pi\hookrightarrow <(a-b)/2,-a_{v}>\times_{j\in [1,v/2[} <(a_{2j+1}-1)/2,-(a_{2j}-1)/2>_{\rho}\times <(a_{1}-1)/2,-(a+b)/2+1>_{\rho}\times \pi'',
$$
o\`u $\psi''$ s'obtient \`a partir de $\psi'$  en enlevant tous les $(\rho,a_{j},1)$ pour $j\in [1,v]$; si $v=0$, on a $\pi\hookrightarrow <(a-b)/2, -(a+b)/2+1>_{\rho}\times \pi''$ tout simplement.
\end{lem}
On fixe un morphisme dominant $\psi$ et plus pr\'ecis\'ement dominant tous les blocs de Jordan de $\psi$ strictement sup\'erieur \`a $(\rho,a,b)$ et contenant les autres. On le note $\tilde{\psi}$. On sait qu'il existe $\tilde{\pi}\in \Pi(\tilde{\psi})$ tel que $\pi=Jac_{z\in \mathcal{Z}}\tilde{\pi}$ pour un ensemble totalement ordonn\'e $\mathcal{Z}$ convenable. On fait d'abord remarquer que le lemme est vrai pour $\tilde{\pi}$ en prenant $v=0$ puisque $\mathcal{E}$ est vide et que l'on peut appliquer  \ref{proprietesimpledejac} (5). Le point est donc de faire redescendre les \'el\'ements de $Jord(\tilde{\psi})$ qui dominent les \'el\'ements de $\mathcal{E}$; si on a le lemme pour cette repr\'esentation interm\'ediaire, avec $\tilde{\pi}''\in \Pi(\tilde{\psi''})$, on peut ensuite faire redescendre les \'el\'ements restant et cela ne touche qu'\`a $\pi''$ pour passer de $\tilde{\psi}''$ \`a $\psi''$ sans modifier les inclusions. On peut donc supposer tout de suite que pour tout $(\rho,a',b')\in Jord(\psi)$ avec $(a'-b')/2\geq (a+b)/2-1$ on a en fait $(a'-b')/2>> (a+b)/2-1$.

On fait la preuve de (i) et on suppose donc que $t=0$. On note $\tilde{\psi}_{0}$ le morphisme qui se d\'eduit de $\tilde{\psi}$ en rempla\c{c}ant $(\rho,a,b)$ par $(\rho, a-b+1,1),(\rho,a,b-2),(\rho,a+b-1,1)$. On note $\alpha$ le plus petit \'el\'ement de $\mathcal{E}$; il existe $T>>0$ tel que $(\rho,\alpha+2T,1)\in Jord(\tilde{\psi}_{0})$ et il faut d'abord calculer:
$$
Jac_{(\alpha-1)/2+T, \cdots, (\alpha+1)/2}\tilde{\pi}.
$$
On va appliquer \ref{approximationjac} pour calculer d'abord $Jac_{(\alpha-1)/2+T, \cdots, (a+b)/2}\tilde{\pi}$; ici on a $\ell=0$ et \`a chaque fois seul (iv) peut se produire. Donc on obtient
$$
Jac_{(\alpha-1)/2+T, \cdots, (a+b)/2}\tilde{\pi} \hookrightarrow <(a+b)/2-1, -(a+b)/2+1>_{\rho}\times \pi''
$$
pour $\pi''\in \Pi(\psi'')$ qui se d\'eduit de $\tilde{\psi}_{0}$ en enlevant $(\rho,\alpha+2T,1)$ et $(\rho,a+b-1,1)$. On applique encore $Jac_{(a+b)/2-1}$ en utilisant \ref{approximationjac}; peuvent se produire (i) et (ii). Supposons que (ii) se produit, cela veut dire que le r\'esultat est dans un paquet qui se d\'eduit de $\tilde{\psi}_{0}$ en rempla\c{c}ant simplement $(\rho,\alpha+2T,1)$ par $(\rho,a+b-3,1)$ et ensuite les $Jac$ suivant ne peuvent plus que faire descendre $(\rho,a+b-3,1)$ en $(\rho,\alpha,1)$; en fait il ne s'est pas pass\'e grand chose et on pourra recommencer avec l'\'el\'ement de $\mathcal{E}$ juste au dessus de $\alpha$. Par contre si le cas (i) se produit, la situation va changer; donc on peut consid\'erer que $\alpha$ est le plus petit \'el\'ement de $\mathcal{E}$ pour lequel dans la proc\'edure pr\'ec\'edente c'est le cas (i) qui se produit. Ici, on a donc
$$
Jac_{(\alpha-1)/2+T, \cdots, (a+b)/2-1}\tilde{\pi} \hookrightarrow <(a+b)/2-3, -(a+b)/2+1>_{\rho}\times \pi''
$$
Et les $Jac$ suivant donnent encore uniquement le cas (i) et finalement on obtient
$$
Jac_{(\alpha-1)/2+T, \cdots, (\alpha+1)/2}\tilde{\pi} \hookrightarrow <\alpha, -(a+b)/2+1>_{\rho}\times \pi''.\eqno(1)
$$
L'\'el\'ement $\alpha$ qui est ici, sera le $\alpha_{v}$ de l'\'enonc\'e bienque l'on ne connaisse par encore $v$. On note $\beta$ l'\'el\'ement de $\mathcal{E}$ qui est juste au dessus de $\alpha$ et on calcule pour $T'$ grand tel que $(\rho,\beta+2T',1)\in Jord(\tilde{\psi}_{0})$, $Jac_{(\beta-1)/2+T', \cdots, (\beta+1)/2}$ du membre de gauche de (1); on commence par calculer $Jac_{(\beta-1)/2+T', \cdots, (a+b)/2}$ et seul (iv) peut se produire. Quand on continue avec $Jac_{(a+b)/2-1}$ soit (ii) soit (iv) se produit et comme ci-dessus, il ne se passe ensuite plus rien, c'est toujours (iv) qui seul peut s'appliquer et le r\'esultat est simplement de faire descendre $(\rho,\beta+2T',1)$ en $(\rho,\beta,1)$; comme ci-dessus, on ne s'int\'eresse donc qu'au premier $\beta$ pour lequel c'est (ii) qui se produit.  Et on trouve
$$
Jac_{(\beta-1)/2+T',\cdots, (\beta+1)/2}Jac_{(\alpha-1)/2+T, \cdots, (\alpha+1)/2}\tilde{\pi} \hookrightarrow <\alpha, -(\beta-1)/2>_{\rho}\times \pi'''
$$
o\`u $\pi'''$ se trouve dans le morphisme qui se d\'eduit de $\tilde{\psi}_{0}$ en enlevant gardant $(\rho,a+b-1,1)$ et en enlevant $(\rho,\alpha+2T,1)$ et $(\rho,\beta+2T,1)$. On remarque que pour tout $\gamma\in \mathcal{E}$ qui doit encore redescendre, on a certainement $(\gamma-1)/2>(\beta-1)/2>(\alpha-1)/2$. Le facteur $<\alpha, -(\beta-1)/2>_{\rho}$ joue le r\^ole muet comme voulu dans \ref{approximationjac} et on est revenu \`a la situation o\`u on faisait redescendre $\alpha$; il est alors facile de finir la preuve du lemme (par exemple par r\'ecurrence sur $\vert \mathcal{E}\vert$). Cela termine la preuve de (i).

Preuve de (ii); on sait d\`es le d\'epart que l'on a une inclusion
$$
\tilde{\pi}\hookrightarrow <(a-b)/2, \cdots, -(a+b)/2+1>\times \tilde{\pi}'.
$$
C'est donc la preuve de (i) qui s'applique telle quelle mais en commen\c{c}ant au cas $\beta$. Cela termine la preuve du lemme.
\subsection{Le cas des blocs n\'egatifs\label{blocnegatif}}
On fixe $\psi$ et on suppose que pour tout $(\rho,a,b)\in Jord(\psi)$ tel que $a\geq b$, on a $b=1$ et qu'un tel \'el\'ement n'intervient qu'avec multiplicit\'e 1. On note $(\rho,a,b)$ un \'el\'ement de $Jord(\psi)$ tel que $(b-a)/2>0$ et minimal avec cette propri\'et\'e. On fixe $\pi\in \Pi(\psi)$ et on suppose que si $a=b+1$,  $\pi\notin \Pi(\psi^{\check{\empty}})$ o\`u $\psi^{\check{\empty}}$ s'obtient en changeant uniquement $(\rho,a,b)$ en $(\rho,b,a)$. On fixe un ordre sur $Jord(\psi)$ tel que $(\rho,a,b)$ soit le plus petit \'el\'ement. On note $\mathcal{E}$ l'ensemble des \'el\'ements de $Jord(\psi)$ de la forme $(\rho,a',1)$ avec $a'\in [(b-a)-1, a+b-1[$.
\begin{lem} Il existe un entier $v$ positif ou nul,  des \'el\'ements de $\mathcal{E}$, $a_{i}$ pour $i\in [1,v]$ rang\'es dans l'ordre d\'ecroissant, un morphisme $\psi''$ et $\pi''\in \Pi(\psi'')$, avec une inclusion
$$
\pi\hookrightarrow <(a-b)/2, \cdots, -(a_{v}-1)/2>_{\rho}\times_{j\in [1,v[; j\equiv v[2]} <(a_{j+1}-1)/2,-(a_{j}-1)/2>_{\rho}$$
$$\times \begin{cases} <(a_{1}-1)/2, -(a+b)/1+1>_{\rho} \times \pi''; \hbox{ si $v$ est pair}\\
\times \pi'', \hbox{ si $v$ est impair}
\end{cases}
$$
o\`u $\psi''$ s'obtient \`a partir de $\psi$ en rempla\c{c}ant $(\rho,a,b)$ par $(\rho,a,b-2)$ et en enlevant les $(\rho,a_{i},1)$ pour $i\in [1,v]$ et en ajoutant $(\rho, a+b-1,1)$ si $v$ est impair; le terme $<(a-b)/2,-(a_{1}-1)/2>_{\rho}$ dispara\^{\i}t si $a_{1}=b-a-1$ et si $v=0$, on n'a que le termes $<(a-b)/2,-(a+b)/2+1>_{\rho}$.
\end{lem}
D'abord on remarque qu'il suffit de d\'emontrer le lemme avec l'hypoth\`ese suppl\'ementaire que pour tout $(\rho,a',b')\in Jord(\psi)$ si $(\rho,a',b')\notin \{(\rho,a,b),(\rho,a'',1); a''<(a+b-1)$ alors $\vert (a'-b')\vert >>0$. En effet, on met un ordre sur $Jord(\psi)$ de tel sorte que ces \'el\'ements sont strictement plus grands que $(\rho,a,b)$ et tout $(\rho,a'',1)$ avec $a''<(a+b-1)$; on consid\`ere un morphisme qui domine $\psi$ et plus pr\'ecis\'ement qui domine tous ces $(\rho,a',b')\in Jord(\psi)$ et au contraire contient $(\rho,a,b)$ et $(\rho,a'',1)$ pour tout $a''\in \mathcal{E}$. On sait qu'il existe $\tilde{\pi}$ dans le paquet associ\'e \`a ce morphisme tel que 
$$
\pi=Jac_{y\in \mathcal{Y}}\tilde{\pi},
$$
o\`u les \'el\'ements de $\mathcal{Y}$ sont soit n\'egatifs strictement plus petits que $(a-b)/2$ soit positifs strictement plus grands que $(a+b)/2-1$. Donc prendre ce module de Jacquet ne perturbe pas les inclusions du lemme si on les conna\^{\i}t pour $\tilde{\pi}$, le module de Jacquet ne s'applique qu'\`a l'analogue de $\pi''$; ici on a remplac\'e $(\rho,a,b)$ par $(\rho,a,b-2)$ il faut v\'erifier que cet \'el\'ement peut toujours \^etre le plus petit \'el\'ement pour un bon ordre. Pour cela on utilise la minimalit\'e de $(b-a)/2$; on a maintenant soit $(a-b+2)/2=1/2$ soit $(b-2-a)/2$ est le plus petit des $(b'-a')/2$ parmi ceux qui sont n\'egatifs. D'o\`u l'assertion. Ainsi quand on applique $Jac_{y\in \mathcal{Y}}$ \`a l'analogue de $\pi''$ on obtient  un $\pi''$ convenable. D\'emonstrons donc le lemme avec cette hypoth\`ese suppl\'ementaire. On fixe $\psi_{>>}$ dominant $\psi$ et il faut uniquement faire redescendre les \'el\'ements de $Jord(\psi_{>>})$ dominant les $(\rho,a'',1)$ avec $a''<a+b-1$. Au d\'epart on sait qu'il existe une repr\'esentation $\pi_{>>}$ tel que 
$$
\pi=\circ_{a''}Jac_{(a''-1)/2+T_{a''}, \cdots, (a''+1)/2}\pi_{>>},
$$o\`u $a''$ parcourt les entiers tels que $(\rho,a'',1)\in Jord(\psi)$ et $a''< a+b-1$ les $T_{a''}$ sont des entiers tels que $T_{a''}>>T_{a'''}$ si $a''>a'''$ et o\`u on prend les \'el\'ements de $\mathcal{E}$ dans l'ordre croissant.

De plus on note $\psi_{>>,-}$ le morphisme qui se d\'eduit de $\psi_{>>}$ en rempla\c{c}ant $(\rho,a,b)$ par $(\rho,a,b-2)$ et on sait (cf. (7) et le lemme de \ref{proprietesimpledejac}) qu'il existe $\pi_{>>,-}\in \Pi(\psi_{>>,-})$ avec une inclusion
$$
\pi_{>>}\hookrightarrow <(a-b)/2, \cdots, -(a+b)/2+1>_{\rho}\times \pi_{>>,-}.
$$ On fait d'abord descendre les \'el\'ements qui dominent les $(\rho,a'',1)$ avec $a''<a-b-1$. Ici il n'y a pas d'autres solutions que d'appliquer le module de Jacquet \`a $\pi_{>>,-}$: a priori quand on calcule $Jac_{z\in [(a''-1)/2+T_{a''},(a''+1)/2]}$ a une induite comme ci-dessus, le segment se coupe en 2, une partie sert \`a prendre le module de Jacquet de $\pi_{>>,-}$ et l'autre s'applique \`a $<(a-b)/2, \cdots, -(a+b)/2+1>_{\rho}$ et \`a son dual. Le d\'ecoupage se fait n\'ecessairement en 2 sous-intervalle et le sous-intervalle qui s'applique \`a $\pi_{>>,-}$ commence n\'ecessairement \`a $(a''-1)/2+T_{a''}$. Puisque $(a''+1)/2< -(a-b)/2$ ce sous-intervalle contient aussi $(a''+1)/2$ et co\"{\i}ncide donc avec tout l'intervalle.
Ceci r\`egle facilement les \'el\'ements $(\rho,a'',1)$ avec $a''$ non dans $\mathcal{E}$; on ne change pas les notations, $\pi_{>>}$ et $\pi_{>>,-}$ alors qu'il le faudrait. On consid\`ere maintenant les \'el\'ements dans $\mathcal{E}$ et on est exactement dans la situation de \ref{blocpositif} (ii), le fait que $(a-b)/2$ est n\'egatif au lieu d'\^etre positif n'a qu'un seul effet que le premier terme $<(a-b)/2, -(a_{1}-1)/2>_{\rho}$ dispara\^{\i}t si $a_{1}=b-a-1$. Cela termine la preuve.

\subsection{Une remarque sur les orbites unipotentes\label{csqsurorbite}}
On consid\`ere un entier $v$ et des $a_{i}\in ]a-b+1,a+b-1[$ pour $i\in [1,v]$; on consid\`ere l'orbite unipotente, $O_{v}$ dont les blocs de Jordan sont $(a-b+1)$ et 2 copies de chacun des $(a_{j}+a_{j+1})/2$ pour $j$ de la parit\'e de $v-1$ auquel on ajoute encore $(a+b-1)$ si $v$ est pair et 2 copies de $(a+b-1+a_{1})/2$ si $v$ est impair.

\begin{lem}L'orbite $O_{v}$ ainsi d\'ecrite contient dans sa fermeture l'orbite dont les blocs de Jordan sont 2 copies de $a$ et tous les $a_{j}$ pour $j\in [1,v]$.
\end{lem}Pour tout $t\in [1,v+2]$ on compare la somme des $t$ plus grands bloc de Jordan. On commence par $t=1$. On a clairement 
$
a+b-1\geq a
$ car $b\geq 1$ et $a+b-1>a_{1}$ par hypoth\`ese, ; d'o\`u $a+b-1\geq sup(a,a_{1})$, ce qui r\`egle $t=1$ et $v$ pair. Mais on a aussi $((a+b-1+a_{1})/2 \geq a_{1}$ et comme $a_{1}>a-b+1$, $((a+b-1+a_{1})/2\geq ((a+b-1+a-b+1)/2=a$ ce qui r\`egle le cas de $t=1$ et $v$ impair.

Le cas $t=v+2$ est simplement le fait que les orbites sont relatives au m\^eme groupe et le cas $t=v+1$ revient \`a comparer les derniers blocs de Jordan (chaque orbite en a le m\^eme nombre); on v\'erifie que $a-b+1 \leq inf(a,a_{v})$ presque par hypoth\`ese.

On suppose donc que $t\in ]1,v]$. Si $t$ est de la parit\'e de $v$, la somme des $t$ plus grands blocs de Jordan de $O_{v}$ vaut
$$
a+b-1+\sum_{j\leq t-2}a_{j}+(a_{t-1}+a_{t})/2
$$
et on doit comparer avec $\sum_{j\leq t-2}a_{j}+sup(a_{t-1},a)+sup(a_{t},a)$. La diff\'erence de ces 2 nombres vaut $\bigl((a+b-1+a_{t-1})/2-sup(a,a_{t-1})\bigr)+\bigl((a+b-1+a_{t})/2-sup(a,a_{t})\bigr)$ et chacune des parenth\`eses est positive. Si $t$ est de la parit\'e oppos\'ee \`a celle de $v$, la somme des $t$ plus grands blocs de Jordan de $O_{v}$ vaut
$
a+b-1+\sum_{j\leq t-1}a_{j}$ que l'on doit comparer \`a $\sum_{j\leq t-2}a_{j}+sup(a_{t-1},a)+sup(a_{t},a)$. La diff\'erence de ces 2 nombres vaut $\bigl((a+b-1+a_{t-1})/2-sup(a,a_{t-1})\bigr)+\bigl((a+b-1+a_{t-1})/2-sup(a,a_{t})\bigr)$ et chacune des parenth\`eses est positive car $a_{t-1}>a_{t}$. Cela termine la preuve.

\

On a un lemme analogue pour le (ii) de \ref{blocpositif}. On consid\`ere l'orbite dont les blocs de Jordan sont 2 copies de $(a-b+1+a_{v})/2$ ainsi que des $(a_{j}+a_{j+1})/2$ pour $j$ de la parit\'e de $v$ auquel on ajoute, si $v$ est pair, 2 copies de $(a+b-1+a_{1})/2$ et si $v$ est impair une copie de $a+b-1$. On note $O'$ cette orbite.
\begin{lem}L'orbite $O'$ contient dans sa fermeture l'orbite, $O$, dont les blocs de Jordan sont les $a_{j}$ pour $j\in [1,v]$ auxquels on ajoute 2 copies de $a$.
\end{lem}
Le plus grand bloc de Jordan de $O'$ est soit $a+b-1$ soit $(a+b-1+a_{1})/2$; dans les 2 cas ce bloc est sup\'erieur (d'ailleurs strictement) \`a $sup(a,a_{1})$. Le plus petit bloc de Jordan de $O'$ vaut $(a-b+1+a_{v})/2$; ceci est inf\'erieur \`a $a_{v}$ car $a-b+1<a_{v}$. Comme $a_{v}<a+b-1$ ceci est aussi inf\'erieur \`a $a$. On a donc montrer que le plus grand bloc de $O'$ est sup\'erieur au plus grand bloc de $O$ et la m\^eme assertion pour la somme des $v+1$ plus grands blocs de chacune de ses orbites. Il faut ensuite faire la m\^eme comparaison pour tout $t\in [2,v]$; supposons $t$ de la parit\'e de $v$. Il faut calculer le signe de
$$
a+b-1 +\sum_{i<t}a_{i}-\sum_{i<t-1}a_{i}-sup(a_{t-1},a)-sup(a_{t},a)$$
$$=a+b-1+a_{t-1}-sup(a,a_{t-1})-sup(a,a_{t})>0
$$
par un calcul d\'ej\`a fait. Si $t$ est de la parit\'e oppos\'ee \`a $v$, il faut calculer le signe de
$$
a+b-1+\sum_{i<t-1}a_{i}+(a_{t-1}+a_{t})/2-\sum_{i<t-1}a_{i}-sup(a_{t-1},a)-sup(a_{t},a)$$
$$=((a+b-1+a_{t-1})/2-sup(a,a_{t-1})+((a+b-1+a_{t})/2-sup(a,a_{t}))>0
$$
par un calcul d\'ej\`a fait. Cela termine la preuve.

\

On a encore un lemme analogue pour \ref{blocnegatif}. On reprend les hypoth\`eses de cette r\'ef\'erence. On note $O'$ l'orbite dont les blocs de Jordan sont 2 copies de $(a_{1}-(b-a-1))/2$, 2 copies des $(a_{j}+a_{j+1})/2$ pour $j\in [1,v]$ de la parit\'e de $v$ et si $v$ est pair, 2 copies de $(a+b-1+a_{1})/2$ tandis que si $v$ est impair une copie de $a+b-1$. On note $O$ l'orbite dont les blocs de Jordan sont 2 copies de $a$ et chaque $a_{j}$ pour $j\in [1,v]$
\begin{lem}L'orbite $O'$ contient dans sa fermeture l'orbite $O$.
\end{lem}
On compare d'abord le plus grand bloc de Jordan de chacune de ces orbites; si $v$ est pair, $(a+b-1+a_{1})/2> a_{1}$ puisque $a+b-1>a_{1}$. De plus $a_{1}\geq b-a-1$, d'o\`u $(a+b-1+a_{1})/2\geq b-1\geq a$ car $b-a>0$ par hypoth\`ese. Si $v$ est impair, on a clairement $a+b-1>sup(a,a_{1})$. L'orbite $O$ a $v+2$ blocs de Jordan et $O'$ aussi sauf si $a_{1}=b-a-1$ o\`u elle en a un de moins. Le $v+2$-i\`eme bloc de Jordan de $O$ est $inf(a,a_{1})$; celui de $O'$ est $(a_{1}-(b-a-1))/2$. Or $(a_{1}-(b-a-1))/2<( a+b-1-(b-a-1))/2=a$ et est aussi inf\'erieur \`a $a_{1}$; donc la somme des $v+1$ plus grands blocs de Jordan de $O'$ est plus grande que son analogue pour $O$. Soit $j\in ]1,v]$, il faut comparer la somme des $j$ plus grands blocs de Jordan pour $O'$ a son analogue pour $O$. On suppose d'abord que $j$ est de la parit\'e de $v$; on doit calculer le signe de
$$
a+b-1+\sum_{k<j}a_{k}-\sum_{k<j-1}a_{k}-sup(a,a_{j-1})-sup(a,a_{j})$$
$$=a+b-1+a_{j-1}-sup(a,a_{j-1})-sup(a,a_{j})$$
$$= ((a+b-1+a_{j-1})/2-sup(a,a_{j-1}))-((a+b-1+a_{j-1})/2-sup(a,a_{j}))$$
et chaque parenth\`ese est positive par l'argument d\'ej\`a donn\'ee pour la premi\`ere et le m\^eme en tenant compte de $a_{j-1}>a_{j}$ pour la deuxi\`eme. On suppose maintenant que $j$ est de la parit\'e oppos\'e \`a $v$, on doit alors calculer le signe de:
$$
a+b-1+\sum_{k<j-1}a_{k}+(a_{j-1}+a_{j})/2-\sum_{k<j-1}a_{k}-sup(a,a_{j-1})-sup(a,a_{j})$$
$$=((a+b-1+a_{j-1})/2-sup(a,a_{j-1}))-((a+b-1+a_{j})/2-sup(a,a_{j}))$$
est encore positif. Et cela termine la preuve.
\section{Comparaison des param\'etrisations de Langlands}
\subsection{D\'efinition des orbites unipotentes \`a la Langlands\label{definitionorbitelanglands}}
Soit $\pi$ une repr\'esentation irr\'eductible et on consid\`ere sa param\'etrisation de Langlands, d'abord d'un point de vue th\'eorique, \`a $\pi$ on associe entre autre un morphisme $\phi_{\pi}$ de $W_{F}\times SL(2,{\mathbb C})$ dans le $L$-groupe de $G$ que l'on voit comme une repr\'esentation de $W_{F}\times SL(2,{\mathbb C})$ dans $GL(m_{G}^*,{\mathbb C})$ comme expliqu\'e dans l'introduction. L'orbite de Langlands associ\'ee \`a $\pi$ est celle qui est d\'efinie par la restriction de cette repr\'esentation \`a $SL(2,{\mathbb C})$. 

On va r\'ecrire cela en termes combinatoires de fa\c{c}on \`a pouvoir, si ce n'est la calculer, du moins l'approximer.

 La param\'etrisation combinatoire de Langlands, associe \`a $\pi$ un sous-groupe parabolique (standard) et une repr\'esentation temp\'er\'ee de ce sous-groupe parabolique tordu\'e n\'egativement, telle que $\pi$ soit un sous-module de l'induite. Etant donn\'ee la forme tr\`es particuli\`ere des sous-groupes paraboliques des groupes consid\'er\'es, les donn\'ees combinatoires associ\'ees \`a la param\'etrisation de Langlands de $\pi$ sont: une repr\'esentation temp\'er\'ee irr\'eductible d'un groupe de m\^eme type que $G$ mais de rang plus petit, not\'ee $\pi_{temp}$ et une collections de repr\'esentations de Steinberg g\'en\'eralis\'ees $St(\rho',a')$ tordues par un caract\`ere de la forme $\vert\,\vert^{-x'}$ avec $x'\in {\mathbb R}_{>0}$ tels que
$$
\pi\hookrightarrow \times_{(\rho,',a',x')\in \mathcal{L}(\pi)}St(\rho',a')\vert\,\vert^{-x'}\times \pi_{temp} (*)
$$
o\`u $\mathcal{L}(\pi)$ est totalement ordonn\'e de sorte que les $x'$ arrivent de fa\c{c}on d\'ecroissante au sens large; ; ici les $\rho'$ sont des repr\'esentations cuspidales unitaires, les $a'$ sont des entiers. Une telle \'ecriture n'est peut \^etre pas telle quelle dans la litt\'erature, voil\`a donc comment on la d\'emontre:  on fixe une telle inclusion avec aucune hypoth\`ese sur les $x'$ sauf que le premier est positif ce qui est \'evidemment possible si $\pi$ n'est pas \'egal \`a $\pi_{temp}$;  en utilisant simplement les propri\'et\'es des induites dans $GL(n)$ on s'arrange pour que les $x'$ arrivent dans l'ordre d\'ecroissant, pour cela on utilise pour $[d,f]$, $[d',f']$ des segments d\'ecroissant la propri\'et\'e que si le milieu de $[d,f]$ est sup\'erieur au milieu de $[d',f']$ les sous-quotients de l'induite $[d,f]_{\rho}\times [d',f']_{\rho}$ sont soit un sous-module de l'induite \'ecrite dans l'autre sens, soit les segments sont li\'es et l'autre sous-quotient est un sous-module de l'induite
$$
[d,f']_{\rho}\times [d',f]_{\rho} \simeq [d',f]_{\rho}\times [d,f']_{\rho}
$$
et les milieus de $[d',f]$ et de $[d,f]$ sont strictement inf\'erieur au milieu de $[d,f]$ (le segment $[d,f']$ n'intervient pas si $d'=f-1$); pour la suite il est utile de remarqu\'e que si le milieu de $[d,f]$ est n\'egatif, cela est aussi vrai pour $[d',f]$ et $[d,f]$.  On a donc am\'elior\'e la situation en gardant le fait que le milieu du premier segment n'a pas augment\'e et la somme des milieux n'a pas augment\'e. On obtient donc l'assertion de d\'ecroissance de proche en proche.

On montre que tous les $x'$ peuvent \^etre pris positifs: puisque l'on a d\'ej\`a ordonn\'e de fa\c{c}on d\'ecroissante, on consid\`ere $\tau:=\times_{(\rho',a',x'); x'\leq 0} \times \pi_{temp}$ et l'inclusion (1) se factorise par un sous-quotient irr\'eductible de $\tau$. On applique \`a $\tau$ l'assertion (on est n\'ecessairement dans un groupe de rang plus petit) et on r\'eordonne comme ci-dessus; on garde la positivit\'e des exposants et on obtient l'ordre d\'ecroissant. D'o\`u l'existence d'un inclusion (*) avec les propri\'et\'es voulues; l'unicit\'e est alors clair car cette inclusion est une inclusion de Langlands moyennant le fait qu'il faut regrouper entre elles les repr\'esentations de Steinberg avec la m\^eme torsion.

 Les donn\'ees sont donc uniques \`a l'ordre pr\`es, les modifications de l'ordre viennent des irr\'eductibilit\'es des induites pour le groupe lin\'eaire convenable de la forme $St(\rho',a')\vert\,\vert^{-x'}\times St(\rho'',a'')\vert\,\vert^{-x''}$ pour $x'=x''$. On associe \`a $\pi_{temp}$ son paquet temp\'er\'e, $\psi_{temp}$; c'est-\`a-dire que $\psi_{temp}$ est un morphisme tel que ceux consid\'er\'es ici pour un groupe de m\^eme type que $G$ tel que $\pi_{temp}\in \Pi(\psi_{temp})$.

Fixons $\rho$ une repr\'esentation cuspidale unitaire; on d\'efinit l'orbite unipotente $O_{\rho,\pi}$ comme une orbite d'un groupe lin\'eaire convenable ayant comme ensemble de blocs de Jordan deux copies de chaque $\cup_{(\rho',a',x')\in \mathcal{L}(\pi); \rho'=\rho}a'$  (compt\'e en plus  avec la multiplicit\'e \'eventuel dans $\mathcal{L}(\pi)$ auxquels on ajoute $\cup_{(\rho',a,1)\in Jord(\psi_{temp}); \rho\simeq \rho' }a$, l\`a aussi en tenant compte des multiplicit\'es \'eventuelles.

Soient $\psi$ un morphisme et $\rho$ comme ci-dessus. On note $O_{\psi,\rho}$ l'orbite unipotente du GL convenable dont les blocs de Jordan sont exactements
$\cup_{(\rho,a,b)\in Jord(\psi)}\underbrace{{a, \cdots, a}}_{b}$.

\begin{rem}Soit $\psi$ un morphisme et on suppose que $\pi\in \Pi(\psi_{L})$. Soit $\rho$ une repr\'esentation cuspidale unitaire; alors $O_{\pi,\rho}=O_{\psi,\rho}$.
\end{rem}
On a d\'ecrit les param\`etres de Langlands des \'el\'ements de $\Pi(\psi_{L})$ dans \ref{description}. La remarque s'en d\'eduit imm\'ediatement.
\subsection{Orbites de Langlands et induction\label{orbitedelanglandsetinduction}}

On fixe une repr\'esentation irr\'eductible $\pi$ de $G$ et une repr\'esentation induite de $G$ de la forme $\times_{i\in {\mathcal I}}St(\rho_{i},a_{i})\vert\,\vert^{-x_{i}}\times \pi_{temp}$ o\`u ${\mathcal I}$ est un ensemble d'indices totalement ordonn\'e et pour tout $i\in {\mathcal I}$, $\rho_{i}$ est une repr\'esentation cuspidale unitaire, $a_{i}$ est un entier et $x_{i}$ un r\'eel positif ou nul et o\`u $\pi'$ est une repr\'esentation irr\'eductible. On suppose que $$\pi\hookrightarrow \times_{i\in {\mathcal I}}St(\rho_{i},a_{i})\vert\,\vert^{-x_{i}}\times \pi'.\eqno(1)$$

\begin{lem} Avec les notations et hypoth\`eses pr\'ec\'edentes, pour toute repr\'esentation cuspidale $\rho$ unitaire irr\'eductible d'un groupe lin\'eaire, l'orbite unipotente ${\mathcal O}_{\pi,\rho}$ contient dans sa fermeture l'orbite unipotente du groupe lin\'eaire convenable dont les blocs de Jordan sont pr\'ecis\'ement l'union des blocs de Jordan de ${\mathcal O}_{\pi',\rho}$ avec $\{(a_{i},a_{i}), i\in {\mathcal I}; \rho_{i}\simeq \rho\}$. 
\end{lem}
On remarque que pour prouver le lemme, on peut tout \`a fait supposer que $\pi'$ est temp\'er\'ee; si ceci n'est pas vrai on remplace $\pi'$ par l'induite qui donne sa param\'etrisation de Langlands. On suppose donc que $\pi'$ est temp\'er\'ee.

On fixe $\rho$ comme dans l'\'enonc\'e et on note ${\mathcal O}_{ind,\rho}$ l'orbite unipotente dont les blocs de Jordan sont l'union de ceux de ${\mathcal O}_{\pi',\rho}$ avec les $(a_{i},a_{i})$ pour tout $i\in {\mathcal I}$ tel que $\rho_{i}\simeq \rho$. On a ${\mathcal O}_{\pi,\rho}={\mathcal O}_{ind,\rho}$ si pour tout $i,j\in {\mathcal I}$ tel que $i<j$ et $\rho_{i}\simeq \rho_{j}\simeq \rho$, $x_{i}>x_{j}$. Si ces in\'egalit\'es ne sont pas v\'erifi\'ees,  il faut \'echanger des facteurs, \'eventuellement  en les modifiant pour passer de l'inclusion donn\'ee en une inclusion comme sous-module de Langlands.

Ceci se fait par \'etape \'el\'ementaire et on est ramen\'e \`a regarder ce qui se passe quand on a deux facteurs cons\'ecutifs $St(\rho,a_{i})\vert\,\vert^{-x_{i}}\times St(\rho,a_{i+1})\vert\,\vert^{-x_{i+1}}$ avec $x_{i}<x_{i+1}$. Deux cas sont alors possibles, soit on peut simplement \'echanger les 2 facteurs de l'induite en gardant l'inclusion de $\pi$ dans cette nouvelle induite soit il faut remplacer $St(\rho,a_{i})\vert\,\vert^{-x_{i}}\times St(\rho,a_{i+1})\vert\,\vert^{-x_{i+1}}$ par son sous-module irr\'eductible. Le premier cas se produit  sauf \'eventuellement si les segments $[(a_{i}-1)/2-x_{i}, -(a_{i}-1)/2-x_{i}]$ et $[(a_{i+1}-1)/2-x_{i+1}, -(a_{i+1}-1)/2-x_{i+1}]$ sont li\'es au sens de Zelevinski. 

Explicitons ce cas: le premier segment est un segment de centre $-x_{i}$ et le deuxi\`eme est de centre $-x_{i+1}$; l'hypoth\`ese $x_{i+1}>x_{i}$, entra\^{\i}ne que les segments sont li\'es si $$(a_{i}-1)/2-x_{i}> (a_{i+1}-1)/2-x_{i+1}\geq -(a_{i}-1)/2-x_{i}-1>-(a_{i+1}-1)/2-x_{i+1}-1.$$
Le sous-module irr\'eductible est alors l'induite irr\'eductible $St(\rho,\alpha)\vert\,\vert^{-x'}\times St(\rho,\beta)\vert\,\vert^{-x''}$ avec $$(\alpha-1)/2-x'=(a_{i}-1)/2-x_{i}; -(\alpha-1)/2-x'=-(a_{i+1}-1)/2-x_{i+1}$$
$$(\beta-1)/2-x''=(a_{i+1}-1)/2-x_{i+1}; -(\beta-1)/2-x'=-(a_{i}-1)/2-x_{i}$$avec $\beta=0$ si $(a_{i+1}-1)/2-x_{i+1}= -(a_{i}-1)/2-x_{i}-1$. Dans tous les cas $\alpha \geq sup (a_{i},a_{i+1})$ et $\beta+\alpha=a_{i}+a_{i+1}$. De plus $x'$ et $x''$ sont des r\'eels; montrons qu'ils sont strictement compris entre $x_{i}$ et $x_{i+1}$; on peut le voir en interpr\'etant $-x'$ et $-x''$ comme des milieux de segments et on voit alors tout de suite $x_{i}\leq x',x'' \leq x_{i+1}$. Par le calcul,  on obtient:
$$
-2x'=(a_{i}-1)/2-x_{i}-(a_{i+1}-1)/2-x_{i+1}=-2x_{i}+(a_{i}-1)/2+x_{i}-(a_{i+1}-1)/2-x_{i+1}<-2x_{i}$$car $-(a_{i+1}-1)/2-x_{i+1}<-(a_{i}-1)/2-x_{i}$, (cf. ci-dessus). De m\^eme
$$
-2x''=(a_{i+1}-1)/2-x_{i+1}-(a_{i}-1)/2-x_{i}=-2x_{i}+((a_{i+1}-1)/2-x_{i+1})- ((a_{i}-1)/2-x_{i})<-2x_{i}.
$$On remarque que $-2(x'+x'')=-2(x_{i}+x_{i+1})$, d'o\`u, ce qui nous servira plus loin:
$$
x'+x''=x_{i}+x_{i+1} \hbox { et } sup(x',x'')\leq x_{i+1},\eqno(1)
$$
puisque l'on a montr\'e que $x_{i}\leq sup(x',x'')$.

Finalement on obtient une inclusion de $\pi$ dans une induite qui est meilleure au sens que l'on se rapproche d'une inclusion de Langlands: pour cela il faut donner un nombre sur lequel faire une r\'ecurrence. On fait en fait une double r\'ecurrence, d'abord sur $\vert {\mathcal I}\vert$ puis sur

$inv(ind):=
\vert \{(i,j); i<j, \rho_{i}\simeq \rho_{j}, x_{i}<x_{j}$ et les segments associ\'es \`a $St(\rho_{i},a_{i})\vert\,\vert^{-x_{i}}$ et $St(\rho_{j},a_{j})\vert\,\vert^{-x_{j}}$ sont li\'es $\}\vert$. 

L'initialisation de la r\'ecurrence est \'evidente.

Il est aussi clair que si $\pi$ est un sous-module de  l'analogue de l'induite (1) mais o\`u on a \'echang\'e 2 facteurs cons\'ecutifs comme ci-dessus avec $x_{i}<x_{i+1}$, on a am\'elior\'e la r\'ecurrence. On suppose donc que $\pi$ est un sous-module de l'analogue de l'induite (1) mais o\`u on a remplac\'e 2 facteurs cons\'ecutifs associ\'es \`a des segments li\'es par leur sous-module irr\'eductible comme d\'ecrit ci-dessus; on reprend les notations d\'ej\`a introduite et on montre que l'on am\'eliore aussi la r\'ecurrence.
On a a consid\'erer uniquement le cas o\`u $\beta\neq 0$ dans les notations ci-dessus. Supposons cela; on transforme les notations d\'ej\`a introduites pour mieux voir les liaisons de segment. Pour tout $j\in \mathcal{I}$, on note $[d_{j},f_{j}]$ le segment associ\'e \`a la repr\'esentation $St(\rho_{j},a_{j})\vert\,\vert^{-x_{j}}$ pour tout $j\in \mathcal{I}$, c'est-\`a-dire, $d_{j}=(\alpha_{j}-1)/2-x_{j}$ et $f_{j}=-(\alpha_{j}-1)/2-x_{j}$. Avec ces notations pour $j=i$ ou $i+1$, on a (puisque $\beta\neq 0$)
$$
d_{i}>d_{i+1}\geq f_{i}>f_{i+1}$$. Et les segments associ\'es aux nouvelles repr\'esentations sont $[d_{i},f_{i+1}]$ et $[d_{i+1},f_{i}]$.  On fait un dessin en mettant en consid\'erant que la d\'ecroissance se fait de gauche \`a droite. On a 
$$
\begin{matrix}
d_{i} & &&f_{i}\\
&d_{i+1}&&&f_{i+1}
\end{matrix} \qquad
\begin{matrix}
&d_{i}&&&&f_{i+1}\\
&&d_{i+1}&&f_{i}
\end{matrix}
$$En particulier le couple $(i,i+1)$ qui contribuait \`a $inv(ind)$ n'y contribue plus. Soit $j\in \mathcal {I}$ et supposons d'abord que $j<i$.
Si $d_{j}\leq d_{i+1}$, il n'y a pas de liaison telles qu'on les cherche avec les nouveaux segments. Si $d_{j}\in [d_{i},d_{i+1}[$, on a une liaison seulement si $f_{j}\in [d_{i+1}+1,f_{i}[$ et uniquement avec le segment $[d_{i+1},f_{i}]$; dans ce cas, on a aussi une liaison avec le segment $[d_{i+1},f_{i+1}]$. On suppose maintenant que $d_{j}<d_{i}$; si $f_{j}\leq f_{j+1}$, il n'y a aucune liaison. Si $f_{j}\in [f_{i},f_{i+1}[$ on a une liaison avec le segment $[d_{i},f_{i+1}]$ et une avec le segment $[d_{i+1},f_{i+1}]$. Si $f_{j }\in [d_{i+1}+1,f_{i}[$ on a une liaisons avec les 4 segments \'ecrits et si $f_{j}\in [d_{i}+1,d_{i+1}+1[$, on a une liaison avec les segments $[d_{i},f_{i+1}]$ et $[d_{i},f_{i}]$. En d'autres termes on n'a pas augmenter la contribution de $[d_{j},f_{j}]$ \`a $inv(ind)$.

Il faut aussi regarder le cas o\`u $j>i+1$ qui est analogue. Ainsi on a am\'elior\'e la r\'ecurrence.

Il reste \`a suivre les orbites unipotentes associ\'es \`a la repr\'esentation induite.
 Dans le cas d'\'echange de 2 facteurs cons\'ecutifs, on ne change pas l'orbite unipotente associ\'ee \`a l'induite. Dans l'autre cas,  on remplace $(a_{i},a_{i},a_{i+1},a_{i+1})$ par $(\alpha,\alpha,\beta,\beta)$. La  nouvelle orbite contient la premi\`ere dans sa fermeture car $\alpha > sup(a_{i},a_{i+1})$ et $\alpha+\beta=a_{i}+a_{i+1}$. L'hypoth\`ese de r\'ecurrence donne alors le lemme.
 
\subsection{Comparaison des orbites de Langlands \`a l'int\'erieur d'un paquet d'Arthur\label{comparaisonorbite}}
 \begin{thm}On fixe un morphisme $\psi$ et $\pi\in \Pi(\psi)$, une repr\'esentation irr\'eductible dans le paquet associ\'e \`a $\psi$. Soit $\rho$ une repr\'esentation cuspidale unitaire; l'orbite $O_{\psi,\rho}$ est une orbite du m\^eme groupe lin\'eaire que l'orbite $O_{\psi,\rho}$ (c'est-\`a-dire $\sum_{(\rho',a',b'); \rho\simeq \rho}a'b'$). Et l'orbite $O_{\psi,\rho}$ est incluse dans la fermeture de l'orbite $O_{\pi,\rho}$.
 \end{thm}
On d\'emontre le th\'eor\`eme par r\'ecurrence sur $\ell(\psi):=\sum_{(\rho,a,b)\in Jord(\psi)}(b-1)$. Si $\ell(\psi)=0$, $\psi$ n'est pas n\'ecessairement temp\'er\'ee \`a cause de $(\rho',a',b')\in Jord(\psi)$ avec $\rho\not\simeq \rho$ mais les constructions ram\`enent ais\'ement au cas temp\'er\'e o\`u il n'y a rien \`a d\'emontrer.

On suppose donc que $\ell(\psi)>0$. On traite d'abord le cas o\`u il existe $(\rho,a,b)\in Jord(\psi)$ avec $a\geq b>1$. On fixe un tel triplet en supposant $(a-b)/2$ maximum avec ces propri\'et\'es. On utilise \ref{blocpositif} pour \'ecrire $\pi\hookrightarrow \delta\times \pi''$ avec $\delta$ et $\pi''$ pr\'ecis\'es en loc.cite. Le lemme de \ref{orbitedelanglandsetinduction} montre que $O_{\pi,\rho}$ contient dans sa fermeture  l'orbite associ\'ee en loc. cite \`a $\delta\times \pi''$. L'hypoth\`ese de r\'ecurrence assure que l'orbite associ\'e \`a $\rho$ et $\delta\times \pi''$ contient dans sa fermeture l'orbite associ\'ee \`a $\rho$ et $\delta\times \psi''$ (pour \^etre plus correct, disons l'orbite associ\'e \`a $\delta\times\pi'''$ avec $\pi''' \in \Pi(\psi''_{L})$). Le fait que cette derni\`ere orbite  contient dans sa fermeture $O_{\psi,\rho}$ a \'et\'e v\'erifi\'e dans les 2 premiers lemmes de \ref{csqsurorbite}. 

On consid\`ere maintenant le cas o\`u pour tout $(\rho,a,b)\in Jord(\psi)$ si $a\geq b$ alors $b=1$. Comme $\ell(\psi)>0$, il existe $(\rho,a,b)$ avec $a<b$, ce qui force $b>1$. On fixe un tel \'el\'ement avec encore $(a-b)/2$ maximum avec cette propri\'et\'e de n\'egativit\'e. On raisonne comme ci-dessus, en utilisant \ref{blocnegatif} au lieu de \ref{blocpositif} et le dernier lemme de \ref{csqsurorbite}. Cela termine la preuve.
\begin{rem}Le m\^eme r\'esultat est vrai sans supposer en rempla\c{c}ant $\pi_{temp}$ par une repr\'esentation irr\'eductible $\pi''$ irr\'eductible mais non n\'ecessairement temp\'er\'ee.
\end{rem}
En effet on commence par plonger $\pi''$ dans un induite dont les exposants sont dans la chambre de Weyl n\'egative (sous-module de Langlands) et on applique le lemme d\'emontr\'e en augmentant $\mathcal{I}$.

\subsection{Paquets d'Arthur contenant une repr\'esentation non ramifi\'ee\label{nonramifie}}
Il est tout \`a fait possible que la proposition ci-dessous ait en fait une preuve plus simple que celle donn\'ee ici.
\begin{prop} Soit $\psi$ un morphisme alors $\Pi(\psi)$ contient une repr\'esentation non ramifi\'ee si et seulement si la restriction de $\psi$ \`a $W_{F}$ est non ramifi\'ee et la restriction de $\psi$ \`a la premi\`ere copie de $SL(2,{\mathbb C})$ est triviale. La repr\'esentation non ramifi\'ee dans $\Pi(\psi)$ est alors unique.
\end{prop}
On fixe $\psi$ et $\pi\in \Pi(\psi)$ et on suppose que $\pi$ est non ramifi\'ee. L'orbite de Langlands associ\'ee \`a $\pi$ est l'orbite triviale; il r\'esulte donc de \ref{comparaisonorbite} que la restriciton de $\psi$ \`a la premi\`ere copie de $SL(2,{\mathbb C})$ est triviale. Le support cuspidal de $\pi$ est le support cuspidal \'etendu de $\pi$ puisque $\pi$ est non ramifi\'ee (cf. la d\'efinition de \ref{definitionsupportcuspidal}) et il se calcule \`a l'aide de $\psi$; cela force le fait que la restriction de $\psi$ \`a $W_{F}$ doit \^etre non ramifi\'e. La fin de la proposition a alors \'et\'e prouv\'ee en \ref{unicite}.

\section{Exposants}
\subsection{D\'efinitions des exposants}
Dans ce paragraphe on donne une majoration des exposants de toutes les repr\'esentations dans un paquet $\Pi(\psi)$. Supposons d'abord que $\pi\in \Pi(\psi_{L})$; on garde la notation $\delta_{b}=1/2$ si $b$ est un entier pair et $1$ si $b>1$ est impair; on a \'evidemment $\delta_{b}=(b-1)/2-[b/2]+1$. Les exposants de $\pi$ sont la collection des demi-entiers $\cup_{(\rho,a,b)\in Jord(\psi); b>1}\cup_{c\in [(b-1)/2,\delta_{b}]}c$. On note $Exp({\psi})$ l'ensemble de ces  demi-entiers strictement positifs, ensemble avec multiplicit\'e.

Soit $\pi$ une repr\'esentation irr\'eductible de $G$; on note $Exp(\pi)$ un ensemble de demi-entier strictement positifs, avec multiplicit\'e, ordonn\'e par l'ordre d\'ecroissant tel que pour tout $x\in Exp(\pi)$, il existe une repr\'esentation de Steinberg $St(\rho_{x},d_{x})$ et il existe une repr\'esentation temp\'er\'ee $\pi_{temp}$ tel que l'on ait une inclusion
$$
\pi\hookrightarrow \times_{x\in Exp(\pi)} St(\rho_{x},d_{x})\vert\,\vert^{-x}\times \pi_{temp}.\eqno(1)
$$En particulier $Exp(\pi)$ est bien uniquement d\'etermin\'e.

On veut d\'efinir $Exp(\pi)\leq Exp(\psi)$; cela se fait comme pour les orbites unipotentes, on dit que $Exp(\pi)\leq Exp(\psi)$ si pour tout nombre entier $t$, la somme des $t$ plus grand \'el\'ements de $Exp(\pi)$ (on ajoute des 0 si n\'ecessaire) est sup\'erieur ou \'egal \`a la somme des $t$ plus grand \'el\'ements de $Exp(\psi)$ (l\`a aussi, on ajoute des 0 si n\'ecessaire).

On peut aussi fixer $\rho$ comme on l'a fait dans tout ce qui pr\'ec\`ede et ne regarder que les exposants tels que le $\rho_{x}$ ci-dessus soit $\rho$, on note alors $Exp_{\rho}(\pi)$ et $Exp_{\rho}(\psi)$. Il est clair que $Exp_{\rho}(\psi)=Exp_{\rho,\pi}$ si $\pi\in \Pi(\psi_{L})$, c'est-\`a-dire est dans le paquet de Langlands \`a l'int\'erieur du paquet d'Arthur.
\subsection{Comparaison des exposants \`a l'int\'erieur d'un paquet d'Arthur\label{comparaisonexposants}}
\begin{prop}Soit $\psi$ un morphisme et soit $\pi\in \Pi(\psi)$. Alors $Exp(\pi)\leq Exp(\psi)$ et plus pr\'ecis\'ement pout tout $\rho$ comme ci-dessus, $Exp_{\rho}(\pi)\leq Exp_{\rho}(\psi)$
\end{prop}
On fait comme dans \ref{comparaisonorbite}  une r\'ecurrence sur $\ell(\psi):=\sum_{(\rho,a,b)\in Jord(\psi)}(b-1)$. Si $\ell(\psi)=0$, le morphisme est temp\'er\'e et  $Exp_{\rho}(\pi)=0=Exp_{\rho}(\psi)$, d'o\`u \`a fortiori la proposition. La d\'emonstration suit les m\^emes lignes que celle de \ref{comparaisonorbite} et il faut donc reprendre du point de vue des exposants plusieurs lemmes. D'abord il faut l'analogue du lemme \ref{orbitedelanglandsetinduction}: soit $(\rho,a_{i},x_{i})$ une collection index\'ee par l'ensemble $\mathcal{I}$ de triplets, o\`u $\rho$ est fix\'e, $a_{i}$ est un entier positif et $x_{i}$ un r\'eel strictement positif; on fixe aussi une repr\'esentation irr\'eductible $\pi''$. On suppose que $\pi\hookrightarrow \times_{i\in \mathcal{I}}St(\rho,a_{i})\vert\,\vert^{-x_{i}}\times \pi''$. Alors $Exp_{\rho}(\pi)\leq (Exp_{\rho}(\pi'')\cup_{i\in \mathcal{I}}x_{i})$; en effet, quitte \`a remplacer $\pi''$ par l'inclusion de Langlands dans une induite ayant ses exposants dans la chambre de Weyl n\'egative, on peut supposer que $\pi''$ est temp\'er\'ee. L'assertion se d\'emontre pas \`a pas en transformant l'inclusion de $\pi$ donn\'ee dans l'\'enonc\'e en une inclusion de Langlands comme cela a \'et\'e expliqu\'ee en \ref{orbitedelanglandsetinduction} et l'assertion cherch\'ee est le (1) de loc.cite.

Il faut aussi l'analogue des lemmes de \ref{csqsurorbite}; le point est dans tous les cas de v\'erifier que $(b-1)/2$ est sup\'erieur ou \'egal aux exposants intervenant dans les induites. Pour les 2 premiers lemmes, les induites sont de la forme:
$$
\times_{j\in [1,v/2} <(a_{2j}-1)/2,-(a_{2j-1}-1)/2>_{\rho} ,$$
  
  ou $\times_{j\in [1,([v-1)/2]]} <(a_{2j+1}-1)/2,-(a_{2j}-1)/2>_{\rho} \times <(a_{1}-1)/2,-(a+b)/1+1>_{\rho}
$

ou bien encore
$
<(a-b)/2,-(a_{v}-1)/2>\times_{j\in [1,v[, j\equiv v[2]} <(a_{j+1}-1)/2, -(a_{j}-1)/2>_{\rho}
$
auquel on ajoute encore le facteur $<(a_{1}-1)/2,-(a+b)/2+1>$ si $v$ est pair.

\noindent
Il suffit de d\'emontrer que la somme des exposants est inf\'erieure ou \'egale \`a $(b-1)/2$; or pour trouver la somme des exposants, il faut prendre l'oppos\'e de la demi-somme des extr\^emit\'es. Par exemple pour la derni\`ere induite \'ecrite, on trouve 
$
1/2 \sum_{j\in [1,v/2]}(a_{2j-1}-a_{2j})/2
$ plus, si $v$ est pair $1/2((a+b)/2-1-(a_{1}-1)/2$.
On utilise le fait que les $a_{j}$ sont rang\'es dans l'ordre d\'ecroissant et que  de plus $a_{1}<a+b-1$ et $a_{v}> a-b+1$. Donc si $v$ est impair la somme des exposants est certainement inf\'erieure \`a $1/2 ((a_{1}-1)/2-(a-b)/2)< 1/2((a+b)/2-1-(a-b)/2)=1/2(b-1)$. Si $v$ est pair, on obtient directement l'in\'egalit\'e avec $1/2((a+b)/2-1-(a-b)/2)=1/2(b-1)$. Pour les premi\`eres induites le calcul est \'evidemment analogue.

Pour le dernier lemme de \ref{csqsurorbite}, l'induite \`a consid\'er\'ee est de la forme 
$$
<(a-b)/2,-(a_{v}-1)/2>_{\rho}\times \times_{j\in [1,v[, j\equiv v[2]} <(a_{j+1}-1)/2, -(a_{j}-1)/2>_{\rho}
$$
auquel on ajoute encore le facteur $<(a_{1}-1)/2,-(a+b)/2+1>$ si $v$ est pair. Les $a_{i}$ sont encore rang\'es dans l'ordre d\'ecroissant mais ici $(a-b)/2<0$ et on a $a_{v}\geq b-a-1$ et $a_{1}< a+b-1$. Que $v$ soit pair ou impair, on obtient que la somme des exposants est inf\'erieur ou \'egal \`a
$$
1/2((a+b)/2-1-(a-b)/2)=1/2(b-1).
$$
La d\'emonstration est alors celle de \ref{comparaisonorbite}.

\end{document}